\documentclass{ifacmtg}
\usepackage{natbib}            
\usepackage{graphicx}          
\usepackage[latin1]{inputenc}
\usepackage[T1]{fontenc}
\usepackage[english]{babel}
\usepackage{graphicx}
\usepackage{psfrag}

\usepackage{epsfig}
\usepackage[hang]{caption2}
\usepackage[normalsize,it,hang]{subfigure}
\usepackage{amsmath}
\usepackage{array}
\usepackage{theorem}
\usepackage{amsmath,amssymb}
\usepackage{amsfonts}
\usepackage{amssymb}
\usepackage{subfigure}
\usepackage{amssymb,amsbsy,amsmath,amsfonts,amssymb,amscd}

\newtheorem{remarque}{\it Remark\/}

\usepackage[dvipsnames,table,svgnames]{xcolor}

\usepackage[dvipsnames,table,svgnames]{xcolor}

\newcommand{\Vect}[1]{ \boldsymbol{#1}}

\begin{document}
\runauthor{Fliess and Join}
\begin{frontmatter}
\title{Model-free control and intelligent PID controllers:
towards a possible trivialization of nonlinear control?}
\author[Baiae,Paestum]{Michel FLIESS}
\author[Baiae,Rome]{\quad C\'{e}dric JOIN}

\address[Baiae]{INRIA-ALIEN}
\address[Paestum]{LIX (CNRS, UMR 7161),
\'Ecole polytechnique \\ 91128 Palaiseau, France \\ {\tt
Michel.Fliess@polytechnique.edu}}
\address[Rome]{CRAN (CNRS, UMR 7039), Nancy-Universit\'{e} \\
BP 239, 54506 Vand\oe{}uvre-l\`es-Nancy, France
\\ {\tt Cedric.Join@cran.uhp-nancy.fr} }

\begin{abstract}.
~We are introducing a model-free control and a control with a
restricted model for finite-dimensional complex systems. This
control design may be viewed as a contribution to ``intelligent''
PID controllers, the tuning of which becomes quite straightforward,
even with highly nonlinear and/or time-varying systems. Our main
tool is a newly developed numerical differentiation. Differential
algebra provides the theoretical framework. Our approach is
validated by several numerical experiments.
\end{abstract}

\end{frontmatter}
\section{Introduction}
Writing down simple and reliable differential equations for
describing a concrete plant is almost always a daunting task. How to
take into account, for instance, frictions, heat effects, ageing
processes, characteristics dispersions due to mass production,
\dots? Those severe difficulties explain to a large extent why the
industrial world is not willing to employ most techniques stemming
from ``modern'' control theory, which are too often based on a
``precise'' mathematical modeling, in spite of considerable advances
during the last fifty years. We try here\footnote{This communication
is a slightly modified and updated version of (\cite{esta}), which
is written in French. Model-free control and control with a
restricted model, which might be useful for hybrid systems
(\cite{bourdais}), have already been applied in several concrete
case-studies in various domains
(\cite{choi,brest,poitiers,vil1,vil}). Other applications on an
industrial level are being developed.} to overcome this unfortunate
situation thanks to recent fast estimation methods.\footnote{See,
{\it e.g.}, (\cite{nl,mboup}) and the references therein.} Two cases
are examined:
\begin{enumerate}
\item {\em Model-free control} is based on an elementary continuously updated {\em local}
modeling via the unique knowledge of the input-output behavior. It
should not be confused with the usual ``black box'' identification
(see, {\it e.g.}, (\cite{aut2,aut1})), where one is looking for a
model which is valid within an operating range which should be as
large as possible.\footnote{This is why we use the terminology
``model-free'' and not ``black box''.} Let us summarize our approach
in the monovariable case. The input-output behavior of the system is
assumed to ``approximatively'' governed within its operating range
by an unknown finite-dimensional ordinary differential equation,
which is not necessarily linear,
\begin{equation}\label{E}
E (y, \dot{y}, \dots, y^{(a)}, u, \dot{u}, \dots, u^{(b)}) = 0
\end{equation}
We replace Eq. \eqref{E} by the following ``pheno\-menological''
model, which is only valid during a very short time interval,
\begin{equation}\label{F}
y^{(\nu)} = F + \alpha u
\end{equation}
The derivation order $\nu$, which is in general equal to $1$ or $2$,
and the constant parameter $\alpha$ are chosen by the practitioner.
It implies that $\nu$ is not necessarily equal to the derivation
order $a$ of $y$ in Eq. (\ref{E}). The numerical value of $F$ at any
time instant is deduced from those of $u$ and $y^{(\nu)}$, thanks to
our numerical differentiators (\cite{nl,mboup}). The desired
behavior is obtained by implementing, if, for instance, $\nu = 2$,
the {\em intelligent PID controller}\footnote{This terminology, but
with other meanings, is not new in the literature (see, {\it e.g.},
\cite{intelligent}).} ({\em i-PID})
\begin{equation}\label{universal}
u = - \frac{F}{\alpha} + \frac{\ddot{y}^\ast}{\alpha} + K_P e + K_I
\int e + K_D \dot{e}
\end{equation}
where
\begin{itemize}
\item $y^\ast$ is the output reference trajectory, which
is determined via the rules of flatness-based control (see, {\it
e.g.}, \cite{flmr,rotella,hsr});
\item $e = y - y^\ast$ is the tracking error;
\item $K_P$, $K_I$, $K_D$ are the usual tuning gains.
\end{itemize}

\item Assume now that a restricted or partial model of the plant
is quite well known and is defined by Eq. (\ref{E}) for instance.
The plant is then governed by the  {\em restricted}, or {\em
incomplete}, {\em modeling},\footnote{Since we are not employing the
terminology ``black box'', we are also here not employing the
terminology ``grey box''.}
$$
E (y, \dot{y}, \dots, y^{(a)}, u, \dot{u}, \dots, u^{(b)}) + G = 0
$$
where $G$ stands for all the unknown parts.\footnote{Any
mathematical modeling in physics as well as in engineering is
incomplete. Here the term $G$ does not need to be ``small'' as in
the classic approaches.} If the known system, which corresponds to
$E = 0$, is flat, we also easily derive an {\em intelligent
controller}, or {\em i-controller}, which gets rid of the unknown
effects.

\end{enumerate}

Differential algebra is briefly reviewed in Sect. \ref{algebra} in
order to
\begin{itemize}
\item derive the input-output differential equations,
\item define minimum and non-minimum phase systems.
\end{itemize}
Numerical differentiation of noisy signals is examined in Sect.
\ref{estder}. Sect. \ref{cs} states the basic principles of our
model-free control. Several numerical experiments\footnote{Those
computer simulations would of course be impossible without precise
mathematical models, which are {\it a priori} known.} are reported
in Sect. \ref{exem}. We deal as well with linear\footnote{The remark
\ref{robuste} underlines the following crucial fact: the usual
mathematical criteria of robust control are becoming pointless in
this new setting.} and nonlinear systems and with monovariable and
multivariable systems. An anti-windup strategy is sketched in Sect.
\ref{antiemb}. Sect. \ref{restricted} studies the control with a
partially known medeling. One of the two examples deals with a
non-minimum phase system.\footnote{Non-minimum phase systems are
today beyond our reach in the case of model-free control. This is
certainly the most important theoretical question which is left open
here.} The numerical simulations of Sect. \ref{list} and Sect.
\ref{spring} show the that our controllers behave much better than
classic PIDs.\footnote{We are perfectly aware that such a comparison
might be objected. One could always argue that an existing control
synthesis in the huge literature devoted to PIDs since
\cite{ziegler} (see, {\it e.g.},
\cite{pid,book,voda,dat,dindeleux,franklin,john,leq,od,zamb,shin,vis,wang,yu})
has been ignored or poorly understood. Only time and the work of
many practitioners will be able to confirm our viewpoint.} Several
concluding remarks are discussed in Sect. \ref{conclusion}.

\begin{remarque}
Only Sect. \ref{algebra} is written in an abstract algebraic
language. In order to understand the sequel and, in particular, the
basic principles of our control strategy, it is only required to
admit the input-output representations \eqref{E} and \eqref{io}, as
well as the foundations of flatness-based control
(\cite{flmr,rotella,hsr}).
\end{remarque}

\section{Nonlinear systems}\label{algebra}
\subsection{Differential fields}
All the fields considered here are commutative and have
characteristic $0$. A {\em differential field}\footnote{See, {\it
e.g.}, (\cite{cl,kolchin}) for more details and, in particular,
(\cite{cl}) for basic properties of usual fields, {\it i.e.},
non-differential ones.} $\mathfrak{K}$ is a field which is equipped
with a {\em derivation} $\frac{d}{dt}$, {\it i.e.}, a mapping
$\mathfrak{K} \rightarrow \mathfrak{K}$ such that, $\forall ~ a, b
\in \mathfrak{K}$,
\begin{itemize}
\item $\frac{d}{dt}(a + b) = \dot{a} + \dot{b}$,
\item $\frac{d}{dt}(ab) = \dot{a}b + a\dot{b}$.
\end{itemize}
\noindent A {\em constant} $c \in \mathfrak{K}$ is an element such
that $\dot{c} = 0$. The set of all constant elements is the {\em
subfield of constants}.

A differential field {\em extension} $\mathfrak{L} / \mathfrak{K}$
is defined by two differential fields $\mathfrak{K}$, $\mathfrak{L}$
such that:
\begin{itemize}
\item ${\mathfrak K} \subseteq {\mathfrak L}$,
\item the derivation of $\mathfrak{K}$ is the restriction
to $\mathfrak{K}$ of the derivation of $\mathfrak{L}$.
\end{itemize}
Write ${\mathfrak K}\langle S \rangle$, $S \subset {\mathfrak L}$,
the differential subfield of $\mathfrak{L}$ generated by
$\mathfrak{K}$ and $S$. Assume that $\mathfrak{L} / \mathfrak{K}$ is
finitely generated, {\it i.e.}, ${\mathfrak L} = {\mathfrak
K}\langle S \rangle$, where $S$ est finite. An element $\xi \in
\mathfrak{L}$ is said to be {\em differentially algebraic} over
$\mathfrak{K}$ if, and only if, it satisfies an algebraic
differential equation $P(\xi, \dots, \xi^{(n)}) = 0$, where $P$ is a
polynomial function over $\mathfrak{K}$ in $n + 1$ variables. The
extension $\mathfrak{L} / \mathfrak{K}$ is said to be {\em
differentially algebraic} if, and only if, any element of
$\mathfrak{L}$ is differentially algebraic over $\mathfrak{K}$. The
next result is important: \\ $\mathfrak{L} / \mathfrak{K}$ is
differentially algebraic if, and only if, its transcendence degree
is finite.

An element of $\mathfrak{L}$, which is non-differentially algebraic
over $\mathfrak{K}$, is said to {\em differentially transcendental}
over $\mathfrak{K}$. An extension $\mathfrak{L} / \mathfrak{K}$,
which is non-differentially algebraic, is said to be {\em
differentially transcendental}. A set $\{ \xi_\iota \in \mathfrak{L}
\mid \iota \in I \}$ is said to be {\em differentially algebraically
independent} over $\mathfrak{K}$ if, and only if, there does not
exists any non-trivial differential relation over $\mathfrak{K}$: \\
$Q(\dots, \xi_\iota^{(\nu_\iota)}, \dots ) = 0$, where $Q$ is a
polynomial function over $\mathfrak{K}$, implies $Q \equiv 0$. \\
Two such sets, which are maximal with respect to set inclusion, have
the same cardinality, {\it i.e.}, they have the same number of
elements: this is the {\em differential transcendence degree} of the
extension $\mathfrak{L} / \mathfrak{K}$. Such a set is a {\em
differential transcendence basis}. It should be clear that
${\mathfrak L} / {\mathfrak K}$ is differentially algebraic if, and
only if, its differential transcendence degree is $0$.

\subsection{Nonlinear systems}
\subsubsection{General definitions}
Let $k$ be a given differential ground field. A {\em
system}\footnote{See also (\cite{delaleau,delaleau2,nl,flmr}) which
provide more references on the use of differential algebra in
control theory.} is a finitely generated differentially
transcendental extension ${\mathcal{K}}/k$. Let $m$ be its
differential transcendence degree. A set of {\em (independent)
control variables} ${\Vect{u}} = (u_1, \dots, u_m)$ is a
differential transcendence basis of ${\mathcal{K}}/k$. The extension
${\mathcal{K}} / k \langle \Vect{u} \rangle$ is therefore
differentially algebraic. A set of {\em output variables}
${\Vect{y}} = (y_1, \dots, y_p)$ is a subset of ${\mathcal{K}}$.

Let $n$ be the transcendence degree of ${\mathcal{K}} / k \langle
\Vect{u} \rangle$ and let ${\Vect{x}} = (x_1, \dots, x_n)$ be a
transcendence basis of this extension. It yields the generalized
state representation:
$$
\begin{array}{l}
A_\iota(\dot{x}_\iota, {\Vect{x}}, {\Vect{u}}, \dots, {\Vect{u}}^{(\alpha)}) = 0 \\
B_\kappa(y_\kappa, {\Vect{x}}, {\Vect{u}}, \dots,
{\Vect{u}}^{(\beta)}) = 0
\end{array}
$$

\noindent where $A_\iota$, $\iota = 1, \dots, n$, $B_\kappa$,
$\kappa = 1, \dots, p$, are polynomial functions over $k$.

The following input-output representation is a consequence from the
fact that $y_1, \dots, y_p$ are differentially algebraic over $k
\langle {\mathbf u} \rangle$: {
\begin{equation} \label{io}
\Phi_j({\Vect{y}}, \dots, {\Vect{y}}^{(\bar{N}_j)}, {\Vect{u}},
\dots, {\Vect{u}}^{(\bar{M}_j)} ) = 0
\end{equation}
where $\Phi_j$, $j = 1, \dots, p$, is a polynomial function over
$k$.

\subsubsection{Input-output invertibility}
\begin{itemize}\item The system is said to be {\em left invertible}
if, and only if, the extension\footnote{This extension $k \langle {
{\Vect{u}}, \Vect{y}} \rangle / k \langle { \Vect{y}} \rangle $ is
called the {\em residual dynamics}, or the {\em zero dynamics}
(compare with \cite{isidori2}).} $k \langle { {\Vect{u}}, \Vect{y}}
\rangle / k \langle { \Vect{y}} \rangle $ is differentially
algebraic. It means that one can compute the input variables from
the output variables via differential equations. Then $m \leq p$.
\item It is said to be {\em right invertible} if, and only if, the
differential transcendence degree of $k \langle {\Vect{y}} \rangle /
k $ is equal to $p$. This is equivalent saying that the output
variables are differentially algebraically independent over $k$.
Then $p \leq m$.
\end{itemize}
The system is said to be {\em square} if, and only if, $m = p$. Then
left and right invertibilities coincide. If those properties hold
true, the system is said to be {\em invertible}.
\subsubsection{Minimum and non-minimum phase systems}
The ground field $k$ is now the field $\mathbb{R}$ of real numbers.
Assume that our system is left invertible. The stable or unstable
behavior of Eq. (\ref{io}), when considered as a system of
differential equations in the unknowns ${\Vect{u}}$ (${\Vect{y}}$ is
given) yields the definition of {\em minimum}, or {\em non-minimum},
{\em phase systems} (compare with \cite{isidori2}).

\section{Numerical differentiation}\label{estder}
The interested reader will find more details and references in
(\cite{nl}). We refer to (\cite{mboup}) for crucial developments
which play an important r\^ole in practical implementions.

\subsection{General principles}
\subsubsection{Polynomial signals} Consider the polynomial time
function of degree $N$
$$
x_{N} (t) = \sum_{\nu = 0}^{N} x^{(\nu)}(0) \frac{t^\nu}{\nu !}
$$
where $t \geq 0$. Its operational, or Laplace, transform (see, {\it
e.g.}, \cite{yosida}) is
\begin{equation}\label{operat}
X_N (s) = \sum_{\nu = 0}^{N} \frac{x^{(\nu)}(0)}{s^{\nu + 1}}
\end{equation}
Introduce $\frac{d}{ds}$, which is sometimes called the {\em
algebraic derivation}. Multiply both sides of Eq. (\ref{operat}) by
$\frac{d^\alpha}{ds^\alpha} s^{N + 1}$, $\alpha = 0, 1, \dots, N$.
The quantities $x^{(\nu)}(0)$, $\nu = 0, 1, \dots, N$, which satisfy
a triangular system of linear equations, with non-zero diagonal
elements,
\begin{equation}\label{triangu}
\frac{d^\alpha s^{N + 1} X_N}{ds^\alpha}  =
\frac{d^\alpha}{ds^\alpha} \left( \sum_{\nu = 0}^{N} x^{(\nu)}(0)
s^{N - \nu} \right)
\end{equation}
are said to be {\em linearly identifiable} (\cite{mfhsr}). One gets
rid of the time derivatives $s^\mu \frac{d^\iota X_N}{ds^\iota}$,
$\mu = 1, \dots, N$, $0 \leq \iota \leq N$, by multiplying both
sides of Eq. (\ref{triangu}) by $s^{- \bar{N}}$, $\bar{N}
> N$.
\begin{remarque}
The correspondence between $\frac{d^\alpha}{ds^\alpha}$ and the
product by $(- t)^\alpha$ (see, {\it e.g.}, \cite{yosida}) permits
to go back to the time domain.
\end{remarque}
\subsubsection{Analytic signals}
A time signal is said to be {\em analytic} if, and only if, its
Taylor expansion is convergent. Truncating this expansion permits to
apply the previous calculations.

\subsection{Noises}\label{noise}
{\em Noises} are viewed here as quick fluctuations around $0$. They
are therefore attenuated by low-pass filters, like iterated
integrals with respect to time.\footnote{See (\cite{ans}) for a
precise mathematical theory, which is based on {\em nonstandard
analysis}.}

\section{Model-free control: general principles}\label{cs}
It is impossible of course to give here a complete description which
would trivialize for practitioners the implementation of our control
design. We hope that numerous concrete applications will make such
an endeavor feasible in a near future.

\subsection{Local modeling}\label{bbi}
\begin{enumerate}
\item Assume that the system is left invertible. If there are more
output variables than input variable, {\it i.e.}, $p \gneqq m$, pick
up $m$ output variables, say the first $m$ ones, in order to get an
invertible square system. Eq. (\ref{F}) may be extended by writing
\begin{equation} \label{mod} \begin{array}{l}y^{(n_1)}_1 = F_1 +
\alpha_{1,1} u_1 + \dots + \alpha_{1,m} u_m
\\ \ldots \\
y^{(n_p)}_p = F_m + \alpha_{m,1} u_1 + \dots + \alpha_{p,m} u_m
\end{array}
\end{equation}
where
\begin{itemize}
\item $n_j \geq 1$, $j = 1, \dots, p$, and, most often, $n_j = 1$ or $2$;
\item $\alpha_{j,i} \in \mathbb{R}$, $i = 1, \dots, m$, $j = 1, \dots, p$,
are {\em non-physical} constant parameters, which are chosen by the
pratitioner such that $\alpha_{j,i}u_i$ and $F_j$ are of the same
magnitude.
\end{itemize}

\item In order to avoid any algebraic loop, the numerical
value of $$F_j = y^{(n_j)}_j - \alpha_{j,1} u_1 - \dots -
\alpha_{j,m} u_m$$ is given thanks to the time sampling
\begin{equation*}\label{Fk} F_{j}(\kappa) =
[y^{(n_j)}_j(\kappa)]_e - \sum_{i = 1}^m \alpha_{j,i} u_i(\kappa -
1)  \end{equation*} where {\bf $[ \bullet (\kappa) ]_e$} stands for
the estimate at the time instant $\kappa$.

\item The determination of the reference trajectories for the output variables $y_j$ is
achieved in the same way as in flatness-based control.
\end{enumerate}

\begin{remarque}\label{boucle}
It should also be pointed out that, in order to avoid algebraic
loops, it is necessary that in Eq. \eqref{io}
\begin{equation*}\label{bo}
\frac{\partial \Phi_j}{\partial y_{j}^{(n_j)}} \not\equiv 0, \quad ~
j = 1, \dots, p
\end{equation*}
It yields
$$n_j \leq \bar{N}_j$$
Numerical instabilities might appear when $\frac{\partial
\Phi_j}{\partial y_{j}^{(n_j)}}$ is closed to $0$.\footnote{This
kind of difficulties has not yet been encountered whether in the
numerical simulations presented below nor in the quite concrete
applications which were studied until now.}
\end{remarque}

\begin{remarque}\label{rema} Our control design lead with non-minimum phase systems
to divergent numerical values for the control variables $u_j$ and
therefore to the inapplicability of our techniques.
\end{remarque}

\begin{remarque}
Remember that the Equations \eqref{E} and \eqref{io} are unknown.
Verifying therefore the properties discussed in the two previous
remarks may only be achieved experimentally within the plant
operating range.
\end{remarque}

\subsection{Controllers}
Let us restrict ourselves for the sake of notations simplicity to
monovariable systems.\footnote{The extension to the multivariable
case is immediate. See Sect. \ref{ex1}.} If $\nu = 2$ in Eq.
\eqref{F}, the intelligent PID controller has already been defined
by Eq. \eqref{universal}. If $\nu = 1$ in Eq. \eqref{F}, replace Eq.
\eqref{universal} by the {\em intelligent PI controller}, or {\em
i-PI},
\begin{equation}\label{ipi}
u = - \frac{F}{\alpha} + \frac{\dot{y}^\ast}{\alpha} + K_P e + K_I
\int e
\end{equation}

\begin{remarque}
Until now we were never obliged to chose $\nu \gvertneqq 2$ in Eq.
\eqref{F}. The previous controllers \eqref{universal} and
\eqref{ipi} might then be easily extended to to the {\em generalized
proportional integral controllers}, or {\em GPIs}, of (\cite{gpi}).
\end{remarque}

\begin{remarque}
In order to improve the performances it might be judicious to
replace in Eq. \eqref{universal} or in Eq. \eqref{ipi} the unique
integral term $K_I \int e$ by a finite sum of iterated integrals
$$
 K_{I_{1}} \int e +  K_{I_{2}} \int \int e + \dots + K_{I_{\Lambda}} \int \dots \int e
$$
where
\begin{itemize}
\item $\int \dots \int e$ stands for the iterated integral of order
$\Lambda$,
\item the $K_{I_{\lambda}}$, $\lambda = 1, \dots, \Lambda$, are gains.
\end{itemize}
We get if $K_{I_{\Lambda}} \neq 0$ an {\em intelligent
PI$^{\Lambda}$D} or {\em PI$^{\Lambda}$} controller. Note also that
setting $\Lambda = 0$ is a mathematical possibility. The lack of any
integral term is nevertheless not recommended from a practical
viewpoint.
\end{remarque}

Let us briefly compare our intelligent PID controllers to classic
PID controllers:
\begin{itemize}
\item We do not need any identification procedure
since the whole structural information is contained in the term $F$
of Eq. \eqref{F}, which is eliminated thanks to Eq.
\eqref{universal}.
\item The reference trajectories, which are chosen thanks to
flatness-based methods, is much more flexible than the trajectories
which are usually utilized in the industry. Overshoots and
undershoots are therefore avoided to a large extent.
\end{itemize}

\section{Some examples of model-free control}\label{exem}
A zero-mean Gaussian white noise of variance $0.01$ is added to all
the computer simulations in order to test the robustness property of
our control design. We utilize a standard low-pass filter with
classic controllers and the principles of Sect. \ref{noise} with our
intelligent controllers.

\subsection{A stable monovariable linear system}\label{list} The
transfer function
\begin{equation}\label{stlin}
\frac{(s+2)^2}{(s+1)^3}
\end{equation}
defines a stable monovariable linear system.
\subsubsection{A classic PID controller}\label{broida} We apply the
well known method due to Bro\"{\i}da (see, {\it e.g.}, \cite{dindeleux})
by approximating the system \eqref{stlin} via the following delay
system
$$
\frac{Ke^{-\tau s}}{(Ts+1)}
$$
$K = 4$, $T = 2.018$, $\tau = 0.2424$ are obtained thanks to
graphical techniques. The gain of the PID controller are then
deduced (\cite{dindeleux}): $K_P
=\frac{100(0.4\tau+T)}{120K\tau}=1.8181$, $K_I =
\frac{1}{1.33K\tau}=0.7754$, $K_D = \frac{0.35T}{K}=0.1766$.

\subsubsection{i-PI.}
We are employing $\dot{y} = F + u$ and the i-PI controller
$$
u=-[F]_e+\dot y^\star + \mbox{\rm PI}(e)$$ where
\begin{itemize}
\item $[F]_e=[\dot y]_e-u$,
\item $y^\star$ is a reference trajectory,
\item $e = y - y^\star$,
\item $\mbox{\rm PI}(e)$ is an usual PI controller.
\end{itemize}

\subsubsection{Numerical simulations}

Fig. \ref{nom} shows that the i-PI controller behaves only slightly
better than the classic PID controller. When taking into account on
the other hand the ageing process and some fault accommodation there
is a dramatic change of situation:
\begin{itemize}
\item Fig. \ref{hold} indicates a clear cut superiority of our i-PI
controller if the ageing process corresponds to a shift of the pole
from $1$ to $1.5$, and if the previous graphical identification is
not repeated.
\item The same conclusion holds, as seen Fig. \ref{def},
if there is a  $50\%$ power loss of the control.
\end{itemize}

\begin{remarque}
This example shows that it might useless to introduce delay systems
of the type
\begin{equation}\label{retard}
T(s) e^{-Ls}, \quad ~ T \in {\mathbb{R}}(s), ~ L \geq 0
\end{equation}
for tuning classic PID controllers, as often done today in spite of
the quite involved identification procedure. It might be reminded
that
\begin{itemize}
\item the structure and the control of systems of type \eqref{retard}
have been studied in (\cite{delay}),
\item their identification with techniques stemming also from
(\cite{mfhsr}) has been studied in (\cite{belkoura2,ollivier,rw}).
\end{itemize}
\end{remarque}

\begin{remarque}\label{robuste}
This example demonstrates also that the usual mathematical criteria
for robust control become to a large irrelevant.  Let us however
point out that our control leads always to a pure integrator of
order $1$ or $2$, for which the classic frequency techniques (see,
{\it e.g.}, \cite{book,franklin,zamb}) might still be of some
interest.
\end{remarque}

\begin{remarque}\label{diag}
As also shown by this example some fault accommodation may also be
achived without having recourse to a general theory of diagnosis.
\end{remarque}

\subsection{A monovariable linear system with a large spectrum}\label{large}
With the system defined by the transfer function
\begin{equation*}\label{transf}
\frac{s^5}{(s + 1)(s + 0.1)(s + 0.01)(s - 0.05)(s - 0.5)(s - 5)}
\end{equation*}
We utilize $\dot{y} = F + u$. A i-PI controller provides the
stabilization around a reference trajectory. Fig. \ref{lx} exhibits
an excellent tracking.

\subsection{A multivariable linear system}\label{ex1} Introduce the
transfer matrix
{\tiny \begin{equation*}\label{matransf}
\begin{pmatrix}\frac{s^3}{(s+0.01)(s+0.1)(s-1)s}& 0\\
\frac{s+1}{(s+0.003)(s-0.03)(s+0.3)(s+3)}&\frac{s^2}
{(s+0.004)(s+0.04)(s-0.4)(s+4)}\end{pmatrix}
\end{equation*}}
\noindent For the corresponding system we utilize after a few
attempts Eq. (\ref{mod}) with the following decoupled form
\begin{equation*}\label{li}
\dot{y}_1=F_1 + 10u_1 \quad ~ \ddot{y}_2=F_2 + 10u_2
\end{equation*}
The stabilization around a reference trajectory $(y_{1}^{\ast},
y_{2}^{\ast})$ is ensured by the multivariable i-PID controller
\begin{equation*}\label{cor}
\begin{array}{l}
u_1 =\frac{1}{10}\left( \dot{y}^{\ast}_1 - F_1 + K_{P1} e_1 +
K_{I1} \int e_1 +K_{D1} \dot{e}_1 \right)\\
u_2=\frac{1}{10}\left(\ddot{y}^{\ast}_2 -F_2 + K_{P2} e_2 +
K_{I2} \int e_2 +K_{D2} \dot{e}_2 \right)
\end{array}
\end{equation*}
where
\begin{itemize}
\item $e_1 = y^{\ast}_1 - y_1$, $e_2 = y^{\ast}_2 - y_2$;
\item  $K_{P1}=1$, $K_{I1}=K_{D1}=0$,
$K_{P2}=K_{I2}=50$, $K_{D2}=10$.
\end{itemize}
The performances displayed on Fig. \ref{fig_Ml1} and \ref{fig_Ml2}
are excellent. Fig. \ref{fig_Ml2}-(b) shows the result if we would
set $F_1=F_2=0$: it should be compared with Fig. \ref{fig_Ml2}-(a).

\begin{remarque}
Model reduction is often utilized for the kind of systems studied in
Sect. \ref{large} and Sect. \ref{ex1} (see, {\it e.g.},
\cite{antoulas,ob}).
\end{remarque}

\subsection{An unstable monovariable nonlinear system \label{NLI1}}

\subsubsection{i-PID}
For $\dot{y} - y = u^3$ we utilize for Eq. (\ref{F}) the local model
$\dot y =F+u$. The stabilization around a reference trajectory
$y^\ast$ is provided by the i-PI controller
\begin{equation}\label{satur}
u=-F+\dot y^\star + K_Pe + K_I\int e
\end{equation}
where $K_P=-2$, $K_I=-1$. The simulations displayed in Fig.
\ref{fig_NLA} are excellent.


\subsubsection{Anti-windup}\label{antiemb}
We now assume that $u$ should satisfy the following constraints $- 2
\leq u \leq 0.4$. The performances displayed by Fig. \ref{fig_NLAsw}
are mediocre if an anti-windup is not added to the classic part of
the i-PI controller. Our solution is elementary\footnote{This is a
well covered subject in the literature (see, {\it e.g.},
\cite{bohn,hippe,peng}.}: as soon as the control variable gets
saturated, the integral $\int e$ in Eq. \eqref{satur} is maintained
constant.\footnote{Better performances would be easily reached, as
flatness-based control is teaching us, with a modified reference
trajectory.}
\subsection{Ball and beam}
Fig. \ref{BB} displays the famous ball and beam example, which obeys
to the equation\footnote{The sine function which appears in that
equation takes us outside of the theory sketched in Sect.
\ref{algebra}. This difficulty may be easily circumvented by
utilizing $\mbox{\rm tg} \frac{u}{2}$ (see \cite{flmr}).} $\ddot y =
By\dot u^2-BG \sin u$, where $u = \theta$ is the control variable.
This monovariable system, which is not linearizable by a static
state feedback, is therefore not flat. It is thus difficult to
handle.\footnote{A quite large literature has been devoted to this
example (see, {\it e.g.}, (\cite{lozano,koko,sastry}) for advanced
nonlinear techniques, and (\cite{zhang}) for neural networks). Let
us add that all numerical simulations in those references are given
without any corrupting noise.}

We have chosen for Eq. (\ref{F}) $\ddot{y} = F + 100 u$. In order to
satisfy as well as possible the experimental conditions,  the
control variable is saturated: $-\pi/3 < u  < \pi/3$ and $- \pi <
\dot u <\pi$. Fig. \ref{poly} and Fig. \ref{sinus} display two types
of trajectories: a B\'{e}zier polynomial and a sine function. We obtain
in both cases excellent trackings thanks to an i-PID controller. In
the Figures \ref{poly}-(b), \ref{sinus}-(b), \ref{poly}-(c),
\ref{sinus}-(c) the control variable and the estimations of $F$ are
presented in the noiseless case. Compare with the Figures
\ref{poly}-(d) and \ref{sinus}-(d)) where a corrupting noise is
added.

\subsection{The three tank example}

The three tank example in Fig. \ref{schem} is quite popular in
diagnosis.\footnote{See (\cite{zeitz}) for details and references.
This paper was presenting  apparently for the first time the
diagnosis, the control and the fault accommodation of a nonlinear
system with uncertain parameters. See also \cite{nl}.} It obeys to
the equations:
$$
    \begin{cases}\begin{array}{ll}
    \dot{x}_1 =&-C_1\text{sign}(x_1 - x_3)\sqrt{|x_1 - x_3|}+u_1 /S\\
    \dot{x}_2 =&C_3\text{sign}(x_3 - x_2)\sqrt{|x_3 - x_2|}\\&-C_2\text{sign}(x_2)\sqrt{|x_2|}
    +u_2 /S\\
    \dot{x}_3 =&C_1\text{sign}(x_1 - x_3)\sqrt{|x_1 - x_3|}\\&-C_3\text{sign}(x_3 - x_2)\sqrt{|x_3 - x_2|}\\
    y_{1} =&x_1\\
    y_{2} =&x_2\\
    y_{3} =&x_3\\
    \end{array}\end{cases}
$$
where
$$
\begin{array}{l}
C_n = (1/S).\mu_n.S_p\sqrt{2g}, n = 1, 2, 3;\\
S=0.0154~m~\text{(tank section)};\\
Sp=5.10^{-5}~m~\text{(pipe section between the tanks)};\\
g=9.81~m.s^{-2}~\text{(gravity)};\\
\mu_1=\mu_3=0.5\text{,} ~ \mu_2=0.675 ~\text{(viscosity
coefficients)}.
\end{array}
$$
As often in industry we utilize a zero-hold control (see Fig.
\ref{fig_3cuves}- (c)). A decoupled Eq. (\ref{mod}) is employed here
as we already did in Sect. \ref{ex1}: $\dot y_i=F_i+ 200 u_i$, $i =
1, 2$. Fig. \ref{fig_3cuves}-(a) displays the trajectories tracking.
The derivatives estimation in Fig. \ref{fig_3cuves}-(b) is excellent
in spite of the additive corrupting noise. The nominal controls
(Fig. \ref{fig_3cuves}-(c)) are not very far from those we would
have computed with a flatness-based viewpoint (see \cite{zeitz}). We
also utilize the following i-PI controllers

{\small
\begin{equation*}\label{pid} u_i =
\frac{1}{200}\left(\dot{y}_{i}^\ast-F_i + 10 e_i + 2.10^{-2} \int
e_i \right) \quad i = 1, 2 \end{equation*}}

\noindent where $y_{i}^\ast$ is the reference trajectory, $e_i =
y_{i}^\ast - y_i$. In order to get a good estimate of $e_i$ we are
denoising $y_i$ (see Fig. \ref{fig_3cuves}-(d)) according to the
techniques of Sect. \ref{estder}.

\section{Control with a restricted model}\label{restricted}
\subsection{General principles}\label{basis}
\subsubsection{Flatness}\label{plat}
The system (\ref{io}) is assumed to be square, {\it i.e.}, $m = p$,
and flat. Moreover ${\Vect{y}} = (y_1, \dots, y_m)$ is assumed to be
a flat output, {\i.e.}, $\bar{M}_j = 0$, $j = 1, \dots, m$. It
yields locally
\begin{equation} \label{flat}
u_j = \Psi_j({\Vect{y}}, \dots, {\Vect{y}}^{(\bar{N}_j)}), \quad ~ j
= 1, \dots, m
\end{equation}
Flatness-based control permits to select easily an efficient
reference trajectory ${\Vect{y}}^\star$ to which corresponds via Eq.
\eqref{io} and Eq. \eqref{flat} the open loop control
${\Vect{u}}^\star$. Let ${\Vect{e}} = {\Vect{y}} - {\Vect{y}}^\star$
be the tracking error. We assume the existence of a feedback
controller ${\Vect{u}}_{\mbox{\rm feedback}} ({\Vect{e}})$ such that
\begin{equation}\label{feed}
{\Vect{u}} = {\Vect{u}}^\star + {\Vect{u}}_{\mbox{\rm feedback}}
({\Vect{e}})
\end{equation}
ensures a stable tracking around the reference trajectory.

\subsubsection{Intelligent controllers}
Replace Eq. (\ref{io}) by
\begin{equation} \label{ion}
\Phi_j({\Vect{y}}, \dots, {\Vect{y}}^{(\bar{N}_j)}, {\Vect{u}},
\dots, {\Vect{u}}^{(\bar{M}_j)} ) + G_j = 0
\end{equation}
where the $G_j$, $j = 1, \dots, m$, stand for the unmodeled parts.
Eq. \eqref{flat} becomes then
\begin{equation} \label{flat-mod}
u_j = \Psi_j({\Vect{y}}, \dots, {\Vect{y}}^{(\bar{N}_j)}) + H_j,
\quad ~ j = 1, \dots, m
\end{equation}
where $H_j \neq G_j$ in general. Thanks to Eq. \eqref{flat-mod},
$H_j$ is estimated in the same way as the fault variables and the
unknown perturbations are in (\cite{nl}). Consider again
${\Vect{y}}^\star$ and ${\Vect{u}}^\star$ as they are defined above.
The {\em intelligent controller}, or {\em i-controller}, follows
from Eq. \eqref{feed}
\begin{equation*}\label{feed-mod}
{\Vect{u}} = {\Vect{u}}^\star + \left( \begin{array}{c} H_1 \\
\vdots \\ H_m \end{array} \right) + {\Vect{u}}_{\mbox{\rm feedback}}
({\Vect{e}})
\end{equation*}
It ensures tracking stabilization around the reference trajectory.

\subsection{Frictions and nonlinearities}\label{spring}

A point mass $m$ at the end of a spring of length $y$ obeys to the
equation
\begin{equation}\label{ressort}m\ddot{y}= -{\mathcal{K}}(y) +
{\mathcal{F}}(\dot{y})-d\dot y+F_{\tiny{\mbox{\rm
ext}}}\end{equation} where
\begin{itemize}
\item $F_{\tiny{\mbox{\rm ext}}}=u$ is the control variable;
\item $d$ and ${\mathcal{F}}(\dot y)$ are due to complex friction phenomena;
\item ${\mathcal{K}}(y) = k_1 y + k_3 y^3$ exhibits a cubic nonlinearity of Duffing type;
\end{itemize}
The mass $m=0.5$ is known;  there is a possible error of $33\%$ for
$k_1 =3$ and we utilize $\hat{k}_1 = 2$; $d$ and $k_3$, which are
unknown, are equal to $5$ and $10$ in the numerical simulations. For
the frictions,\footnote{There is a huge literature in tribology
where various possible friction models are suggested (see, {\it
e.g.}, in control (\cite{ols,lille}). Those modelings are bypassed
here.} we have chosen for the sake of computer simulations the well
known model due to \cite{tustin}. Fig. \ref{fig_SLF}-(a) exhibits
its quite wild behavior when the sign of the speed is changing.

\subsubsection{A classic PID controller}
The PID controller is tuned only thanks to the restricted model $
m\ddot y=-\hat{k}_1 y+u$. Its gains are determined in such a way
that all the poles of the closed-loop system are equal to $-3$: $K_P
= -\hat{k}_1+27m$, $K_I = -27m$, $K_D = 9m$.
\subsubsection{The corresponding i-PID controller}
Pick up a reference trajectory $y^\star$. Set
$$u^\star=m\ddot y^\star+\hat{k}_1 y^\star$$
Our i-PID controller is given by
\begin{equation}\label{equfrot}
F_{\tiny{\mbox{\rm ext}}} = u =u^\star-[{\mathcal{G}}]_e + \mbox{\rm
PID}(e)
\end{equation}
where
\begin{itemize}
\item ${\mathcal{G}} ={\mathcal{F}}(\dot y)-(k_1-\hat k_1)y-k_3y^3-d\dot
y$, which stands for the whole set of unknown effects, is estimated
via
$$[{\mathcal{G}}]_e = m[\ddot y]_e+\hat{k} [y]_e- F_{\tiny{\mbox{\rm
ext}}}$$ which follows from Eq. \eqref{ressort} ($[y]_e$ and $[\ddot
y]_e$ are the denoised output variable and its denoised
$2^{nd}$-order derivative -- see Fig. \ref{fig_SLF}-(d,f));
\item  $\mbox{\rm PID}(e)$, $e = y - y^\star$, is the above classic PID controller.
\end{itemize}

\subsubsection{Numerical simulations}
The performances of our i-PID controller \eqref{equfrot}, which are
displayed in Fig. \ref{fig_SLF}-(c,d), are excellent. When compared
to the Figures
\begin{itemize}
\item \ref{fig_SLF}-(e,f),  where
\begin{itemize}
\item flatness-based control is employed for determined the open-loop
output and input variables,
\item the loop is closed via a classic PID controller, which does not take into
account the unknown effects;
\end{itemize}
\item \ref{fig_SLF}-(g,h), where only a classic PID controller
is used, without any flatness-based control;
\end{itemize}
the superiority of our control design is obvious. This superiority
is increasing with the friction.


\subsection{Non-minimum phase systems}
Consider the transfer function
$$\frac{s-a}{s^2 -
(b+c)s + bc}$$ $a, b, c \in \mathbb{R}$ are respectively its zero
and its two poles. The corresponding input-output system is
non-minimum phase if $a > 0$. The controllable and observable
state-variable representation
\begin{equation}\label{deph}
\begin{cases}\dot x_1=x_2\\
\dot x_2=(b+c)x_2-bc x_1+u\\ y=x_2-ax_1\end{cases}
\end{equation}
shows that $z=x_1$ is a flat output. The flat output is therefore
not the measured output, as we assumed in Sect. \ref{plat}. Our
control design has therefore to be modified.
\subsubsection{Control of the exact model}\label{exact} To a nominal
flat output\protect\footnote{This Section, which is based on
previous studies (\cite{mfrm,gpi}), should make the reading of Sect.
\ref{inexact} easier.} $z^\star$ corresponds a nominal control
variable
$$u^\star=\ddot z^\star-(b+c)\dot z^\star+bc
z^\star$$ and a nominal output variable
\begin{equation}\label{z}
y^\star=\dot{z}^\star-a z^\star
\end{equation}
Introduce the {\em GPI} controller (\cite{gpi})
\begin{equation}\begin{array}{lc}u=& u^\star+\gamma\int(u-u^\star)+K_P (y-y^\star) \\ &+
K_I\int(y-y^\star)+K_{II}\int \!\! \int(y-y^\star)
\end{array}
\label{pnmu}\end{equation} where the coefficients $\gamma, K_P, K_I,
K_{II} \in \mathbb{R}$ are chosen in order to stabilize the error
dynamics $e= z - z^\star$. Excellent performances are displayed in
Figures \ref{pnm0}-(a) and (b), where $a=1$, $b=-1$, $c=-0.5$, even
with an additive corrupting noise. See Figures \ref{pnm0}-(c) and
(d) for $y^\star$ and $z^\star$ which is calculated by integrating
Eq. (\ref{z}) back in time.

\subsubsection{Unmodeled effects}\label{inexact}
The second line of Eq. (\ref{deph}) may be written again as
$$ \dot{x}_2=(b+c)x_2-bc x_1+u+ \varpi$$ where $\varpi$ stands for the unmodeled effects,
like frictions or an actuator's fault. Replace the nominal control
variable $u^\star$ of Sect. \ref{exact} by
$$u_{\tiny{\mbox{\rm pert}}}^\star=u^\star-[\varpi]_e$$
where $[\varpi]_e$ is the estimated value of $\varpi$, which is
given by
$$[\varpi]_e=-\left(\frac{[\ddot y]_e-(b+c)[\dot y]_e+bc[y]_e-[\dot u]_e}{a}+u\right)$$
Moreover,
$$\begin{array}{lc}{[\dot u]_e} =& \dot u^\star+\gamma (u-u^\star)+
K_P([\dot y]_e-\dot y^\star)\\ &+K_I([y]_e-y^\star)
+K_{II}\int([y]_e-y^\star)\end{array}$$ follows from Eq.
(\ref{pnmu}).

Start with the Figures \ref{pnm1} and \ref{pnm2}, where
$\varpi=-0.5$, $a=1$, $b=-1$, $c=-0.5$.  Fig. \ref{pnm1}-(b), where
the nominal  control is left unmodified, and Fig. \ref{pnm1}-(e),
where it is modified, demonstrate a clear-cut superiority of our
approach, even with an additive corrupting noise.

In Fig. \ref{pnm3}, $\varpi$ is no more assumed to be constant, but
equal to $-0.1 \dot{y}$. In the numerical simulations, $a=2$, $b=-1$
$c=1$. If $\varpi$ is not estimated, it influences the tracking
quite a lot even if its amplitude is weak. When $\varpi$ is
estimated on the other hand, the results are excellent, even with an
additive corrupting noise.

\section{Conclusion}\label{conclusion}
The results which were already obtained with our intelligent PID
controllers lead us to the hope that they will greatly improve the
practical applicability and the performances of the classic PIDs, at
least for all finite-dimensional systems which are known to be
non-minimum phase within their operating range:
\begin{itemize}
\item the tuning of the gains of i-PIDs is straightforward since
\begin{itemize}
\item the unknown part is eliminated,
\item the control design boils down to a pure integrator of order
$1$ or $2$;
\end{itemize}
\item the identification techniques for implementing classic
PID regulators, which are often imprecise and difficult to handle,
are becoming obsolete.
\end{itemize}

Model-free control and the control with a restricted model seem to
question the very principles of modeling in applied sciences, at
least when one wishes to control some concrete plant. This might be
a fundamental ``epistemological'' change, which needs of course to
be further discussed and analyzed. A natural extension to
uncontrolled systems is being developed via various questions in
financial engineering: see already the preliminary studies in
(\cite{coventry,fes,malo}).



\begin{thebibliography}{xx}

\bibitem[Antoulas(2005)]{antoulas}
A.C. Antoulas.
\newblock Approximation of Large-Scale Dynamical Systems.
\newblock SIAM, 2005.

\bibitem[{\AA}str\"om \& H\"agglund(2006)]{pid}
K.J. {\AA}str\"om, T. H\"agglund.
\newblock Advanced PID Control.
\newblock Instrument Soc. Amer., 2006.

\bibitem[{\AA}str\"om, Persson \& Hang(1992)Hang, Persson and Ho]{intelligent}
K.J. {\AA}str\"om, C.C. Hang, P. Persson, W.K. Ho.
\newblock Towards intelligent PID control.
\newblock \emph{Automatica}, 28:\penalty0 1--9, 1992.

\bibitem[{\AA}str\"om \& Murray(2008)]{book}
K.J. {\AA}str\"om, R.M. Murray.
\newblock Feedback Systems: An Introduction for Scientists and Engineers.
\newblock Princeton University Press, 2008.



\bibitem[Belkoura, Richard \& Fliess(2009)Belkoura, Richard and Fliess]{belkoura2}
L. Belkoura, J.-P. Richard, M. Fliess.
\newblock Parameters estimation of systems with delayed and structered entries.
\newblock \emph{Automatica}, 2009. Available at {\tt
http://hal.inria.fr/inria-00343801/en/}.

\bibitem[Besan\c{c}on-Voda \& Gentil(1999)Besan\c{c}on-Voda and Gentil]{voda}
A. Besan\c{c}on-Voda, S. Gentil.
\newblock R\'{e}gulateurs PID analogiques et num\'{e}riques.
\newblock \emph{Techniques de l'Ing\'{e}nieur}, R7416, 1999.

\bibitem[Bohn \& Atherton(1995)]{bohn}
C. Bohn, D.P. Atherton.
\newblock An analysis package comparing PID anti-windup strategies.
\newblock \emph{IEEE Control Syst. Magaz.}, 15:\penalty0 34--40, 1995.

\bibitem[Bourdais, Fliess, Join \& Perruquetti(2007)Bourdais, Fliess, Join and Perruquetti]{bourdais}
R. Bourdais, M. Fliess, C. Join, W. Perruquetti.
\newblock Towards a model-free output tracking of switched nonlinear systems
\newblock \emph{Proc. 7$^{th}$ IFAC Symp. Nonlinear Control Systems},
Pretoria, 2007. Available at {\tt
http://hal.inria.fr/inria-00147702/en/}.

\bibitem[Chambert-Loir(2005)]{cl}
A. Chambert-Loir.
\newblock Alg\`{e}bre corporelle.
\newblock \emph{\'Editions \'Ecole Polytechnique}, 2005. English translation:
{\em A Field Guide to Algebra}, Springer, 2005.


\bibitem[Choi, d'Andr\'{e}a-Novel, Fliess \& Mounier(2009)Choi, Andr\'{e}a-Novel,
Choi, Fliess and Mounier]{choi} S. Choi, B. d'Andr\'{e}a-Novel, M.
Fliess, H. Mounier.
\newblock  Model-free control of automotive engine and brake for stop-and-go
scenario. \emph{Proc. 10$^{th}$ IEEE Conf. Europ. Control Conf.},
Budapest, 2009. Soon available at {\tt http://hal.inria.fr/}.

\bibitem[Dattaet, Ho \& Bhattacharyya(2000)Datta, Ho and Bhattacharyya]{dat}
A. Datta, M.T. Ho, S.P. Bhattacharyya.
\newblock Structure and Synthesis of PID Controllers.
\newblock Springer, 2000.

\bibitem[Delaleau(2002)]{delaleau}
E. Delaleau.
\newblock Alg\`ebre diff\'erentielle.
\newblock In J.-P. Richard, editor, \emph{Math\'ematiques pour les syst\`emes dynamiques},
vol. 2, chapter 6,  pp. 245--268, Herm\`es, 2002.

\bibitem[Delaleau(2008)]{delaleau2}
E. Delaleau.
\newblock Classical electrical engineering questions in the light of Fliess's differential algebraic framework of non-linear control systems.
\newblock \emph{Int. J. Control}, 81:\penalty0 380--395, 2008.



\bibitem[Dindeleux(1981)]{dindeleux}
D. Dindeleux.
\newblock Technique de la regulation industrielle.
\newblock Eyrolles, 1981.


\bibitem[Fantoni \& Lozano(2002)]{lozano}
I. Fantoni, R. Lozano.
\newblock Non-linear control for underactuated mechanical systems.
\newblock Springer, 2002.

\bibitem[Fliess(2006)]{ans}
M. Fliess.
\newblock Analyse non standard du bruit.
\newblock \emph{C.R. Acad. Sci. Paris Ser. I}, 342:\penalty0 797--802, 2006.


\bibitem[Fliess \& Join(2008a)]{esta}
M. Fliess, C. Join.
\newblock Commande sans mod\`{e}le et commande \`{a} mod\`{e}le restreint.
\newblock \emph{e-STA}, 5 (n$^\circ$ 4):\penalty0 1-23, 2008a. Available at {\tt
http://hal.inria.fr/inria-00288107/en/}.

\bibitem[Fliess \& Join(2008b)]{coventry}
M. Fliess, C. Join.
\newblock Time series technical analysis via new fast estimation
methods: a preliminary study in mathematical finance.
\newblock {\em Proc. 23$^{rd}$ IAR Workshop Advanced Control
Diagnosis (IAR-ACD08)}, Coventry, 2008b. Available at {\tt
http://hal.inria.fr/inria-00338099/en/}.

\bibitem[Fliess \& Join(2009a)]{fes}
M. Fliess, C. Join.
\newblock A mathematical proof of the existence of trends in
financial time series.
\newblock \emph{Proc. Int. Conf. Systems Theory: Modelling,
Analysis and Control}, Fes, 2009a. Available at {\tt
http://hal.inria.fr/inria-00352834/en/}.

\bibitem[Fliess \& Join(2009b)]{malo}
M. Fliess, C. Join.
\newblock Towards new technical indicators for trading systems and
risk management.
\newblock \emph{Proc. 15$^{th}$ IFAC Symp. System Identif. (SYSID 2009)},
Saint-Malo, 2009b. Available at {\tt
http://hal.inria.fr/inria-00370168/en/}.



\bibitem[Fliess, Join \& Sira-Ram\'{\i}rez(2005)Fliess, Join and Sira-Ram\'{\i}rez]{zeitz}
M. Fliess, C. Join, H. Sira-Ram\'{\i}rez.
\newblock Closed-loop fault-tolerant control for uncertain nonlinear systems.
\newblock In T. Meurer, K. Graichen and E.D. Gilles editors, \emph{Control and Observer Design for
Nonlinear Finite and Infinite Dimensional Systems}, Lect. Notes
Control Informat. Sci., volume 322, pages 217--233, Springer, 2005.


\bibitem[Fliess, Join \& Sira-Ram\'{\i}rez(2008) Fliess, Join and Sira-Ram\'{\i}rez]{nl}
M. Fliess, C. Join, H. Sira-Ram\'{\i}rez.
\newblock Non-linear estimation is easy.
\newblock \emph{Int. J. Modelling Identification Control}, vol.
4:\penalty0 12--27 (2008). Available at {\tt
http://hal.inria.fr/inria-00158855/en/}.

\bibitem[Fliess, L\'{e}vine, Martin \& Rouchon(1995) Fliess, L\'evine, Martin and Rouchon]{flmr}
M. Fliess, J. L\'evine, P. Martin, P. Rouchon.
\newblock Flatness and defect of non-linear systems: introductory theory and examples.
\newblock \emph{Int. J. Control}, 61:\penalty0 1327--1361, 1995.

\bibitem[Fliess \& Marquez(2000)Fliess and Marquez]{mfrm}
M. Fliess, R. Marquez.
\newblock Continuous-time linear predictive control and flatness:
a module-theoretic setting with examples.
\newblock \emph{Int. J. Control}, 73:\penalty0 606--623, 2000.

\bibitem[Fliess, Marquez, Delaleau \& Sira-Ram\'{\i}rez(2002)Fliess, Marquez, Delaleau and Sira-Ram\'{\i}rez]{gpi}
M. Fliess, R. Marquez, E. Delaleau, H. Sira-Ram\'{\i}rez.
\newblock Correcteurs proportionnels-int\'egraux g\'en\'eralis\'es ",
\newblock \emph{ESAIM Control Optim. Calc. Variat.}, 7:\penalty0 23--41, 2002.

\bibitem[Fliess, Marquez \& Mounier(2002)Fliess, Marquez and Mounier]{delay}
M. Fliess, R. Marquez, H. Mounier.
\newblock An extension of predictive control, PID regulators and Smith predictors to some delay systems.
\newblock \emph{Int. J. Control}, 75:\penalty0 728--743, 2002.

\bibitem[Fliess \& Sira-Ram\'{\i}rez(2003)Fliess and Sira-Ram\'{\i}rez]{mfhsr}
M. Fliess, H. Sira-Ram\'{\i}rez.
\newblock An algebraic framework for linear identification.
\newblock \emph{ESAIM Control Optim. Calc. Variat.}, 9:\penalty0 151--168, 2003.


\bibitem[Franklin, Powell \& Emami-Naeini(2002)Franklin, Powell and Emami-Naeini]{franklin}
G.F. Franklin, J.D. Powell, A. Emami-Naeini.
\newblock Feedback Control of Dynamic Systems.
\newblock 4$^{th}$ edition, Prentice Hall, 2002.

\bibitem[G\'{e}douin, Join, Delaleau, Bourgeot, Chirani \& Calloch(2008)G\'{e}douin, Join, Delaleau, Bourgeot, Chirani and Calloch]{brest}
P.-A. G\'{e}douin, C. Join, E. Delaleau, J.-M. Bourgeot, S.A. Chirani,
S. Calloch.
\newblock Model-free control of shape memory alloys antagonistic actuators.
\newblock \emph{Proc. 17$^{th}$ IFAC World Congress (WIFAC-2008)}, Seoul, 2008.
Available at {\tt http://hal.inria.fr/inria-00261891/en/}.



\bibitem[Hippe(2006)Hippe]{hippe}
P. Hippe.
\newblock Windup in Control - Its Effects and Their Prevention.
\newblock Springer, 2006.

\bibitem[Isidori(1999)Isidori]{isidori2}
A. Isidori.
\newblock Nonlinear Control Systems II.
\newblock Springer, 1999.

\bibitem[Johnson \& Moradi(2005)Johnson and Moradi]{john}
M.A. Johnson, M.H. Moradi.
\newblock PID Control: New Identification and Design Methods.
\newblock Springer, 2005.

\bibitem[Join, Masse \& Fliess(2008)Join, Masse and Fliess]{poitiers}
C. Join, J. Masse, M. Fliess.
\newblock \'Etude pr\'{e}liminaire d'une commande sans mod\`{e}le pour papillon de moteur.
\newblock \emph{J. europ. syst. automat.}, 42:\penalty0 337--354,
2008. Available at {\tt http://hal.inria.fr/inria-00187327/en/}.

\bibitem[Kerschen, Worden, Vakakis \& Golinval(2006)Kerschen, Worden,
Vakakis and Golinval]{aut2}
G. Kerschen, K. Worden, A.F. Vakakis, J.-C. Golinval.
\newblock Past, present and future of nonlinear system identification in structural dynamics.
\newblock \emph{Mech. Systems Signal Process.}, 20:\penalty0 505--592, 2006.

\bibitem[Hauser, Sastry \& Kokotovic(1992)Hauser, Sastry and Kokotovi\'c]{koko}
J. Hauser, S. Sastry, P. Kokotovi\'c.
\newblock Nonlinear control via approximate input-output linearization:
The ball and beam example.
\newblock \emph{ IEEE Trans. Automat. Control}, 37:\penalty0 392--398, 1992.

\bibitem[Kolchin(1973)Kolchin]{kolchin}
E.R. Kolchin.
\newblock Differential Algebra and Algebraic Groups.
\newblock Academic Press, 1973.


\bibitem[Lequesne(2006) Lequesne]{leq}
D. Lequesne.
\newblock R\'{e}gulation P.I.D. : analogique - num\'{e}rique - floue.
\newblock Herm\`{e}s, 2006.


\bibitem[Mboup, Join \& Fliess(2009) Mboup, Join and Fliess]{mboup}
M. Mboup, C. Join, M. Fliess.
\newblock Numerical differentiation with
annihilators in noisy environment.
\newblock \emph{Numer. Algo.}, 50, 2009. DOI:
{ \tt 10.1007/s11075-008-9236-1}.

\bibitem[Nuninger, Perruquetti \& Richard(2006)Nuninger, Perruquetti and Richard]{lille}
W. Nuninger, W. Perruquetti, J.-P. Richard.
\newblock Bilan et enjeux des mod\`{e}les de frottements: tribologie
et contr\^{o}le au service de la s\'{e}curit\'{e} des transports, \newblock
\emph{Actes 5$^e$ Journ\'{e}es Europ. Freinage (JEF'2006)}, Lille, 2006.
Available at {\tt http://hal.inria.fr/inria-00192425/en/}.

\bibitem[Obinata \& Anderson(2001)Obinata and Anderson]{ob}
G. Obinata, B.D.O. Anderson.
\newblock Model Reduction for Control Systems.
\newblock Springer, 2001.


\bibitem[O'Dwyer(2006)O'Dwyer]{od}
A. O'Dwyer.
\newblock Handbook of PI and PID Controller Tuning Rules.
\newblock 2$^{nd}$ ed., Imperial College Press, 2006.

\bibitem[Ollivier,  Moutaouakil \& Sadik(2007)Ollivier, Moutaouakil and Sadik]{ollivier}
F. Ollivier, S. Moutaouakil, B. Sadik.
\newblock Une m\'{e}thode d'identification pour un syst\`{e}me lin\'{e}aire \`{a} retards.
\newblock \emph{C.R. Acad. Sci. Paris Ser. I}, 344:\penalty0 709--714, 2007.

\bibitem[Olsson, {\AA}str\"{o}m, Canudas de Wit, G\"{a}fvert \& Lischinsky(1998)Olsson, {\AA}str\"{o}m, Canudas de Wit, G\"{a}fvert and Lischinsky]{ols}
H. Olsson, K. J. {\AA}str\"{o}m, C. Canudas de Wit, M. G\"{a}fvert, P.
Lischinsky.
\newblock Friction models and friction compensation.
\newblock \emph{Europ. J. Control}, 4:\penalty0 176--195, 1998.

\bibitem[Peng, Vrancic \& Hanus(1996)Peng, Vrancic and Hanus]{peng}
Y. Peng, D. Vrancic, R. Hanus.
\newblock Anti-windup, bumpless, and conditioned transfer transfer techniques
for PID controllers.
 \newblock \emph{IEEE Control Syst. Magaz.}, 16:\penalty0 48--57,1996.


\bibitem[Rotella \& Zambettakis(2007)Rotella and Zambettakis]{rotella}
F. Rotella, I. Zambettakis.
\newblock Commande des syst\`{e}mes par platitude.
\newblock \emph{Techniques de l'ing\'{e}nieur}, S7450, 2007.

\bibitem[Rotella \& Zambettakis(2008)Rotella and Zambettakis]{zamb}
F. Rotella, I. Zambettakis.
\newblock Automatique \'{e}l\'{e}mentaire.
\newblock Herm\`{e}s, 2008.


\bibitem[Rudolph \& Woittennek(2007)Rudolph and Woittennek]{rw}
J. Rudolph, F. Woittennek.
\newblock Ein algebraischer Zugang zur Parameteridentifkation in
linearen unendlichdimensionalen Systemen.
\newblock \emph{at--Automatisierungstechnik}, 55:\penalty0 457--467, 2007.

\bibitem[Sastry(1999) Sastry]{sastry}
S. Sastry.
\newblock Nonlinear Systems.
\newblock Springer, 1999.

\bibitem[Shinskey(1996)Shinskey]{shin}
F.G. Shinskey.
\newblock Process Control Systems - Application, Design, and Tuning.
\newblock 4$^{th}$ ed., McGraw-Hill, 1996.

\bibitem[Sira-Ram\'{\i}rez \& Agrawal(2004)Sira-Ram\'{\i}rez and Agrawal]{hsr}
H. Sira-Ram\'{\i}rez, S. Agrawal.
\newblock Differentially Flat Systems.
\newblock Marcel Dekker, 2004.


\bibitem[Sj\"{o}berg, Zhang, Ljung, Benveniste, Delyon,
Glorennec, Hjalmarsson \& Juditsky(1995)Sj\"{o}berg, Zhang, Ljung,
Benveniste, Delyon, Glorennec, Hjalmarsson and Juditsky]{aut1} J.
Sj\"{o}berg, Q. Zhang, L. Ljung, A. Benveniste, B. Delyon, P.-Y.
Glorennec, H. Hjalmarsson, A. Juditsky.
\newblock Nonlinear black-box modeling in system identification: a unified overview.
\newblock \emph{Automatica}, 31:\penalty0 1691--1724, 1995.


\bibitem[Tustin(1947)Tustin]{tustin}
A. Tustin.
\newblock The effect of backlash and of speed dependent friction on the stability of closed-cycle control systems.
\newblock \emph{J. Instit. Elec. Eng.}, 94:\penalty0 143--151, 1947.


\bibitem[Villagra, d'Andr\'{e}a-Novel, Fliess \& Mounier(2008a)Villagra, Andr\'{e}a-Novel,
Fliess and Mounier]{vil1}
J. Villagra, B. d'Andr\'{e}a-Novel,  M. Fliess, H. Mounier.
\newblock Estimation of longitudinal and lateral vehicle velocities$:$ an algebraic approach.
\newblock \emph{Proc. Amer. Control Conf. (ACC-2008)}, Seattle, 2008a. Available at {\tt
http://hal.inria.fr/inria-00263844/en/}.

\bibitem[Villagra, d'Andr\'{e}a-Novel, Fliess \& Mounier(2008b)Villagra, Andr\'{e}a-Novel,
Choi, Fliess and Mounier]{vil} J. Villagra, B. d'Andr\'{e}a-Novel, M.
Fliess, H. Mounier.
\newblock  Robust grey box closed-loop stop-and-go control. \emph{Proc. 47$^{th}$ IEEE
Conf. Decision Control}, Cancun, 2008. Available at {\tt
http://hal.inria.fr/inria-00319591/en/}.

\bibitem[Visioli(2006)Visioli]{vis}
A. Visioli.
\newblock Practical PID Control.
\newblock Springer, 2006.

\bibitem[Wang, Ye, Cai \& Hang(2008)Wang, Ye, Cai and Hang]{wang}
Q.-J. Wang, Z. Ye, W.-J. Cai and C.-C. Hang.
\newblock PID Control for Multivariable Processes.
\newblock Springer, 2008.

\bibitem[Yosida(1984)Yosida]{yosida}
K. Yosida.
\newblock Operational Calculus$:$ A Theory of Hyperfunctions.
\newblock Springer, 1984 (translated from the Japanese).

\bibitem[Yu(1999)Yu]{yu}
C.C. Yu.
\newblock Autotuning of PID Controllers.
\newblock Springer, 1999.

\bibitem[Zhang, Jiang \& Wang(2002)Zhang, Jiang and Wang]{zhang}
Y. Zhang, D. Jiang, J. Wang.
\newblock A recurrent neural network for solving Sylvester equation with
time-varying coefficients.
\newblock \emph{IEEE Trans. Neural Networks}, 13:\penalty0 1053--1063, 2002.

\bibitem[Ziegler \& Nichols(1942)Ziegler and Nichols]{ziegler}
J.G. Ziegler, N.B. Nichols.
\newblock Optimum settings for automatic controllers.
\newblock \emph{Trans. ASME}, 64:\penalty0 759--768, 1942.

\end{thebibliography}

\begin{figure*}[H]
\centering {\subfigure[\footnotesize i-PI control]{
\rotatebox{-90}{\resizebox{!}{5cm}{%
   \includegraphics{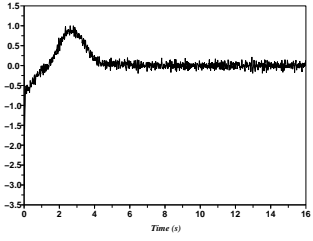}}}}}
 %
    %
\centering {\subfigure[\footnotesize Output (--); reference (- -);
denoised output (. .)]{
\rotatebox{-90}{\resizebox{!}{5cm}{%
   \includegraphics{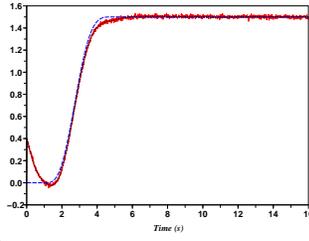}}}}}
   \centering {\subfigure[\footnotesize PID control]{
\rotatebox{-90}{\resizebox{!}{5cm}{%
   \includegraphics{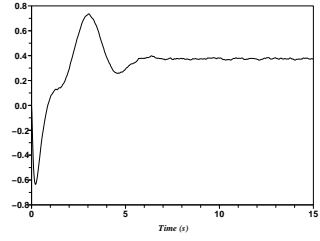}}}}}
 %
    %
\centering {\subfigure[\footnotesize Output (--); reference (- -);
denoised output (. .)]{
\rotatebox{-90}{\resizebox{!}{5cm}{%
   \includegraphics{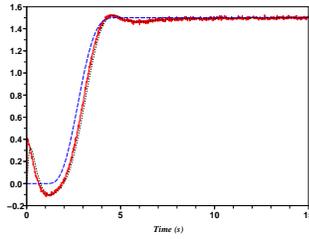}}}}}
\caption{Stable linear monovariable system \label{nom}}
\end{figure*}
\begin{figure*}[H]
\centering {\subfigure[\footnotesize i-PI control]{
\rotatebox{-90}{\resizebox{!}{5cm}{%
   \includegraphics{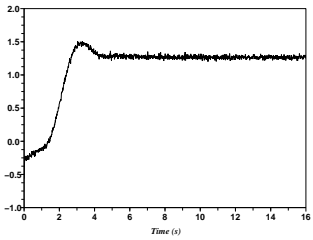}}}}}
 %
    %
\centering {\subfigure[\footnotesize Output (--); reference (- -);
denoised output (. .)]{
\rotatebox{-90}{\resizebox{!}{5cm}{%
   \includegraphics{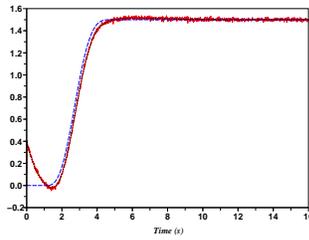}}}}}
   \centering {\subfigure[\footnotesize PID control]{
\rotatebox{-90}{\resizebox{!}{5cm}{%
   \includegraphics{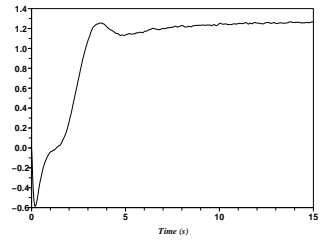}}}}}
 %
    %
\centering {\subfigure[\footnotesize Output (--); reference (- -);
denoised output (. .)]{
\rotatebox{-90}{\resizebox{!}{5cm}{%
   \includegraphics{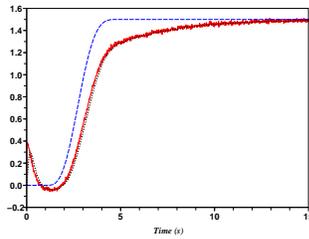}}}}}
\caption{Modified stable linear monovariable system \label{hold}}
\end{figure*}

\begin{figure*}[H]
\centering {\subfigure[\footnotesize i-PI control]{
\rotatebox{-90}{\resizebox{!}{5cm}{%
   \includegraphics{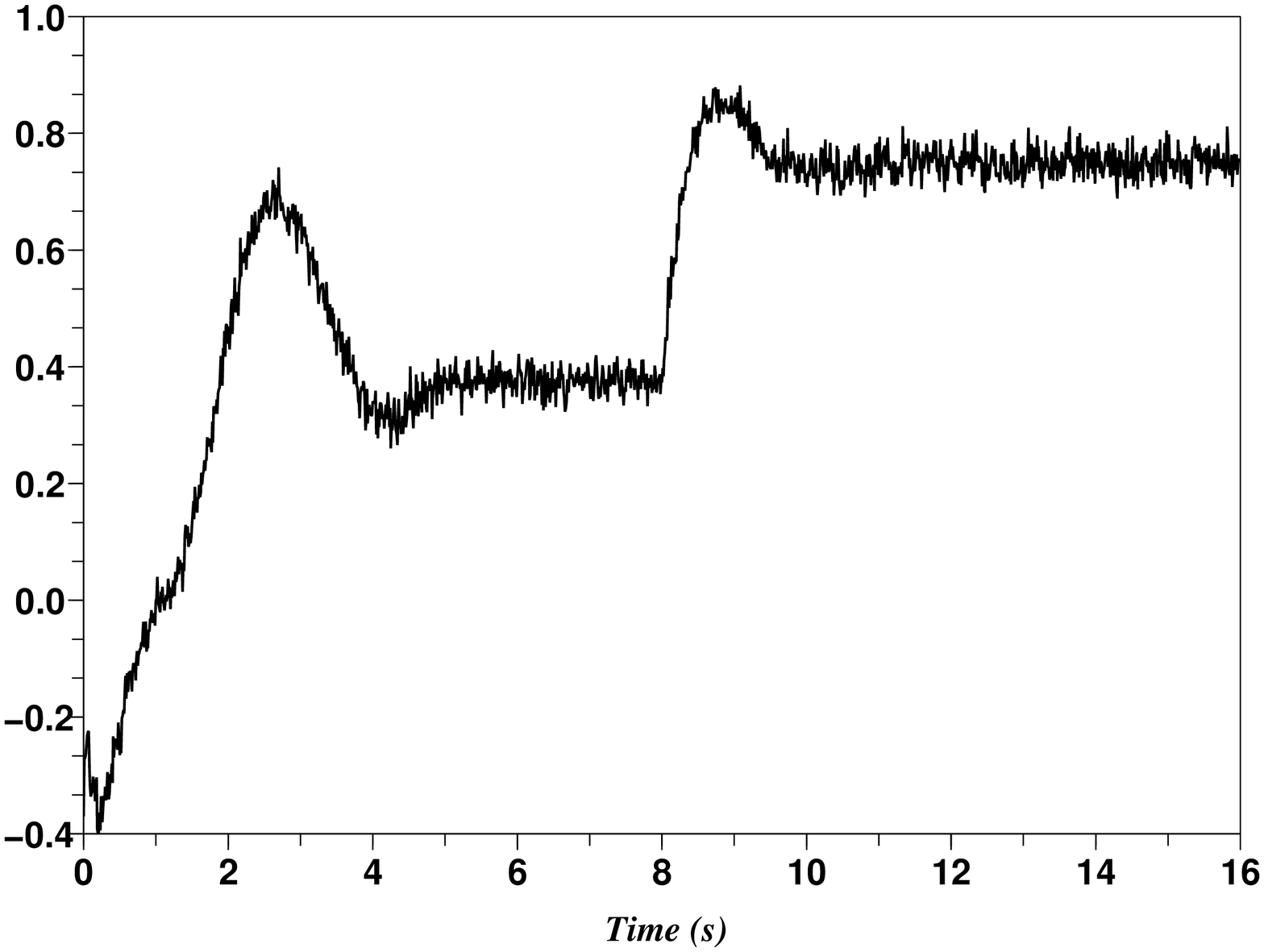}}}}}
 %
    %
\centering {\subfigure[\footnotesize Output (--); reference (- -);
denoised output (. .)]{
\rotatebox{-90}{\resizebox{!}{5cm}{%
   \includegraphics{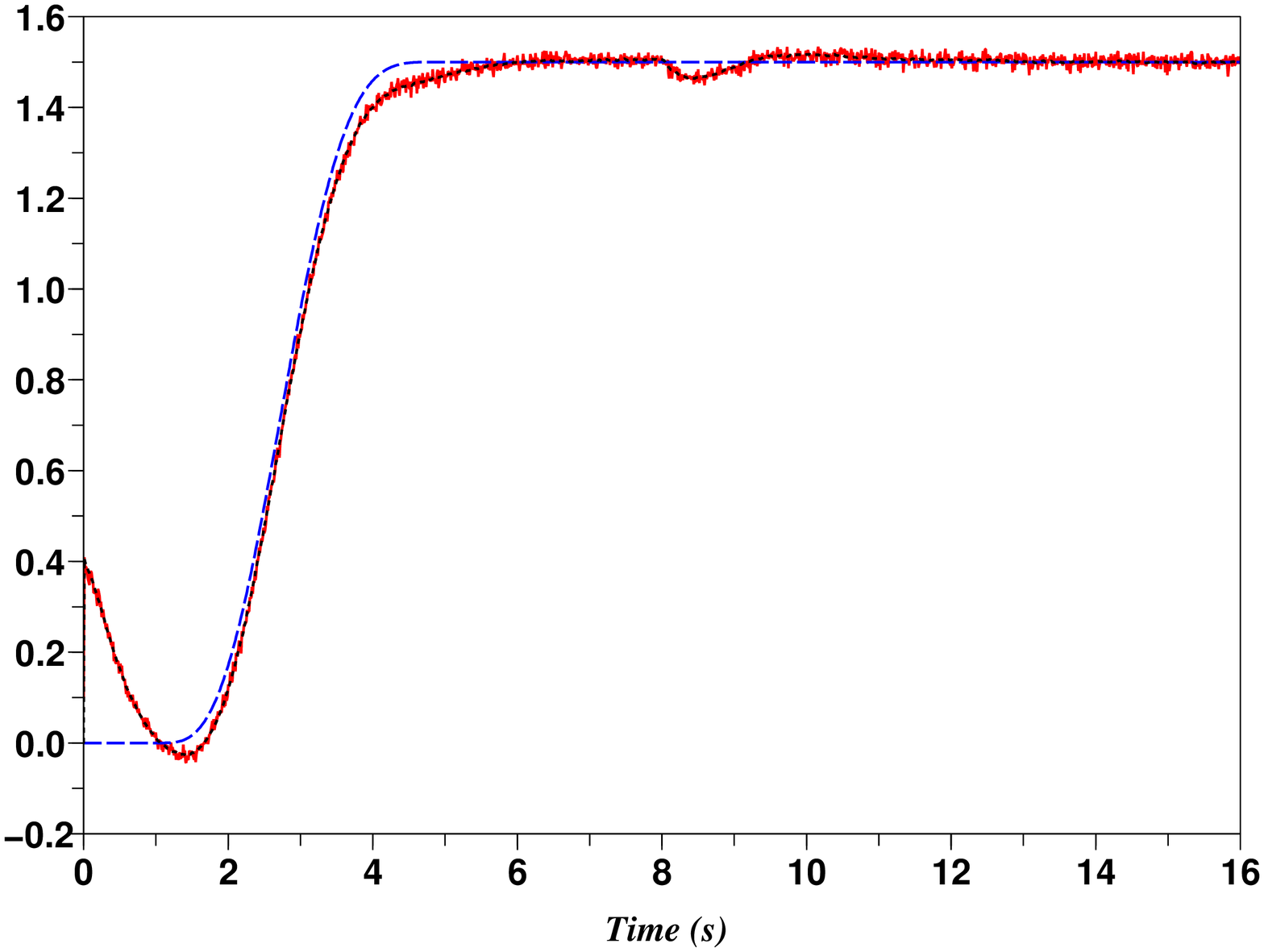}}}}}
   \centering {\subfigure[\footnotesize PID control]{
\rotatebox{-90}{\resizebox{!}{5cm}{%
   \includegraphics{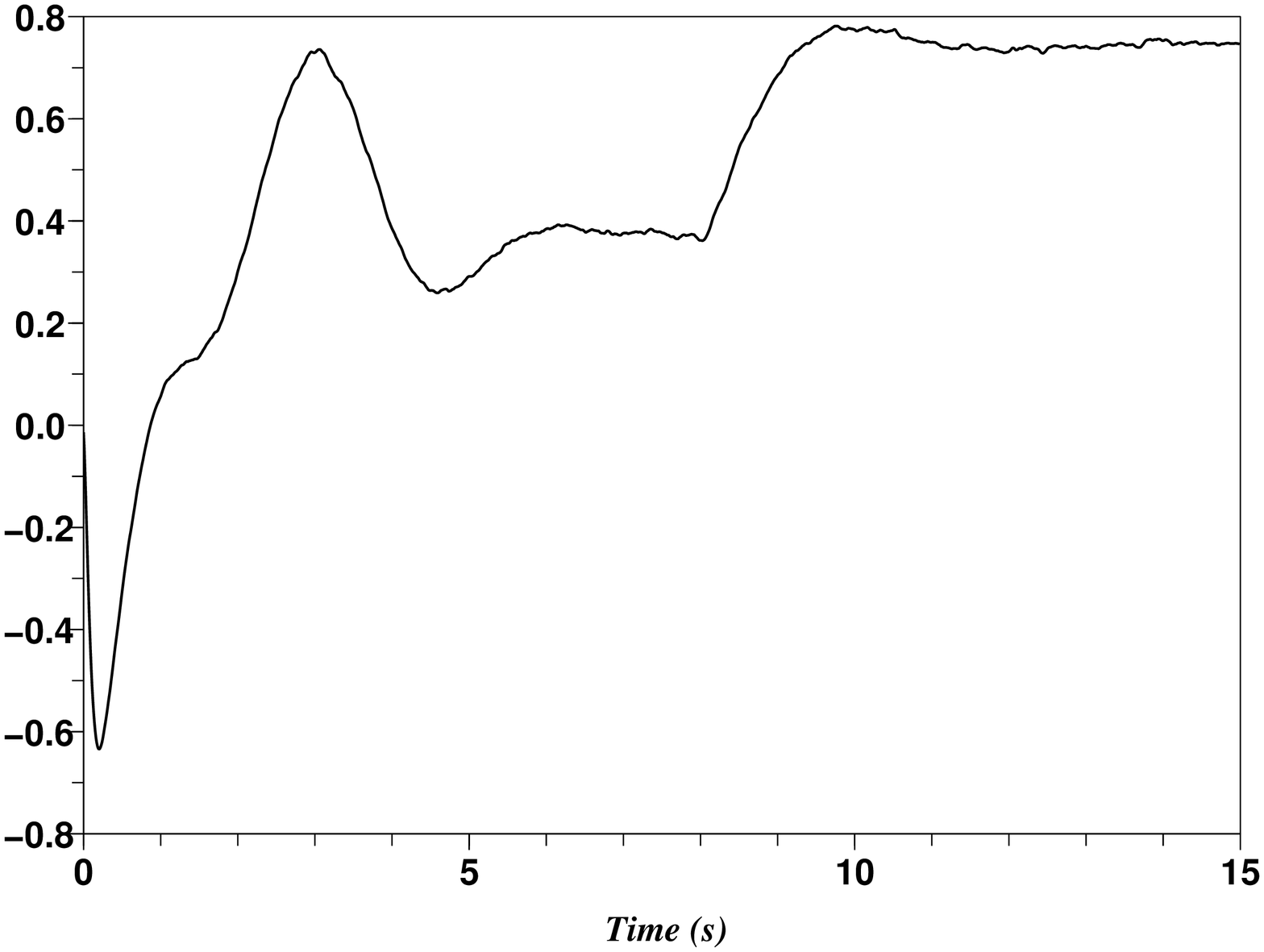}}}}}
 %
    %
\centering {\subfigure[\footnotesize Output (--); reference (- -);
denoised output (. .)]{
\rotatebox{-90}{\resizebox{!}{5cm}{%
   \includegraphics{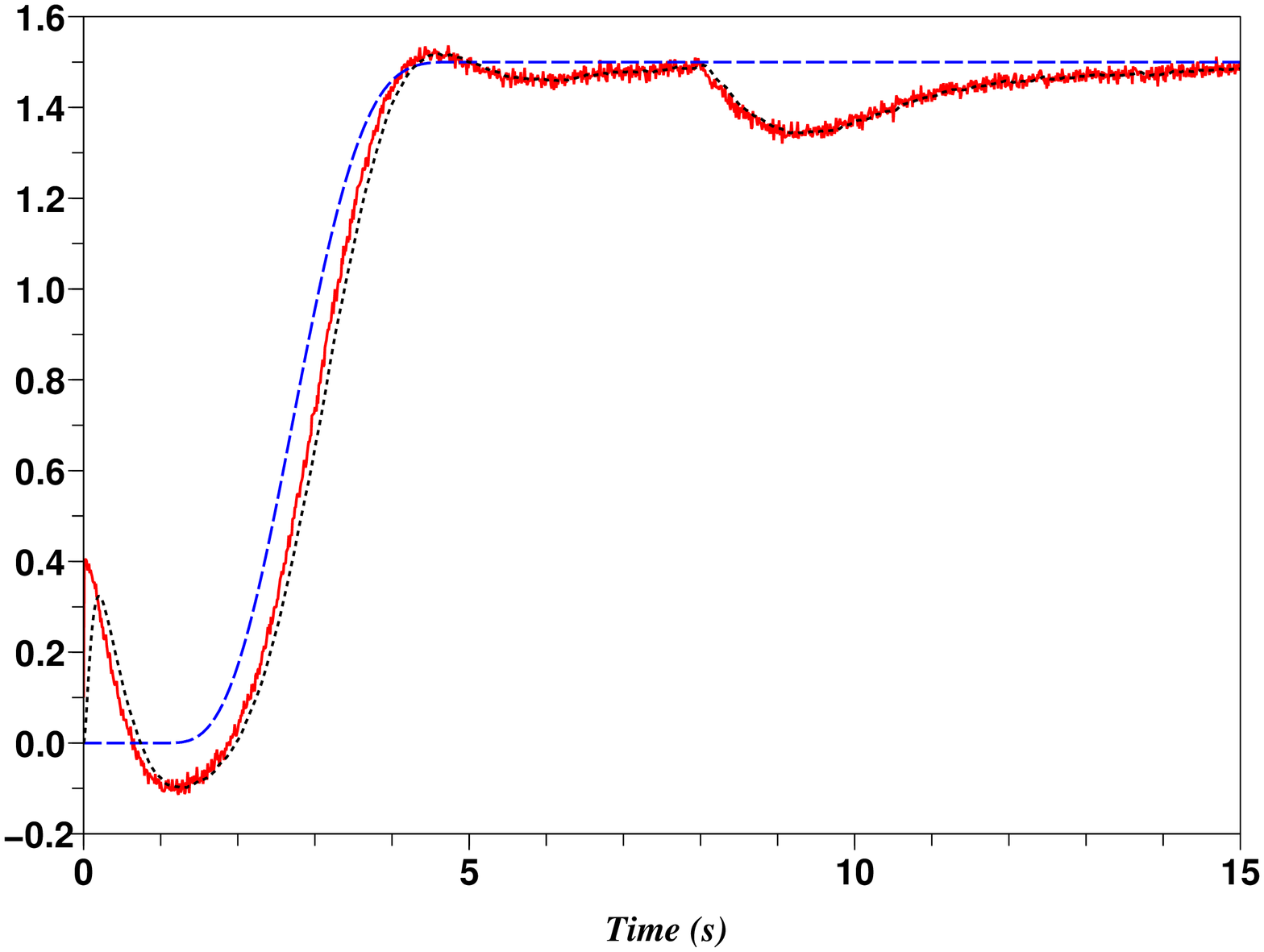}}}}}
\caption{Stable linear monovariable system, with an actuator's fault
\label{def}}
\end{figure*}

\begin{figure*}[H]
\centering%
\subfigure[\footnotesize Reference (- -) and output]
{\epsfig{figure=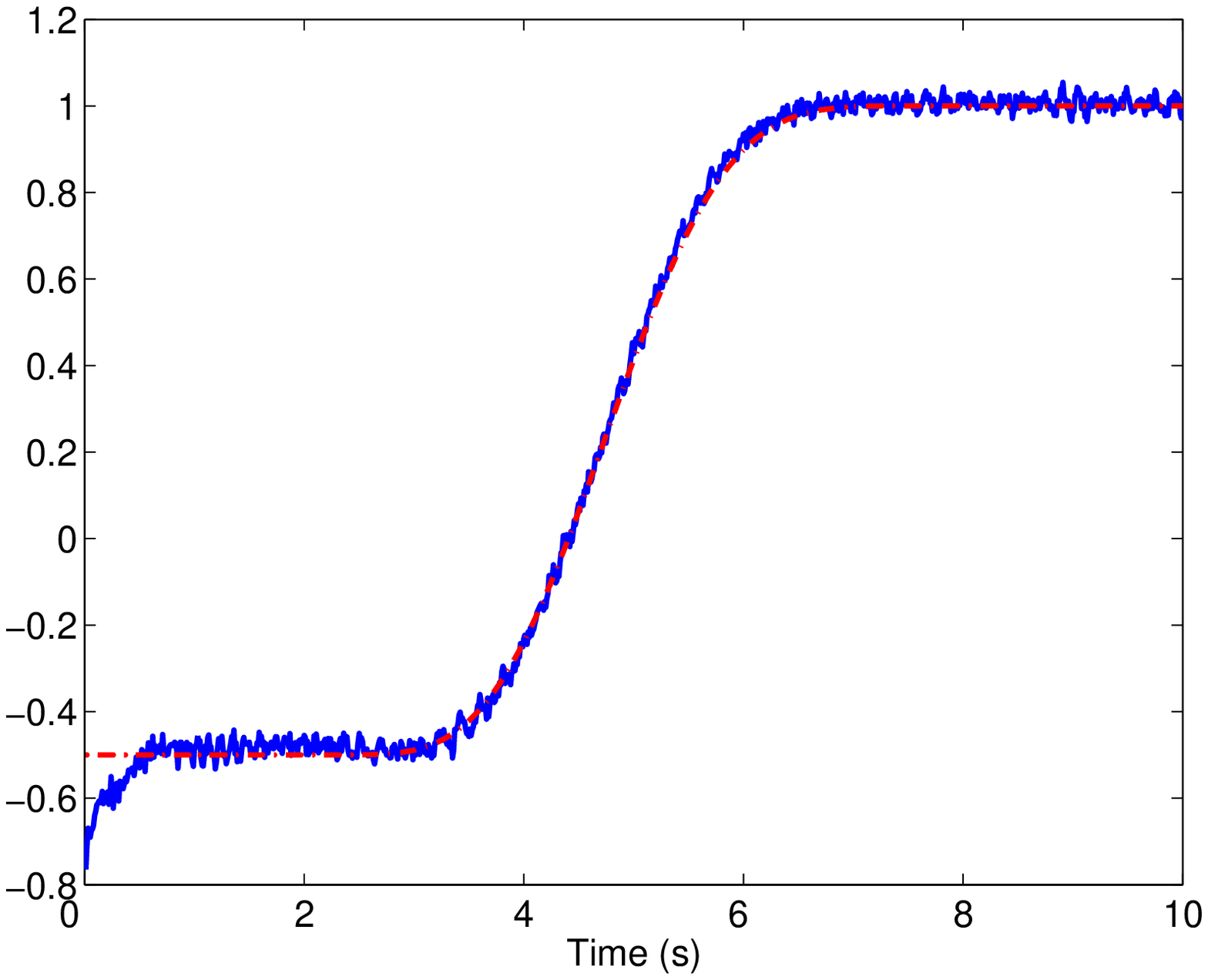,width= 0.4\textwidth}}
\subfigure[\footnotesize Control]
{\epsfig{figure=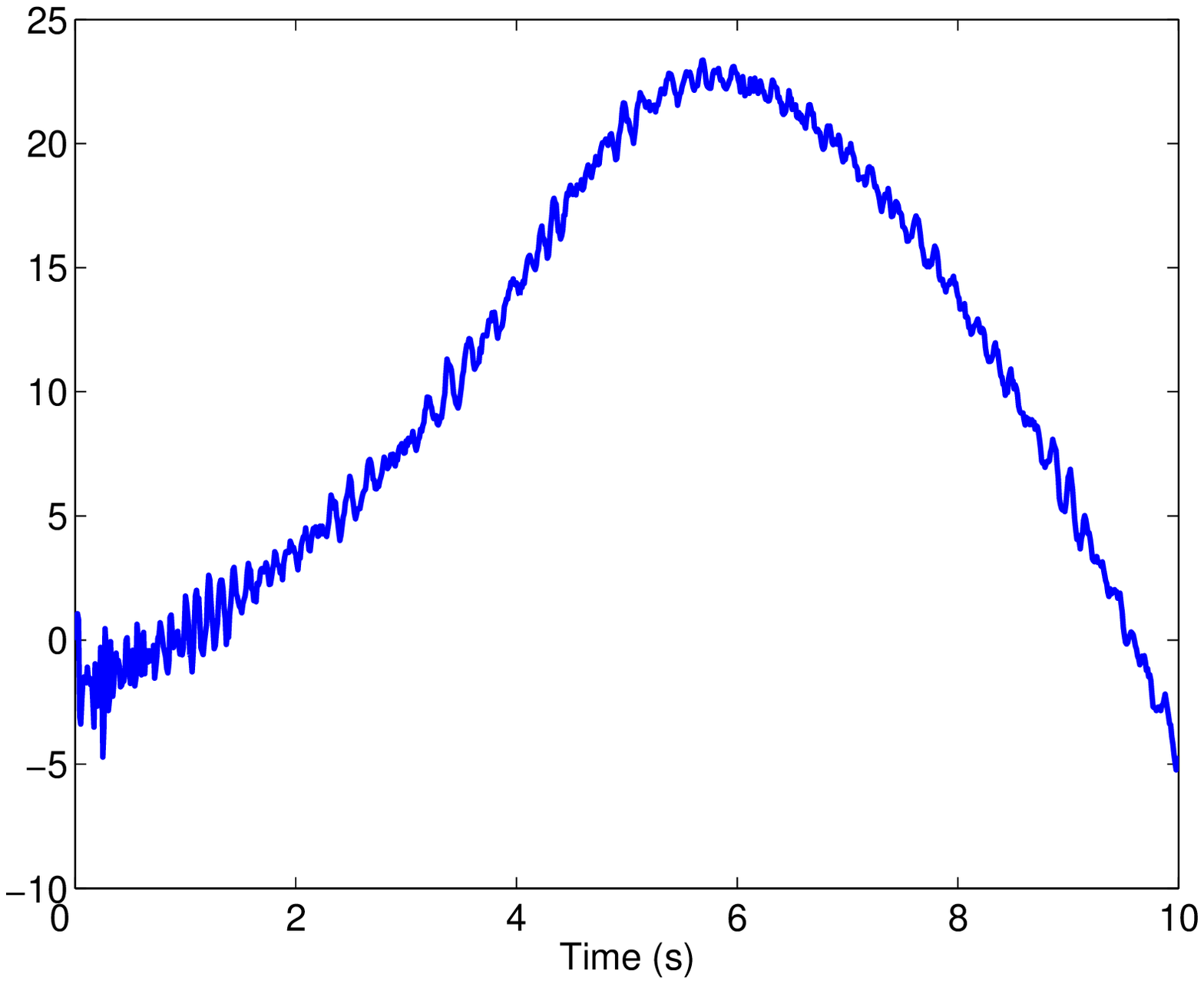,width= 0.4\textwidth}}\\
\subfigure[\footnotesize Estimation of $F$]
{\epsfig{figure=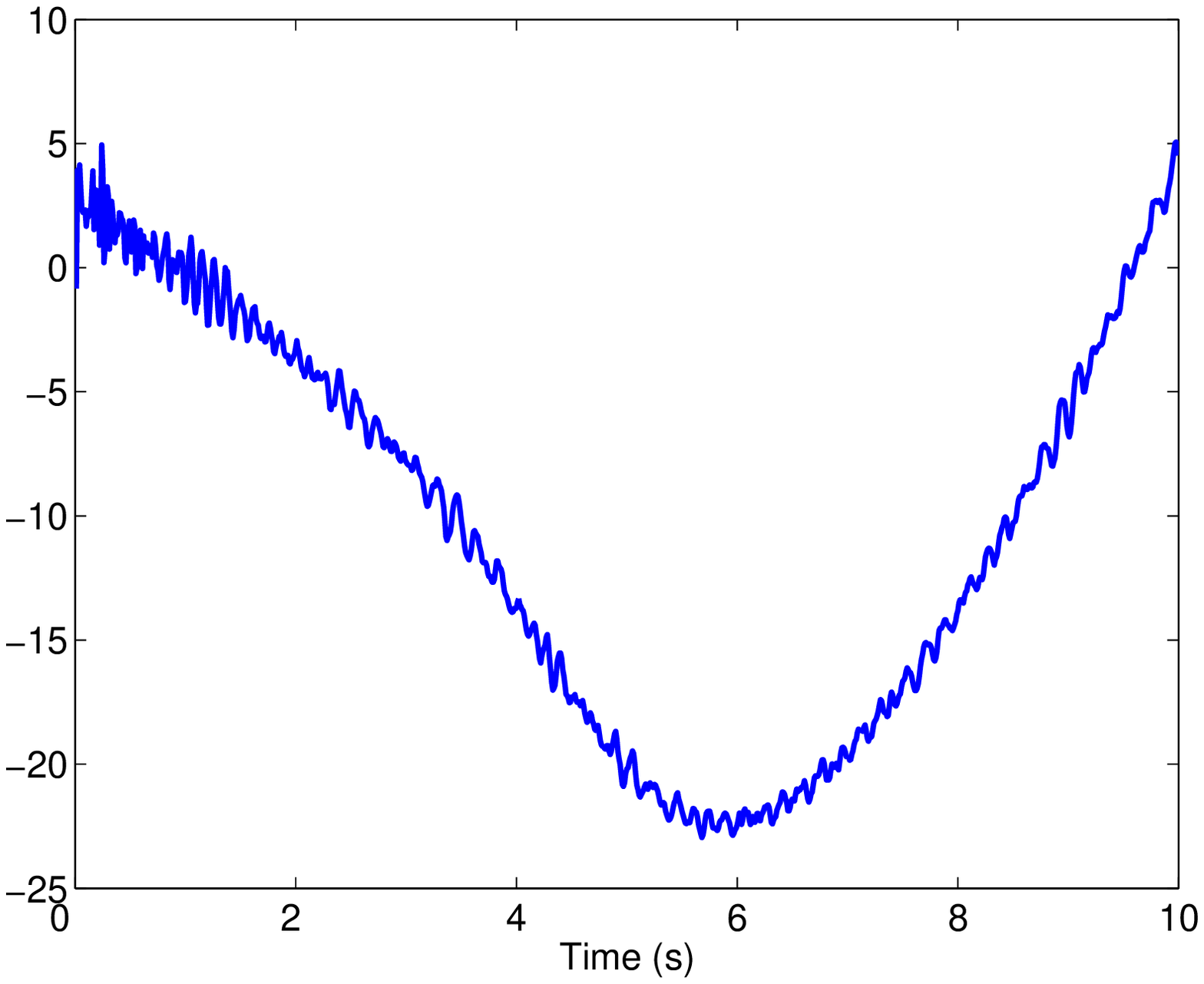,width= 0.4\textwidth}}
\subfigure[\footnotesize Estimation of $\dot y$]
{\epsfig{figure=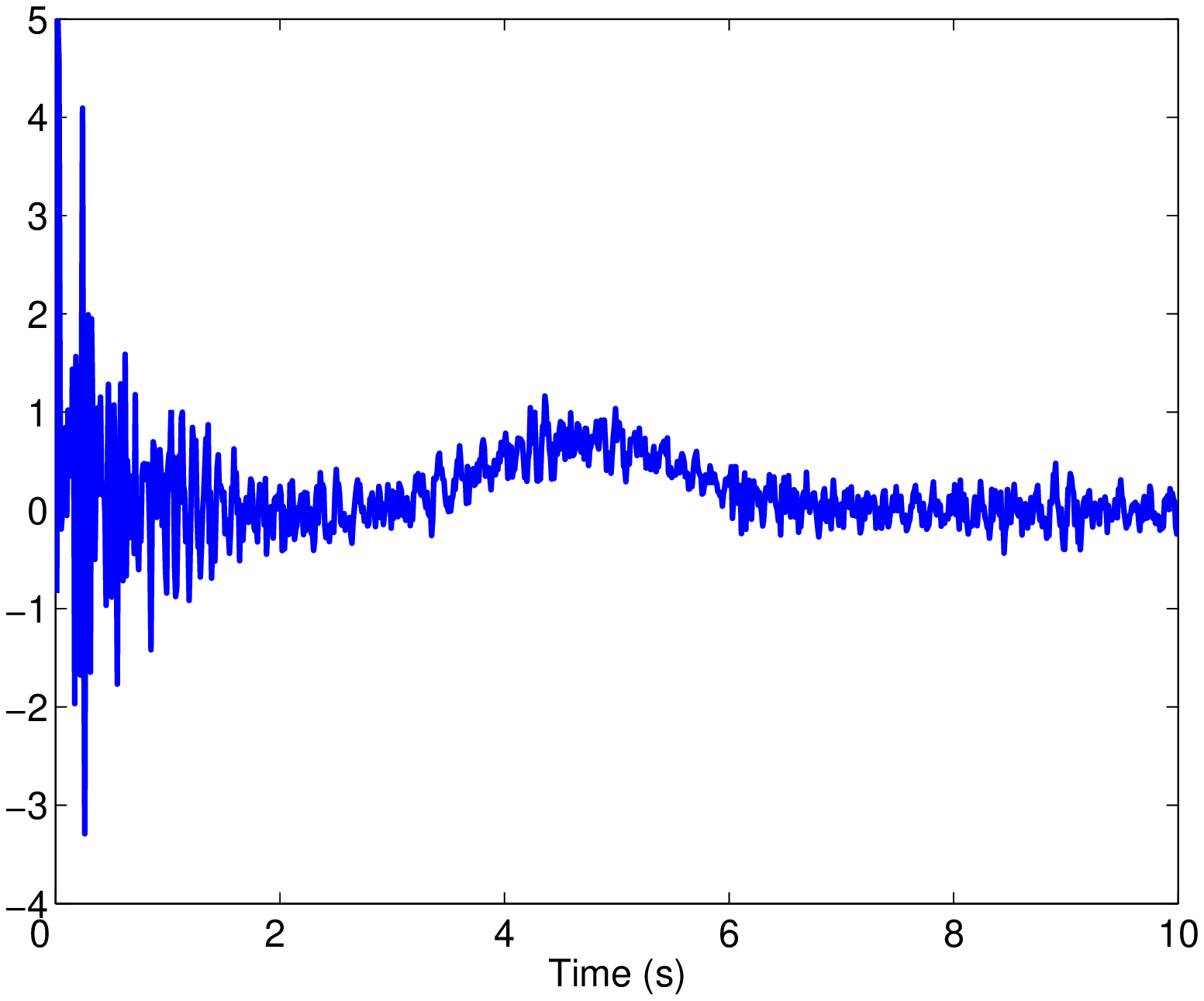,width= 0.4\textwidth}}
\caption{Linear monovariable system with a large spectrum}\label{lx}
\end{figure*}
\begin{figure*}[h!t]
\centering
%
\vspace{-0.5cm} 
{\subfigure[\footnotesize Output's
noises]{\rotatebox{-90}{\resizebox{!}{5cm}{\includegraphics{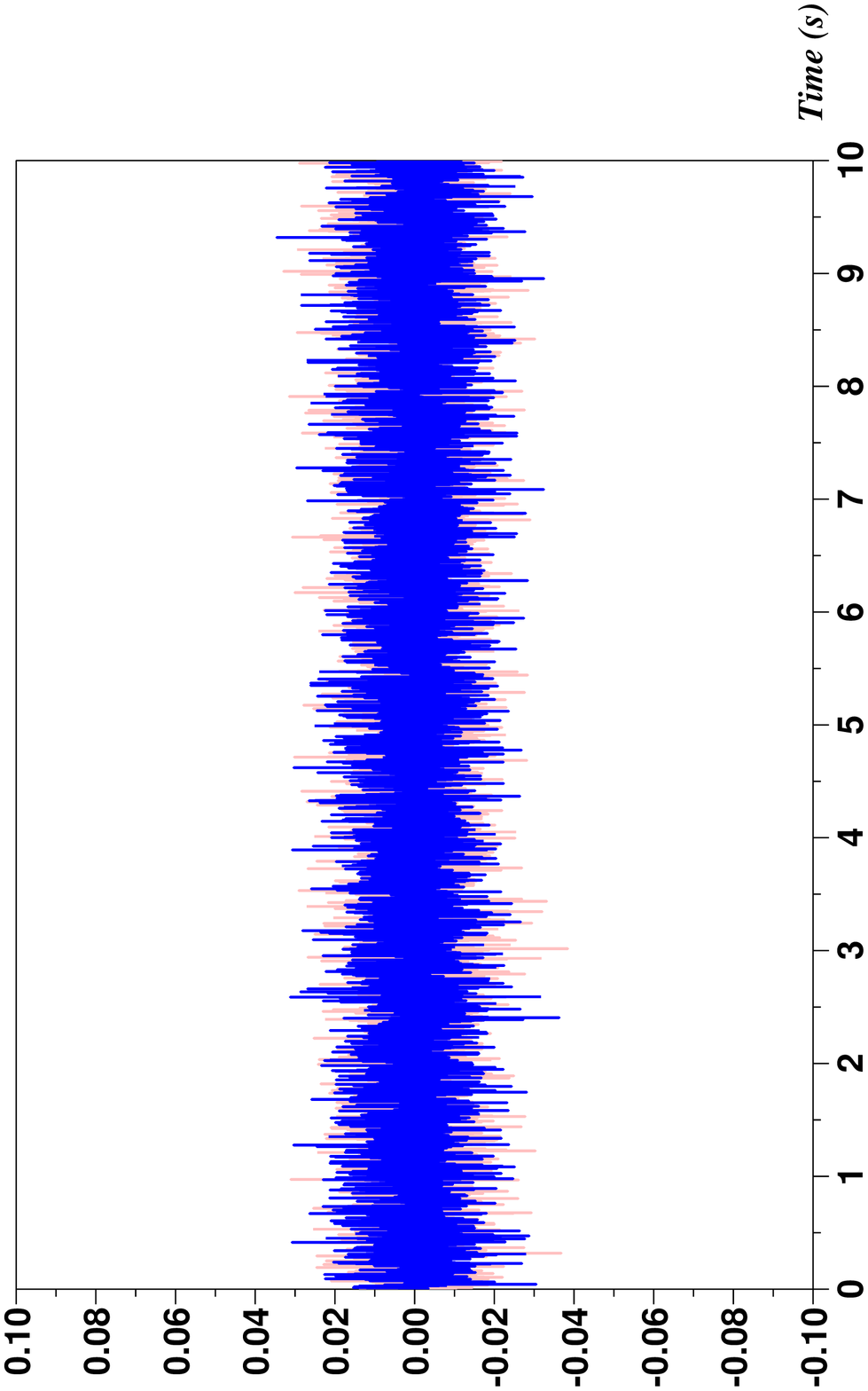}}}}}
%
 {\subfigure[\footnotesize Estimation of $\dot y_1$]{
\rotatebox{-90}{\resizebox{!}{5cm}{\includegraphics{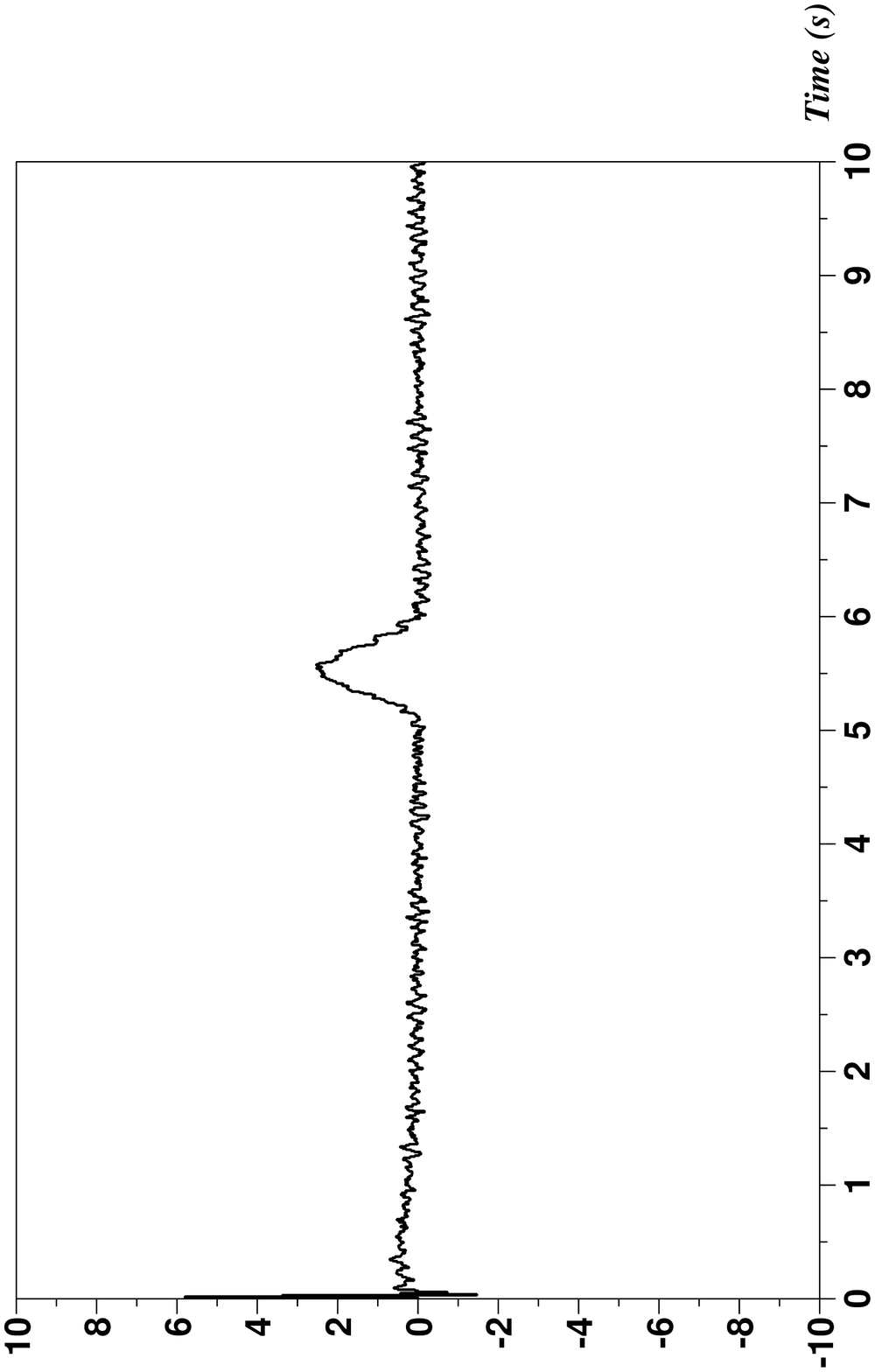}}}}}
\vspace{-0.3cm} {\subfigure[\footnotesize Estimation of $\dot y_2$]{
\rotatebox{-90}{\resizebox{!}{5cm}{\includegraphics{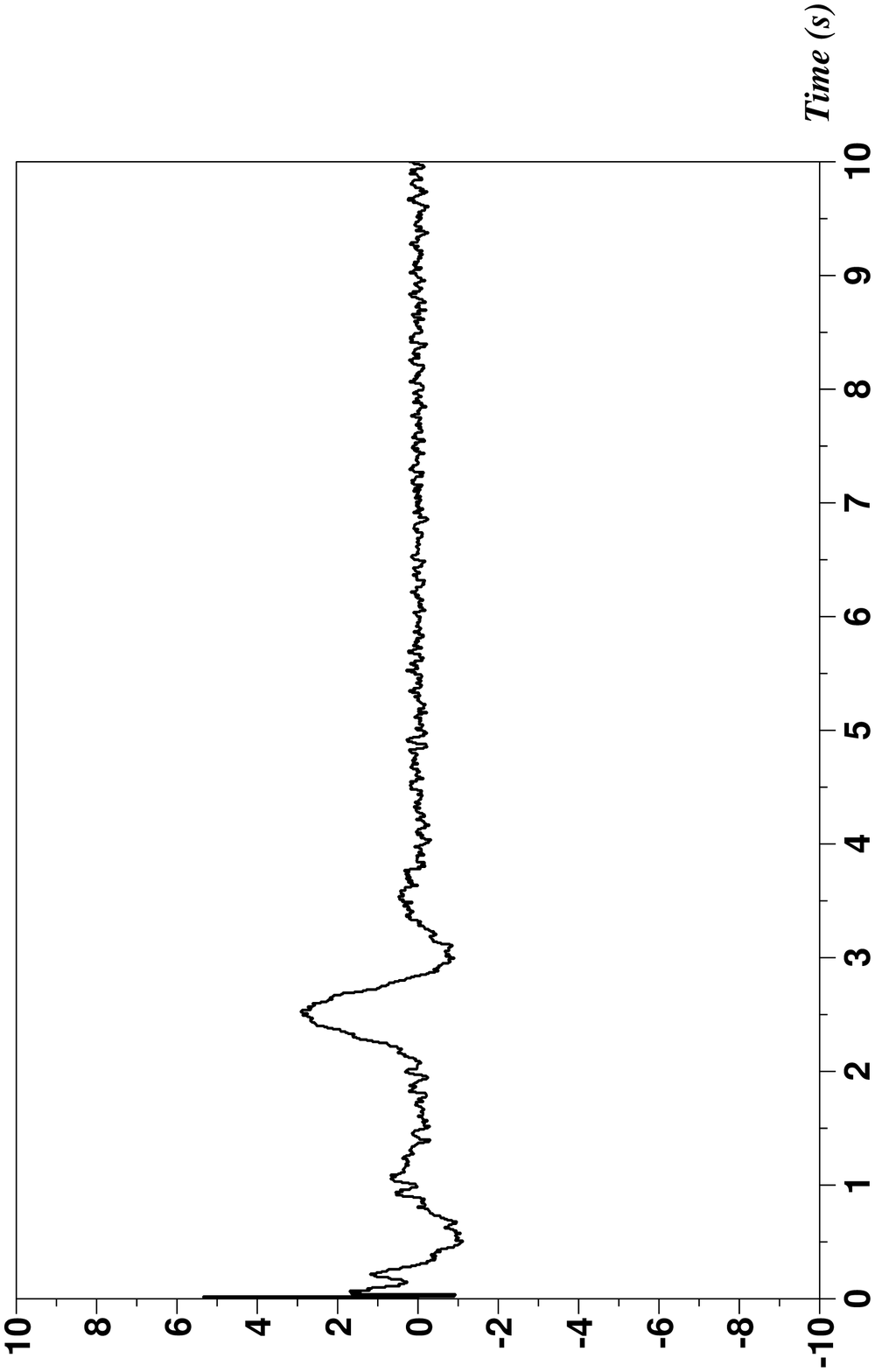}}}}}
%
{\subfigure[\footnotesize Estimation of $\ddot y_2$]{
\rotatebox{-90}{\resizebox{!}{5cm}{\includegraphics{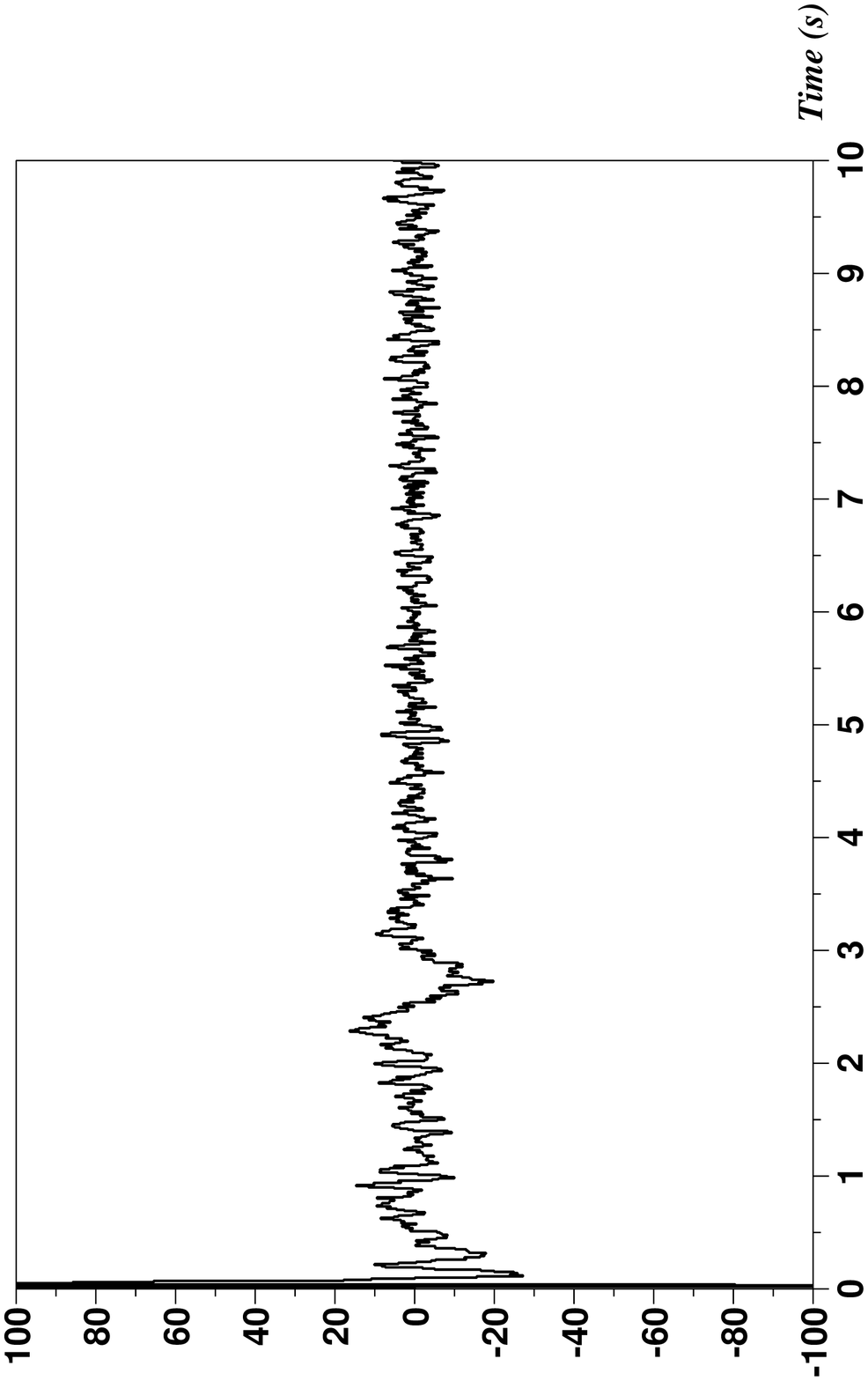}}}}}
 \vspace{-0.3cm}
{\subfigure[\footnotesize Estimation of $F_1$]{
\rotatebox{-90}{\resizebox{!}{5cm}{\includegraphics{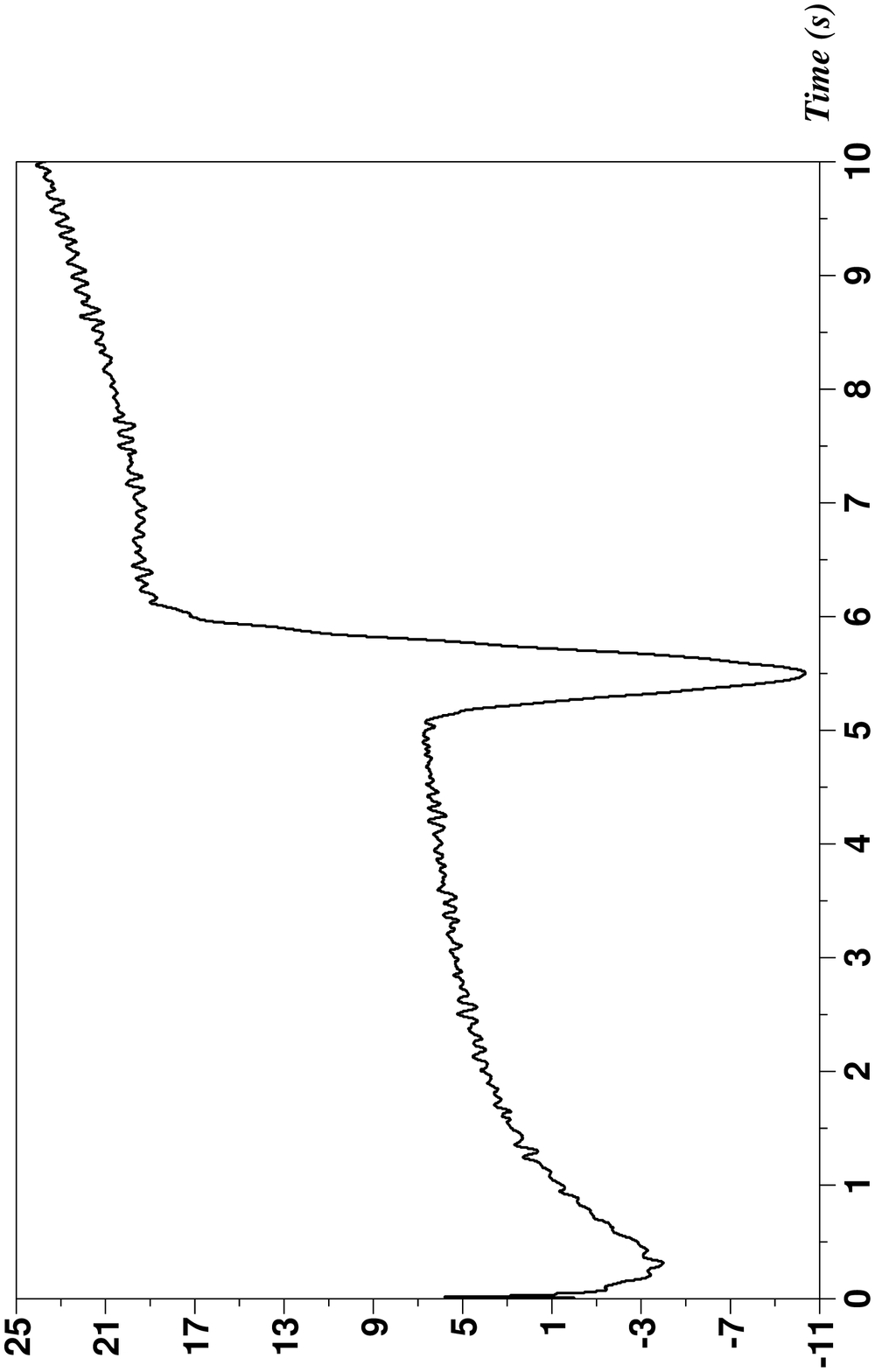}}}}}
{\subfigure[\footnotesize Estimation of $F_2$]{
\rotatebox{-90}{\resizebox{!}{5cm}{\includegraphics{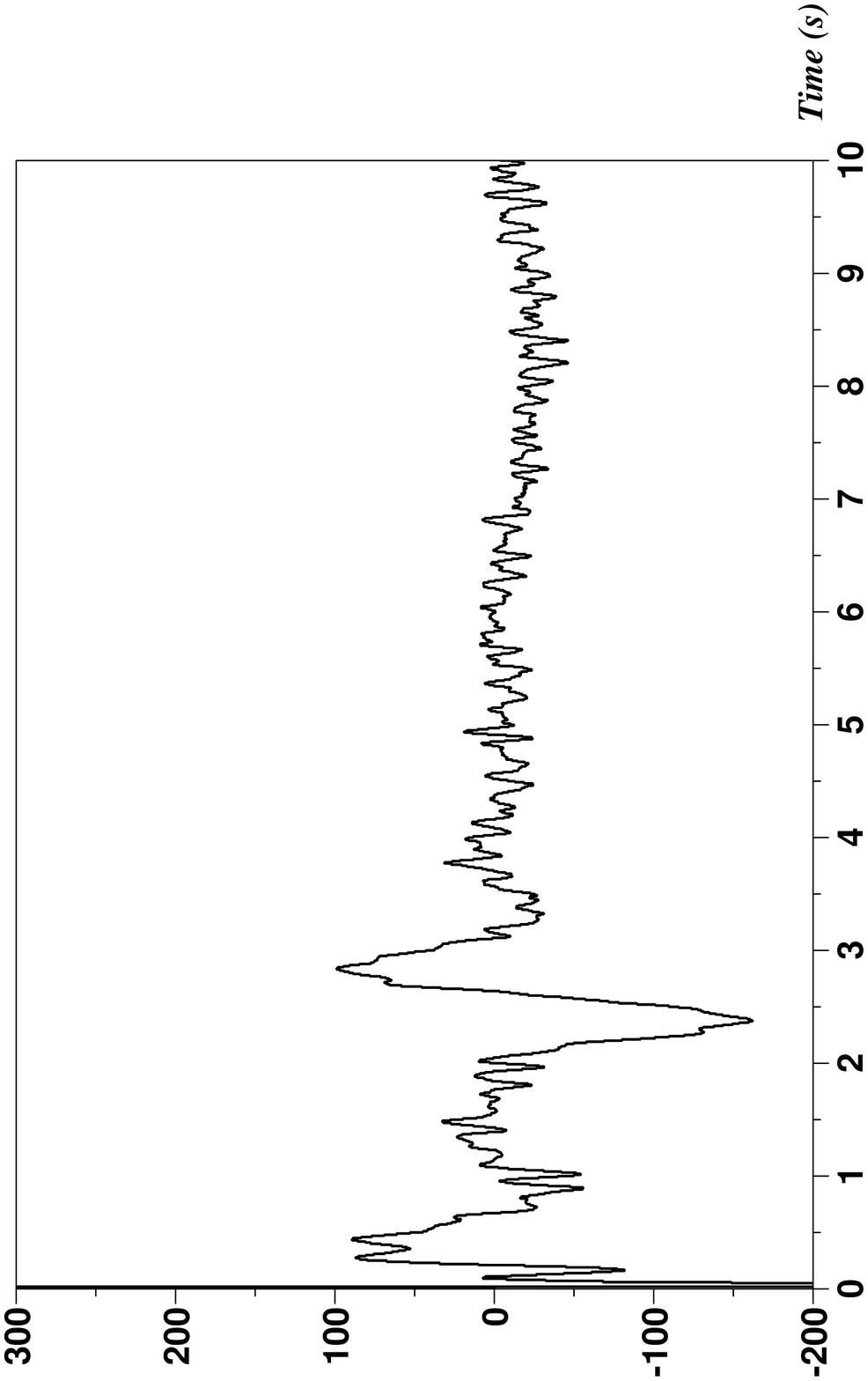}}}}}
 \vspace{-0.3cm}
{\subfigure[\footnotesize Control $u_1$]{
\rotatebox{-90}{\resizebox{!}{5cm}{\includegraphics{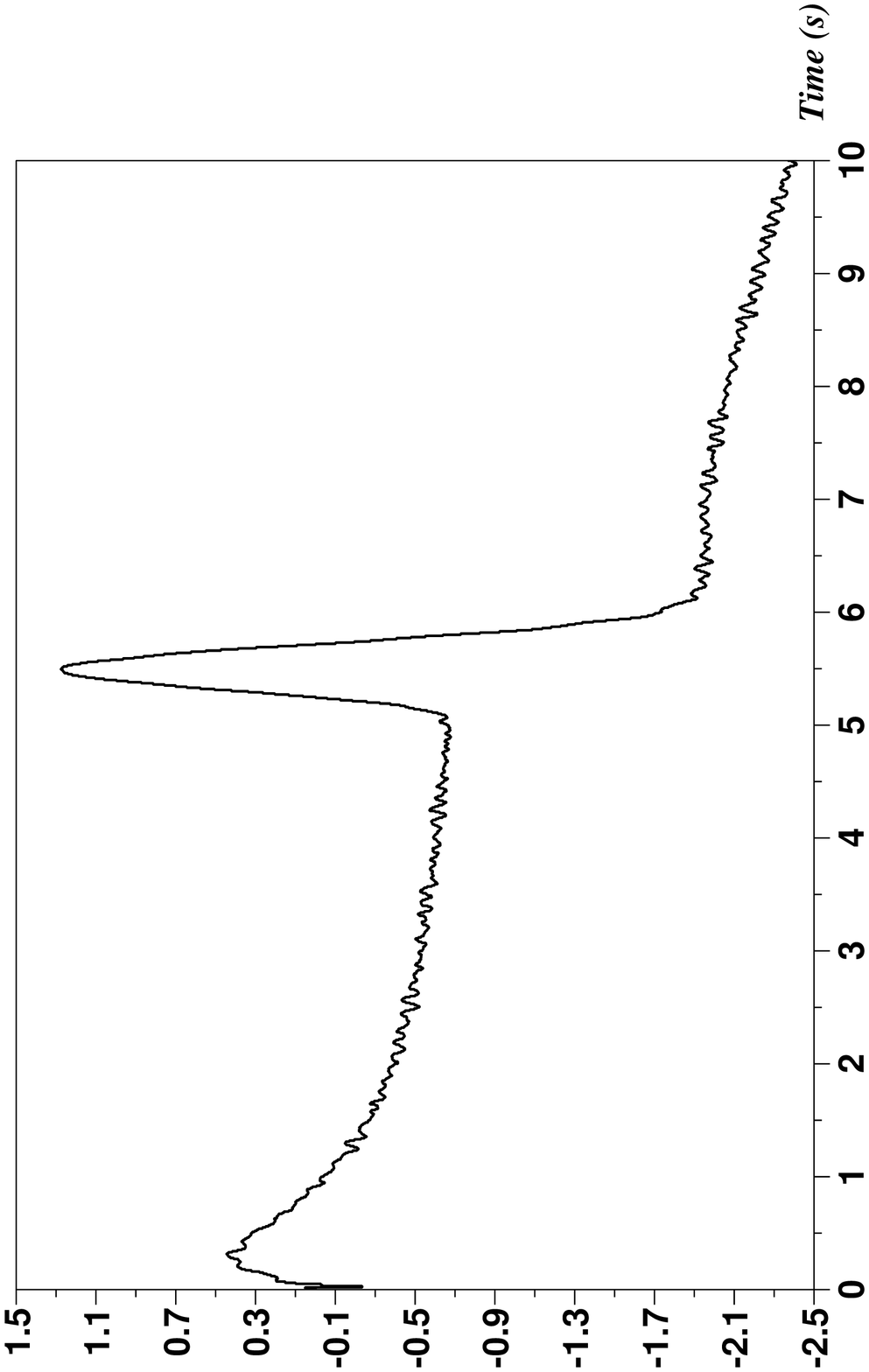}}}}}
{\subfigure[\footnotesize Control $u_2$]{
\rotatebox{-90}{\resizebox{!}{5cm}{\includegraphics{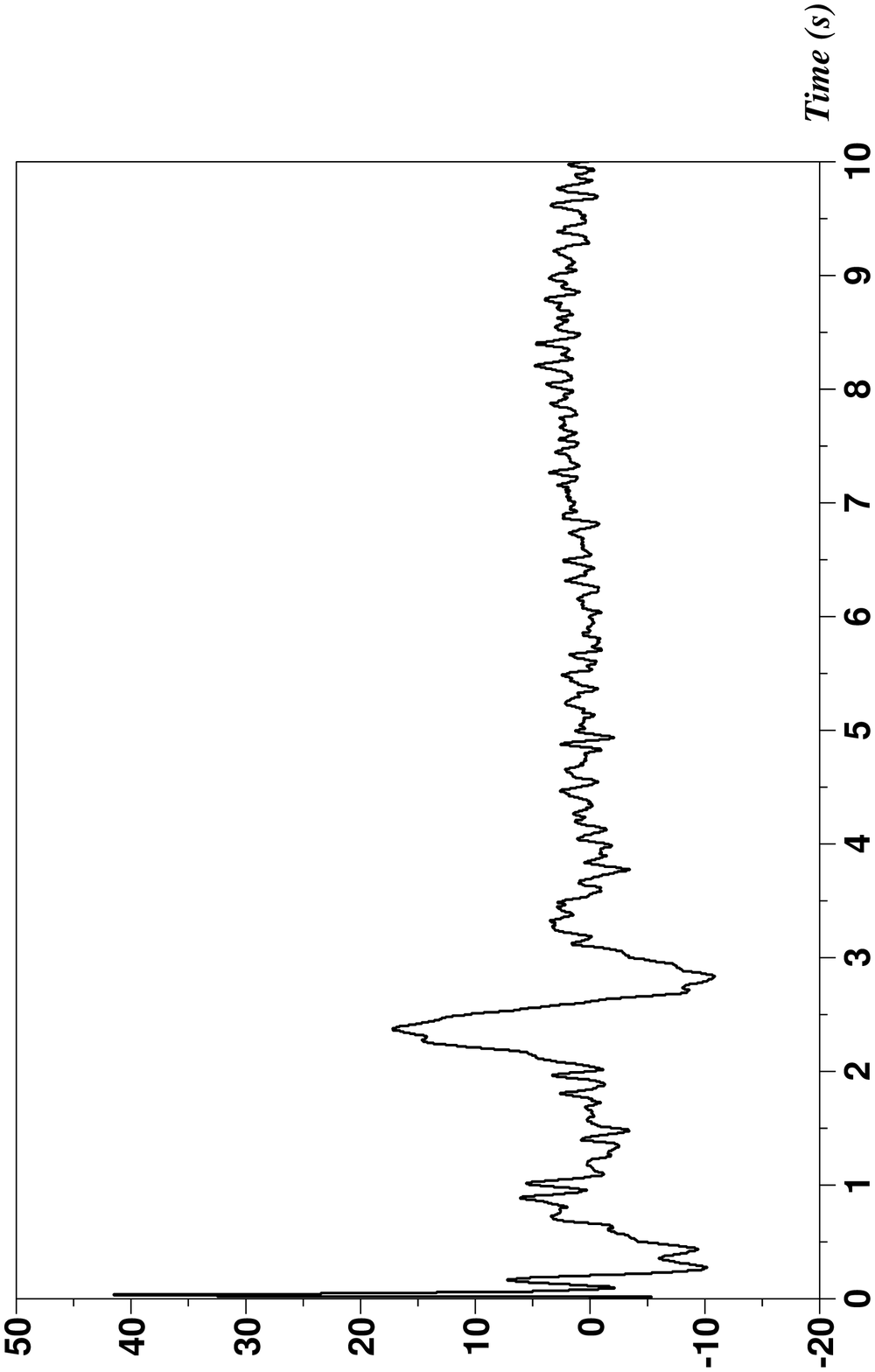}}}}}
 %
%
 \caption{Linear multivariable system \label{fig_Ml1}}
\end{figure*}

\begin{figure*}[h!t]
\centering
 {\subfigure[\footnotesize i-PID: references (- -) and outputs]{
\rotatebox{-90}{\resizebox{!}{5cm}{%
   \includegraphics{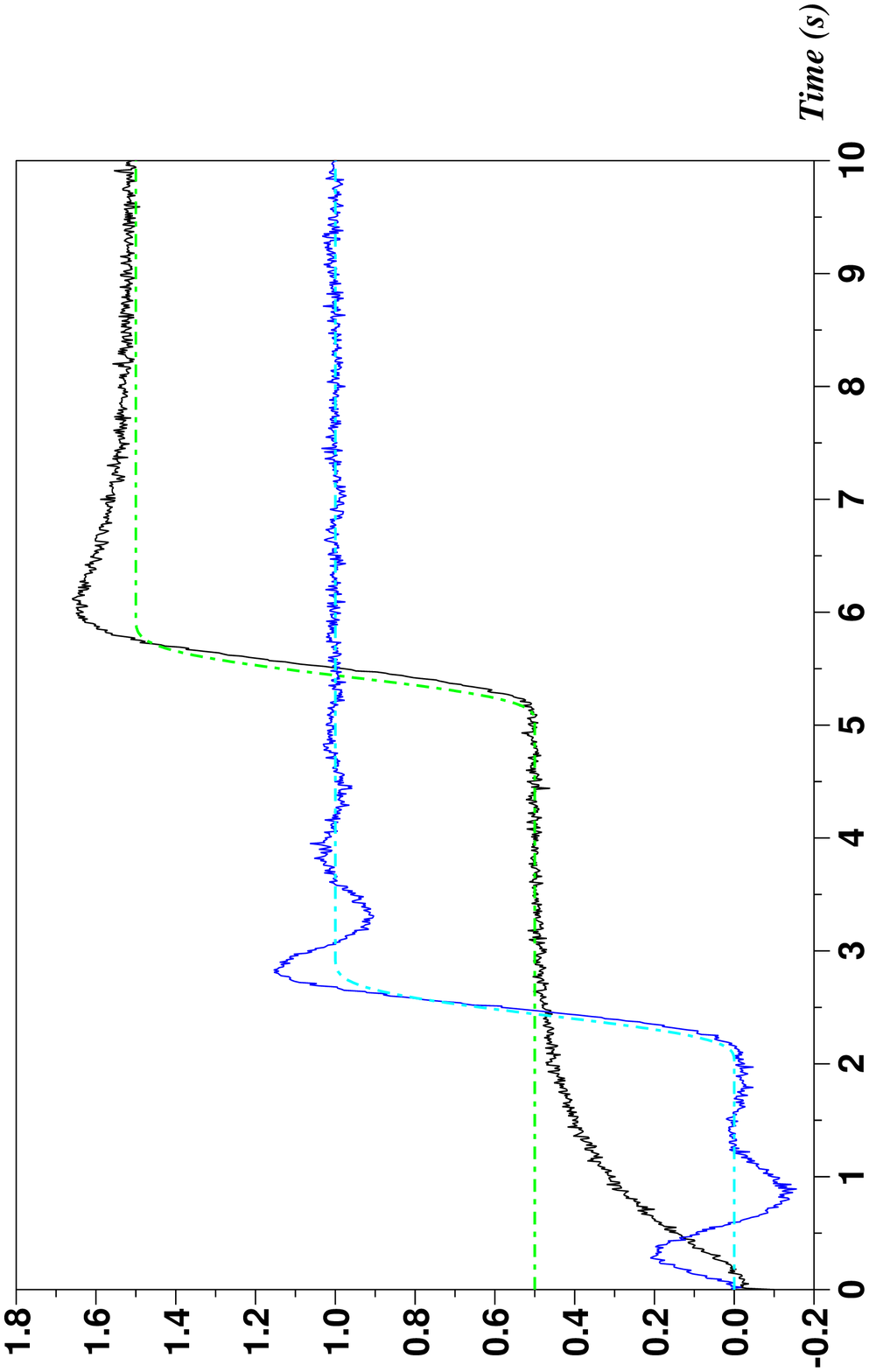}}}}}
    {\subfigure[\footnotesize PID: references (- -) and
outputs]{
\rotatebox{-90}{\resizebox{!}{5cm}{%
   \includegraphics{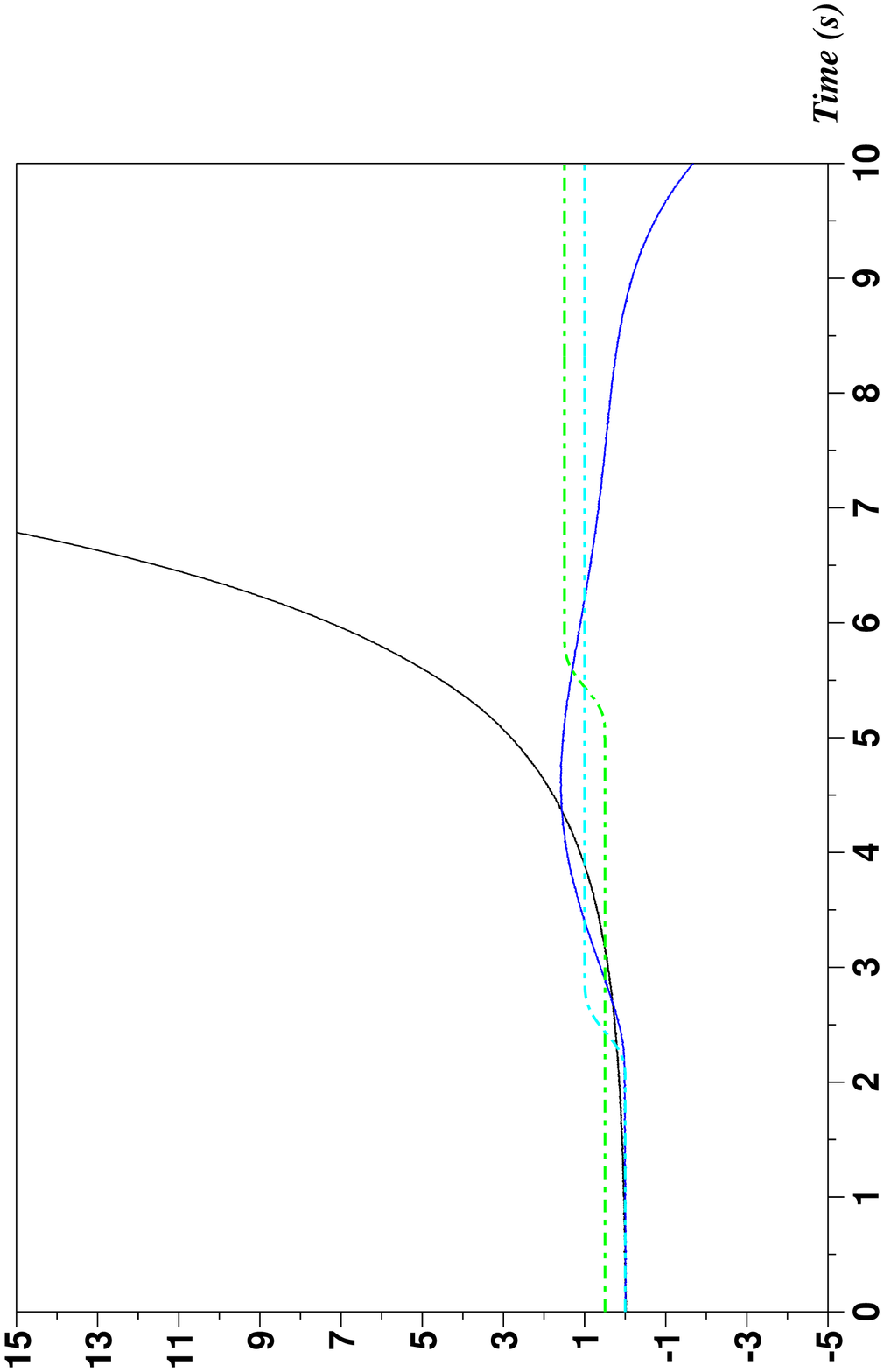}}}}}
 \caption{Linear multivariable system\label{fig_Ml2}}
\end{figure*}

\begin{figure*}[H]
\centering {\subfigure[\footnotesize Estimation of $\dot
y$]{\rotatebox{-90}{\resizebox{!}{5cm}{\includegraphics{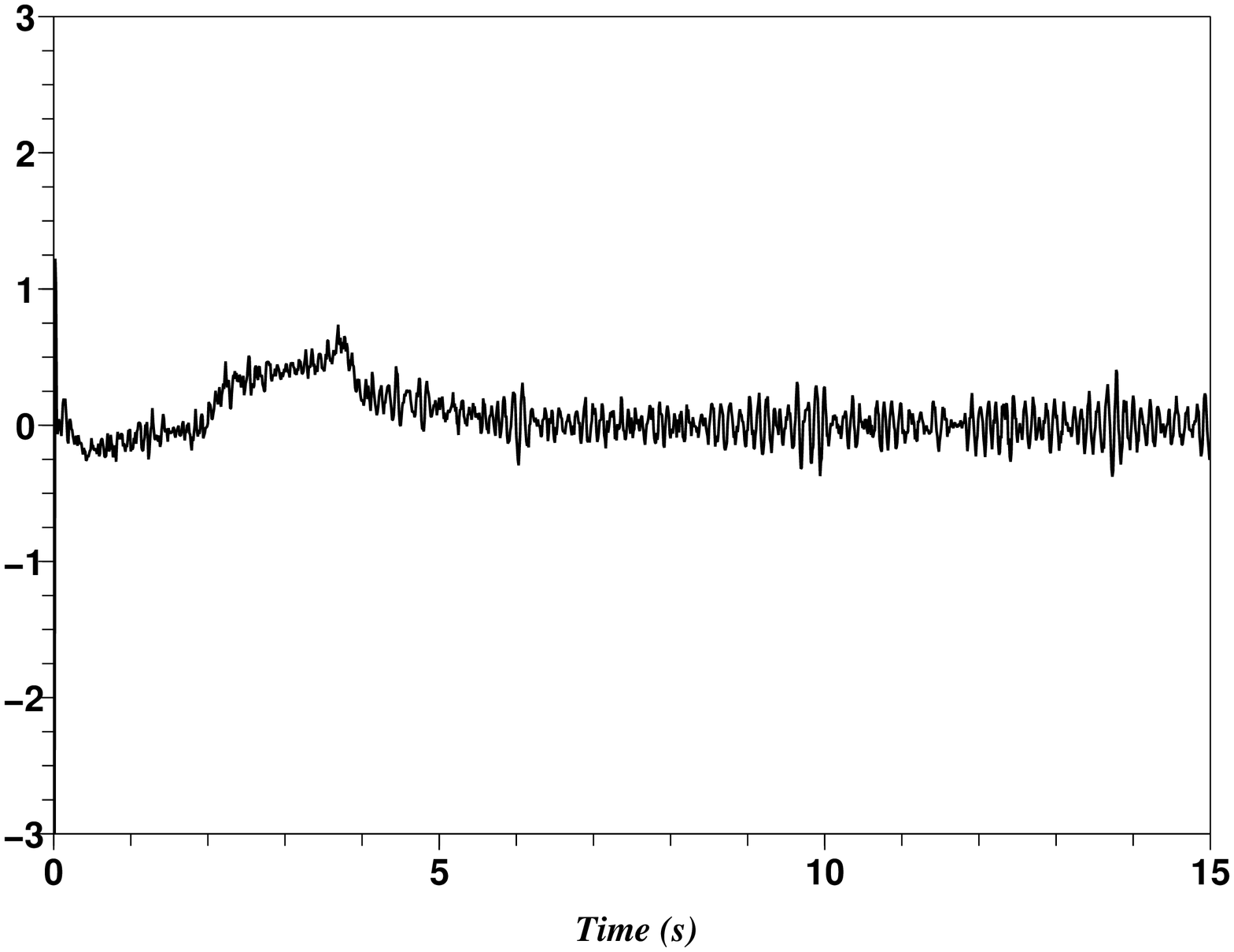}}}}}
{\subfigure[\footnotesize
Commande]{\rotatebox{-90}{\resizebox{!}{5cm}{\includegraphics{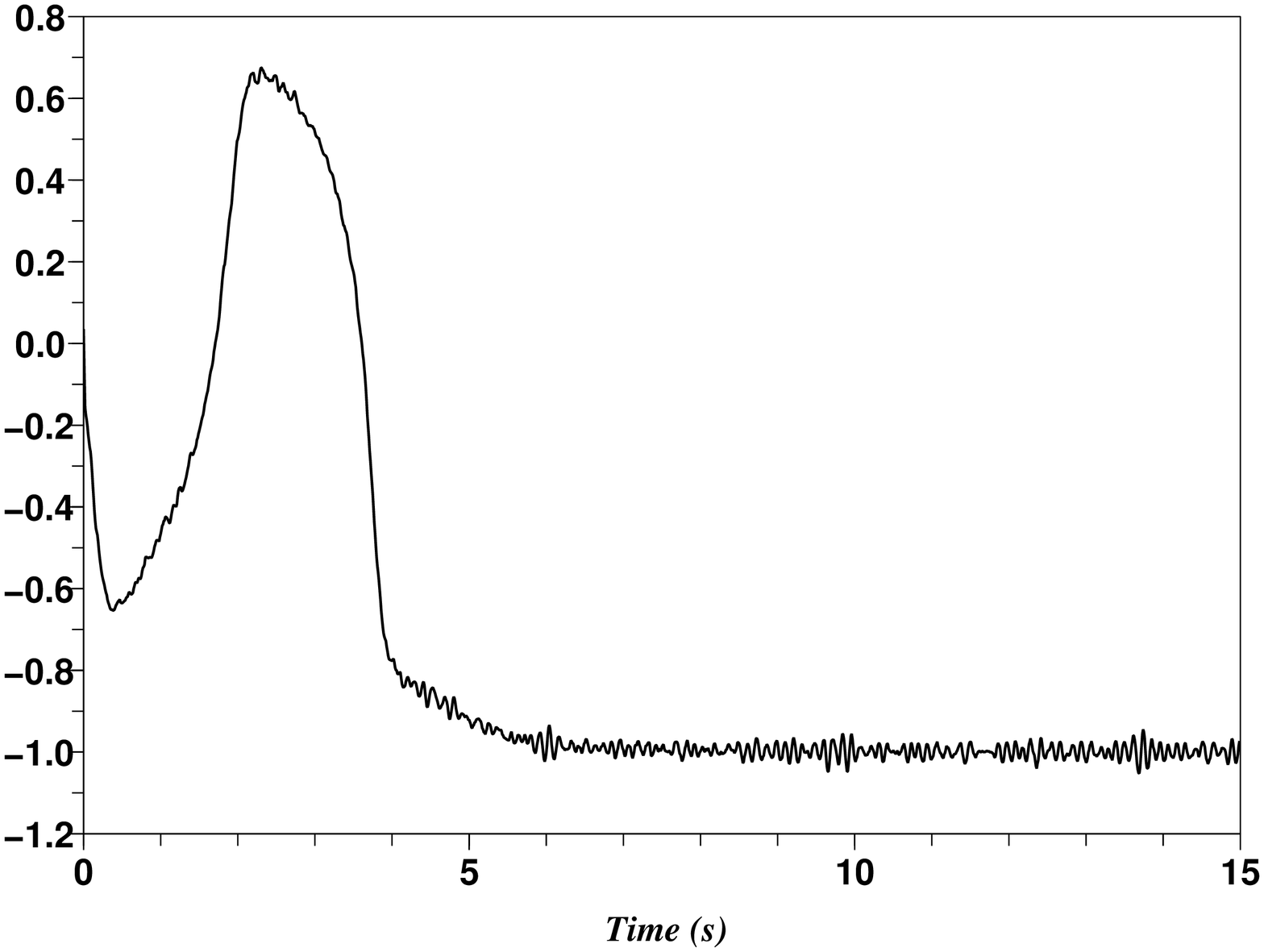}}}}}
{\subfigure[\footnotesize Output (--); reference (- -); denoised
output (.
.)]{\rotatebox{-90}{\resizebox{!}{5cm}{\includegraphics{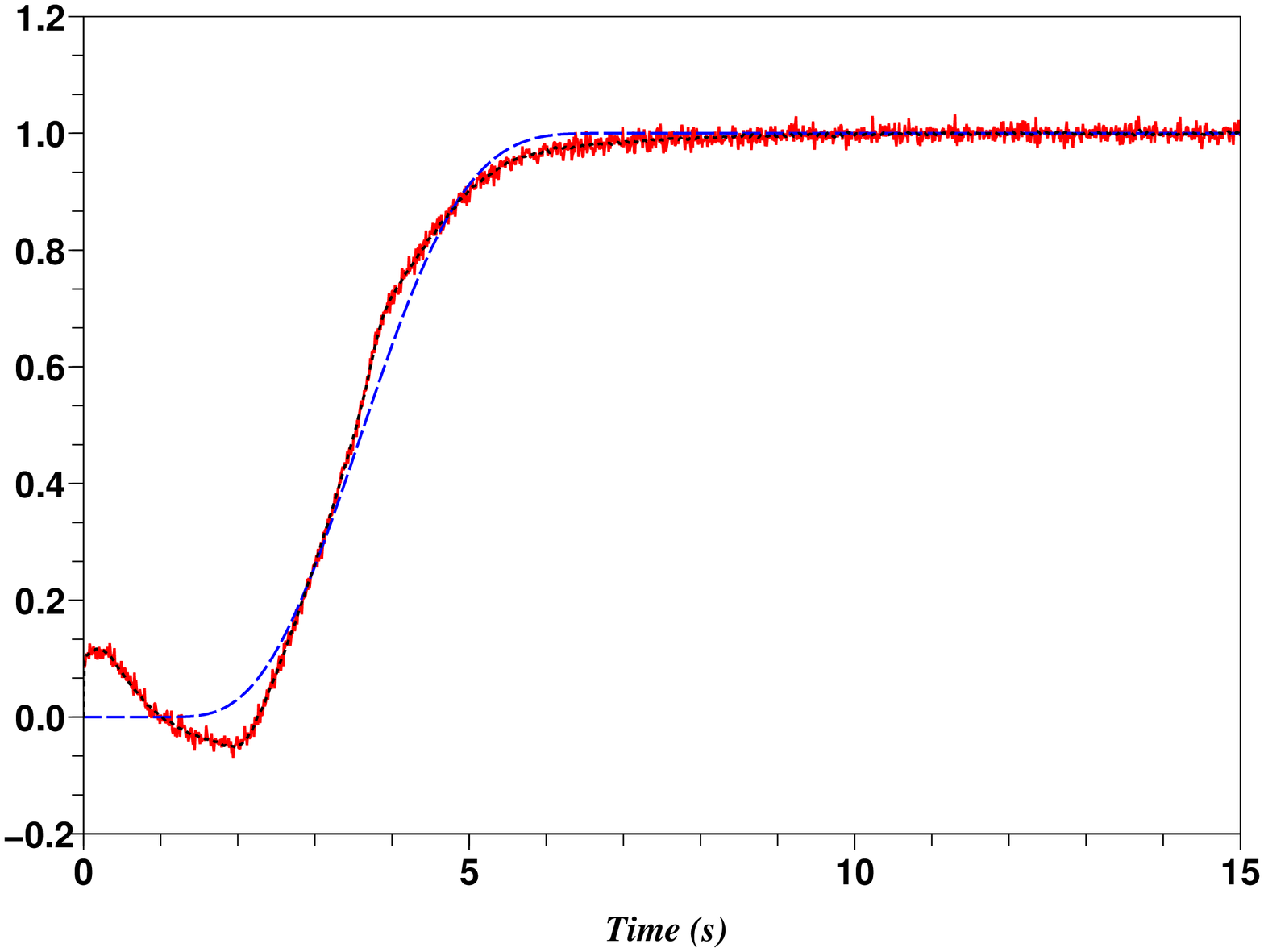}}}}}
 \caption{Instable nonlinear monovariable system\label{fig_NLA}}
\end{figure*}
\begin{figure*}[H]
\centering {\subfigure[\footnotesize
Commande]{\rotatebox{-90}{\resizebox{!}{5cm}{\includegraphics{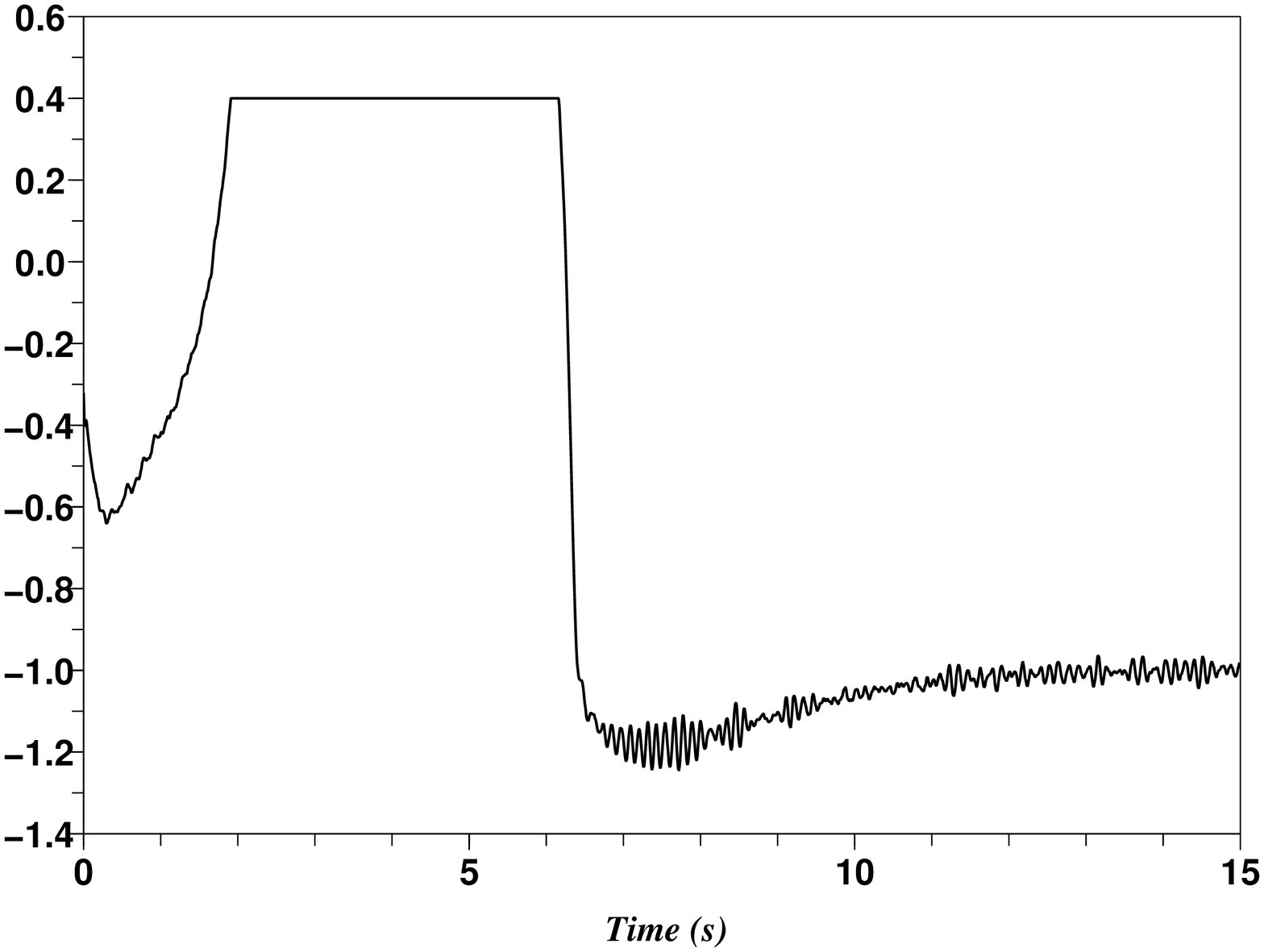}}}}}
{\subfigure[\footnotesize Output (--); reference (- -); denoised
output (.
.)]{\rotatebox{-90}{\resizebox{!}{5cm}{\includegraphics{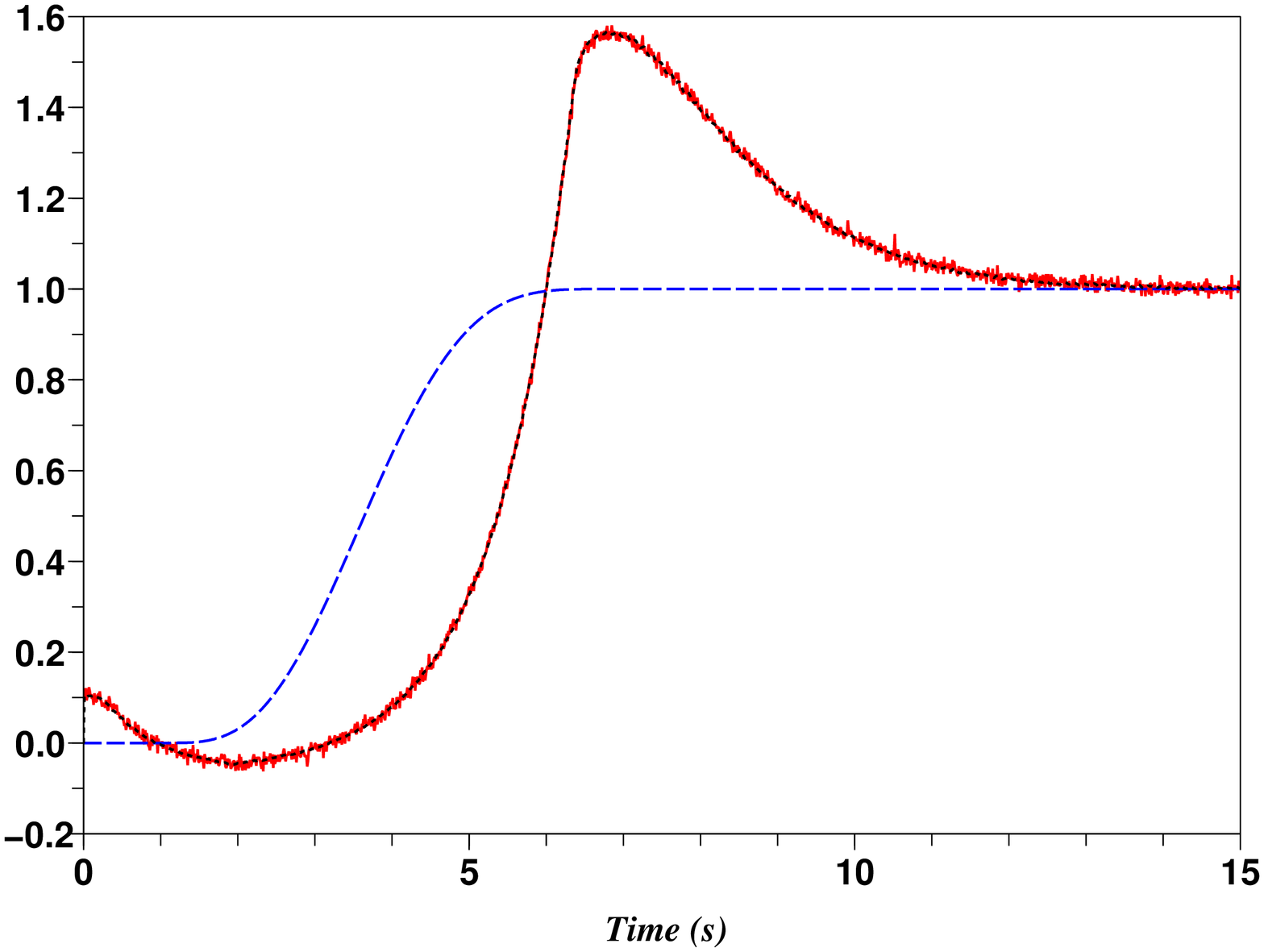}}}}}
 \caption{Instable nonlinear system: saturated control without anti-windup\label{fig_NLAsw}}
\end{figure*}

\begin{figure*}[H]
\centering {\subfigure[\footnotesize
Commande]{\rotatebox{-90}{\resizebox{!}{5cm}{\includegraphics{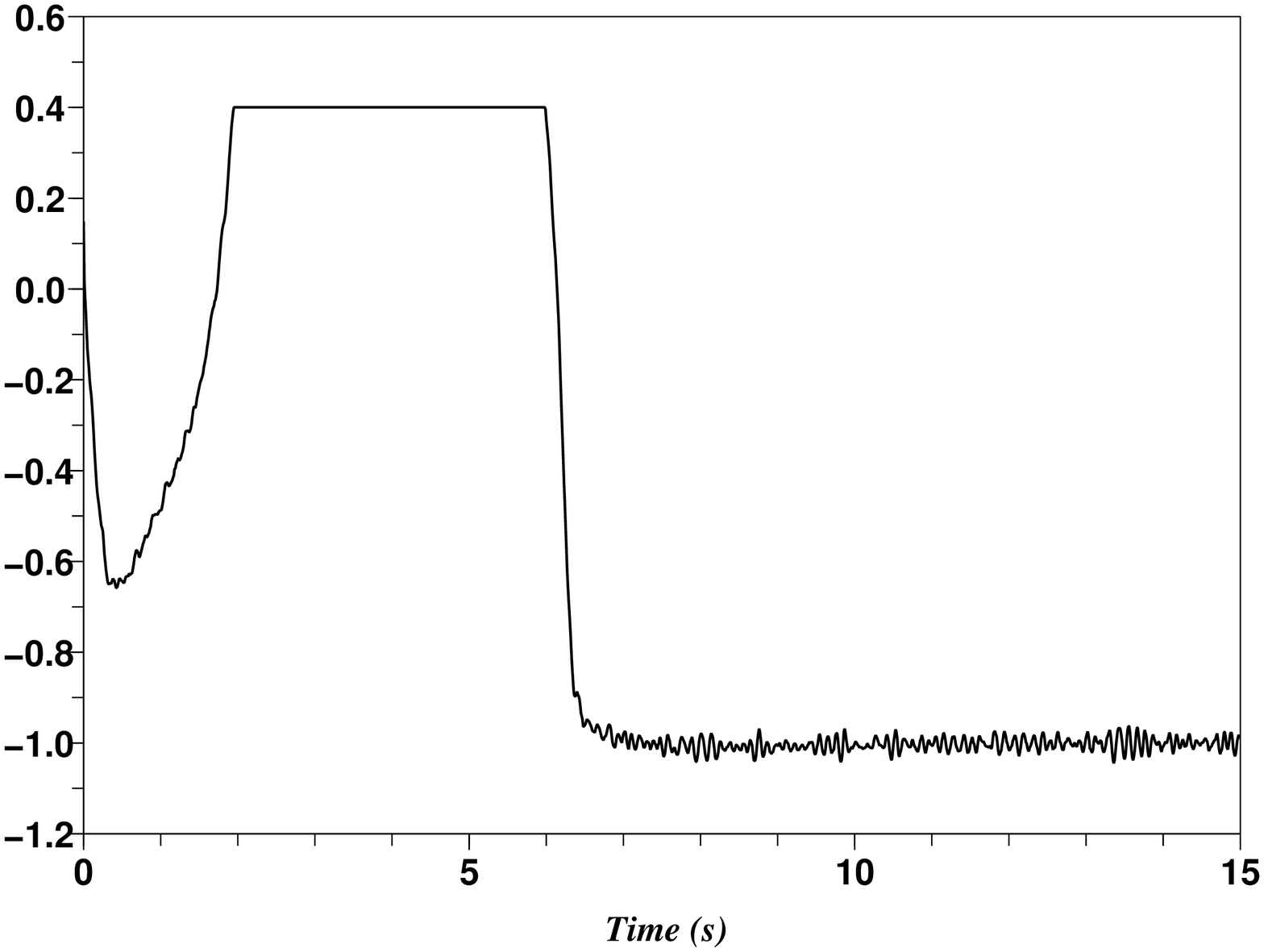}}}}}
{\subfigure[\footnotesize Output (--); reference (- -); denoised
output (.
.)]{\rotatebox{-90}{\resizebox{!}{5cm}{\includegraphics{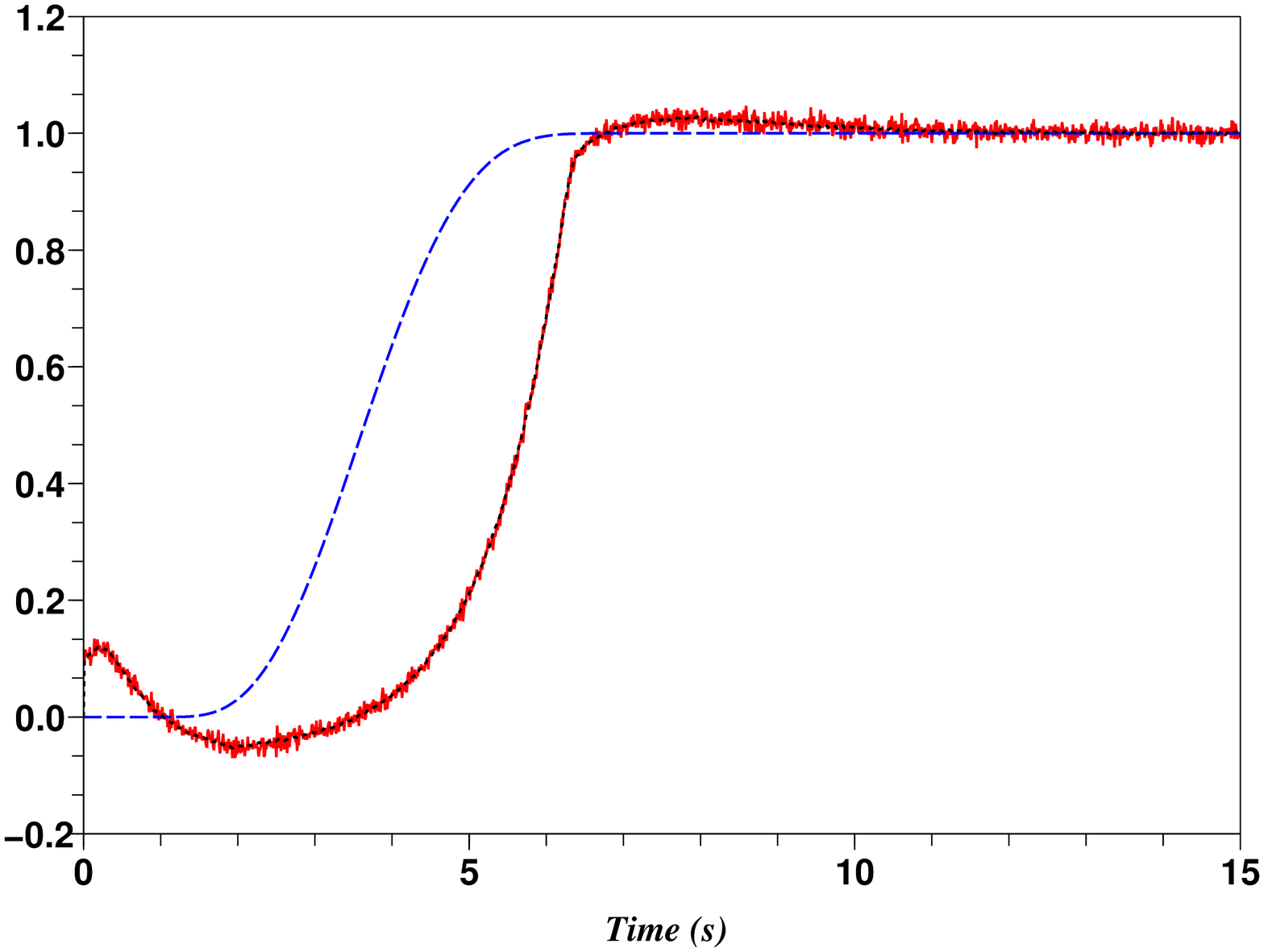}}}}}
 \caption{Instable nonlinear system: saturated control with anti-windup\label{fig_NLAaw}}
\end{figure*}

\begin{figure*}[H]
\centering%
\epsfig{figure=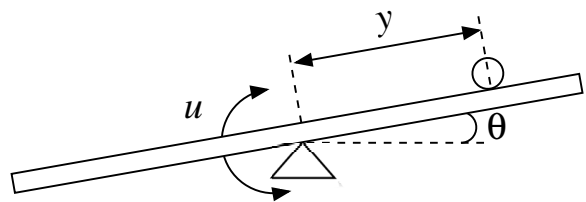,width= 0.5\textwidth}
\caption{The ball and beam example}\label{BB}
\end{figure*}

\begin{figure*}[H]
\centering%
\subfigure[\footnotesize Reference (- -) and output]
{\epsfig{figure=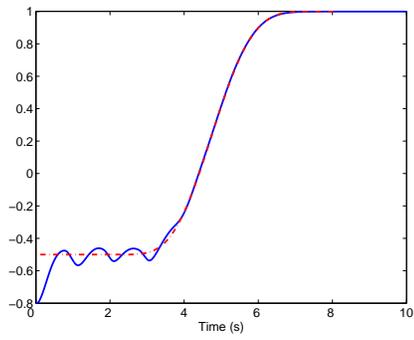,width= 0.4\textwidth}}
\subfigure[\footnotesize Control]
{\epsfig{figure=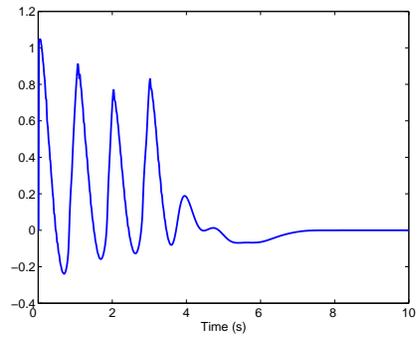,width= 0.4\textwidth}}\\%
\subfigure[\footnotesize Estimation of $F$]
{\epsfig{figure=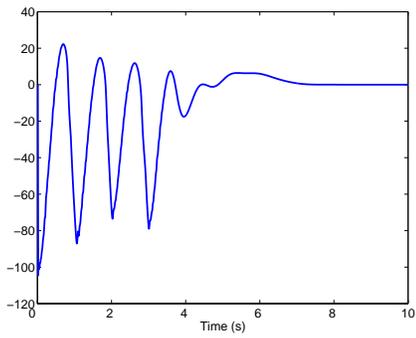,width= 0.4\textwidth}}
\subfigure[\footnotesize Reference (- -); output in the noisy case]
{\epsfig{figure=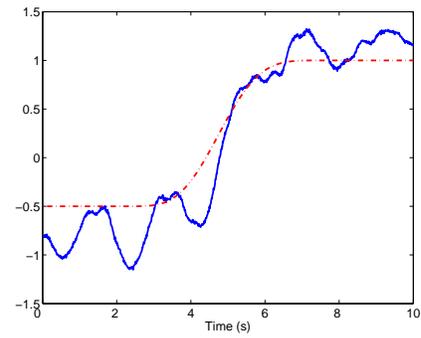,width= 0.4\textwidth}}
\caption{Polynomial trajectory for the ball and beam}\label{poly}
\end{figure*}
\begin{figure*}[H]
\centering%
\subfigure[\footnotesize Reference (- -); output]
{\epsfig{figure=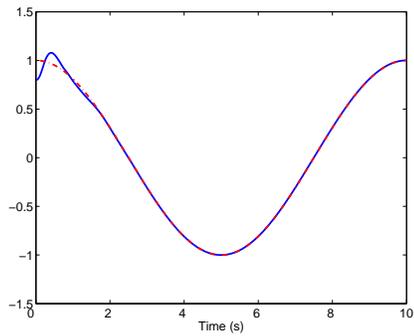,width= 0.4\textwidth}}
\subfigure[\footnotesize Control]
{\epsfig{figure=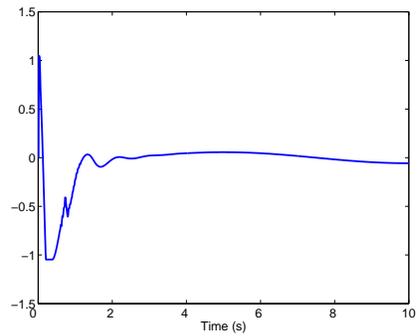,width= 0.4\textwidth}}\\
\subfigure[\footnotesize Estimation of $F$]
{\epsfig{figure=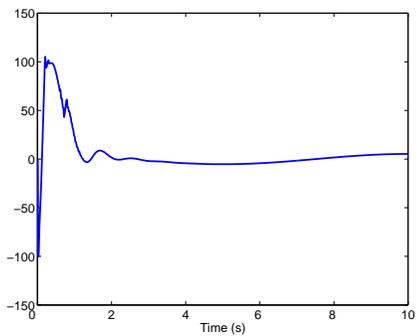,width= 0.4\textwidth}}
\subfigure[\footnotesize Reference (- -); output in the noisy case]
{\epsfig{figure=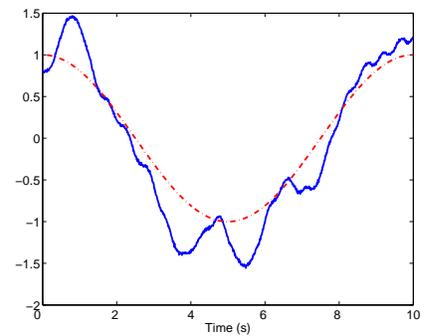,width= 0.4\textwidth}}
\caption{Sinusoidal trajectory for the ball and beam}\label{sinus}
\end{figure*}

\begin{figure*}[H]
\centering
%
%
{\includegraphics[scale=0.7]{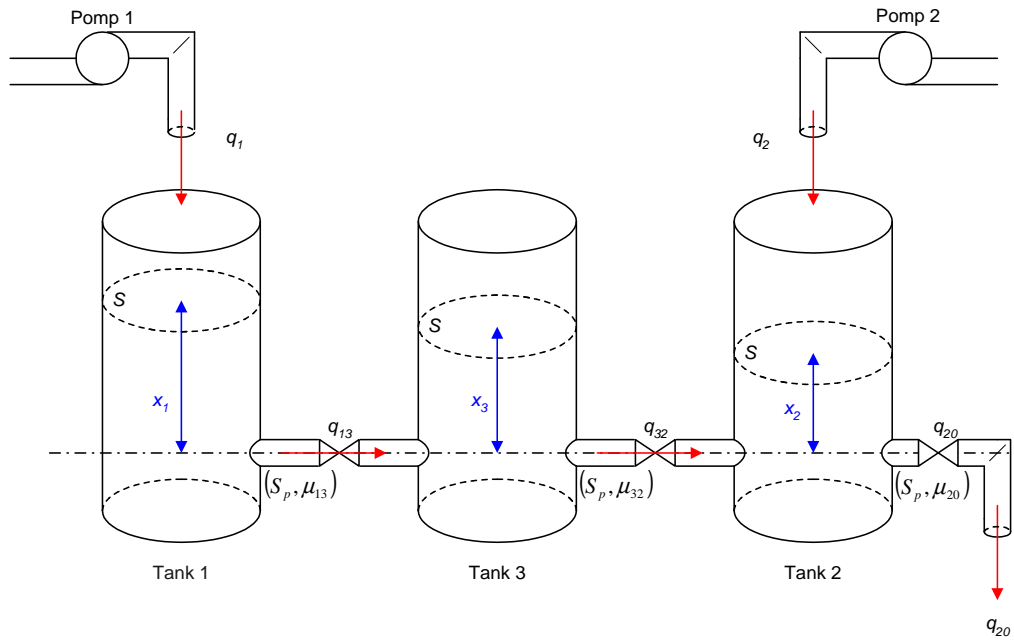}} \caption{The 3 tank
system \label{schem}}
\end{figure*}

\begin{figure*}[H]
\centering {\subfigure[\footnotesize References (- -); outputs]{
\rotatebox{-90}{\resizebox{!}{5cm}{%
   \includegraphics{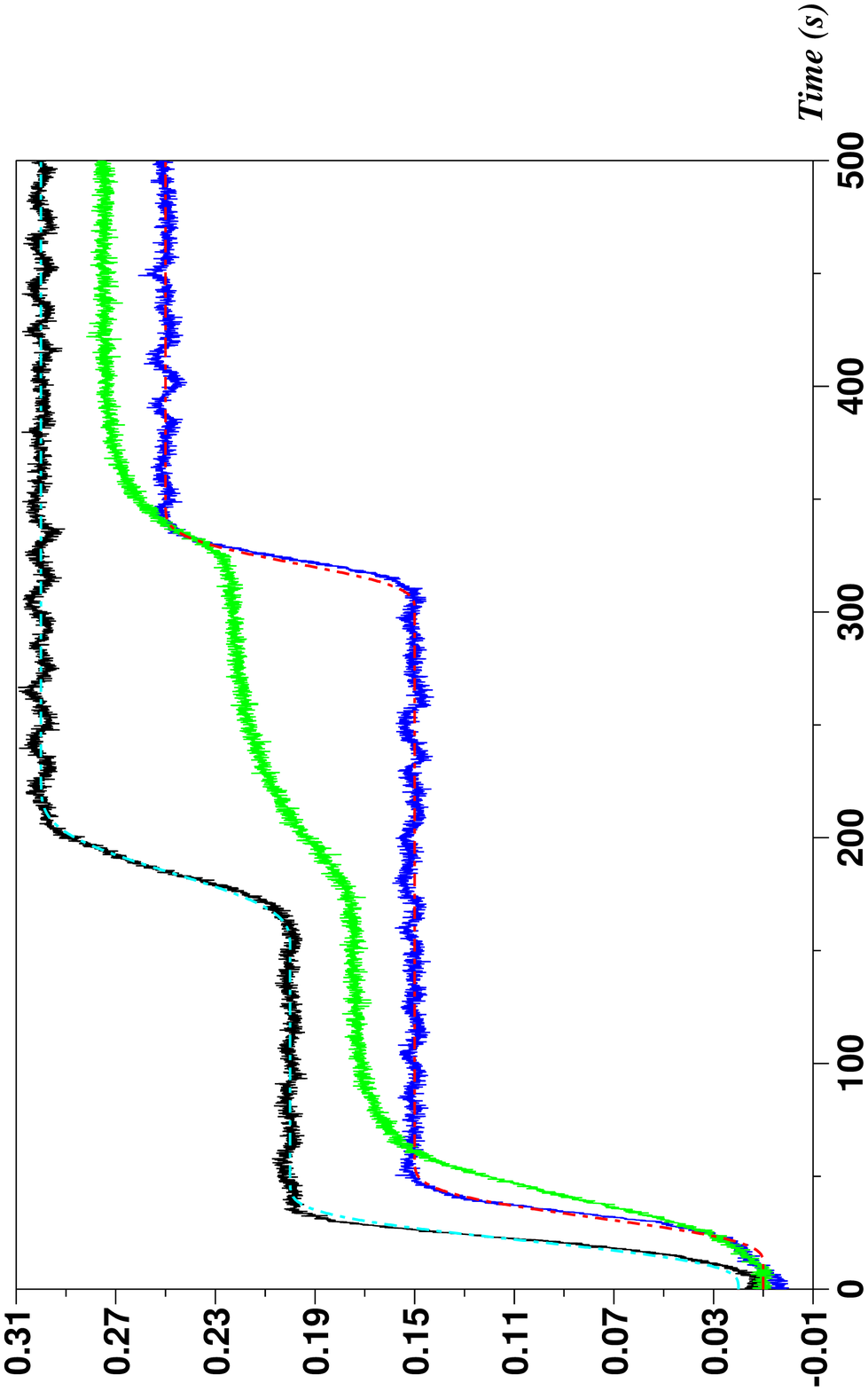}}}}}
{\subfigure[\footnotesize Estimation of the derivatives $\dot y_1$
(-) and $\dot y_2$ (- -), which is shifted back $0.02$]{
\rotatebox{-90}{\resizebox{!}{5cm}{\includegraphics{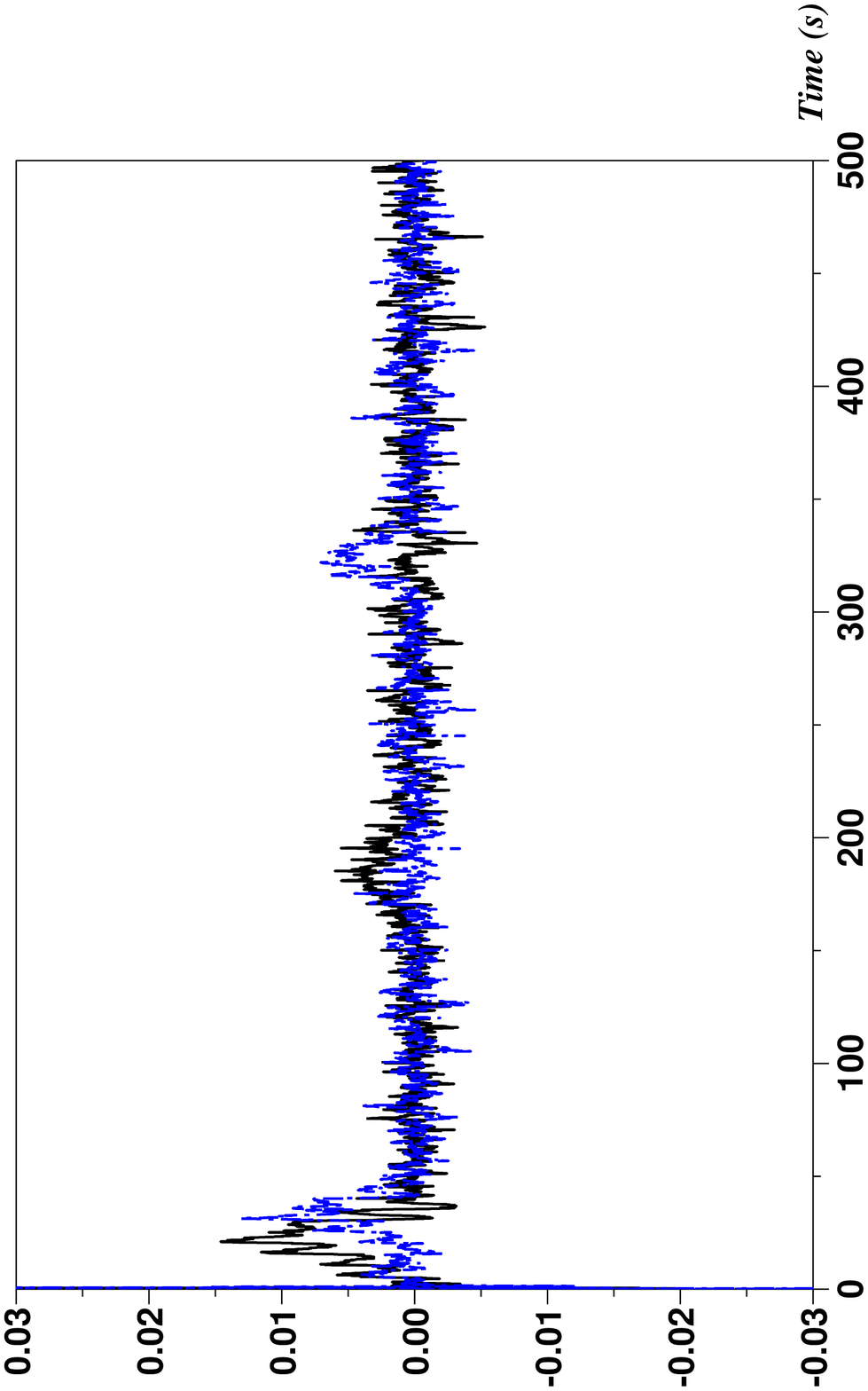}}}}}
{\subfigure[\footnotesize Control $u_1$ (-) and $u_2$ (- -)]{
\rotatebox{-90}{\resizebox{!}{5cm}{%
   \includegraphics{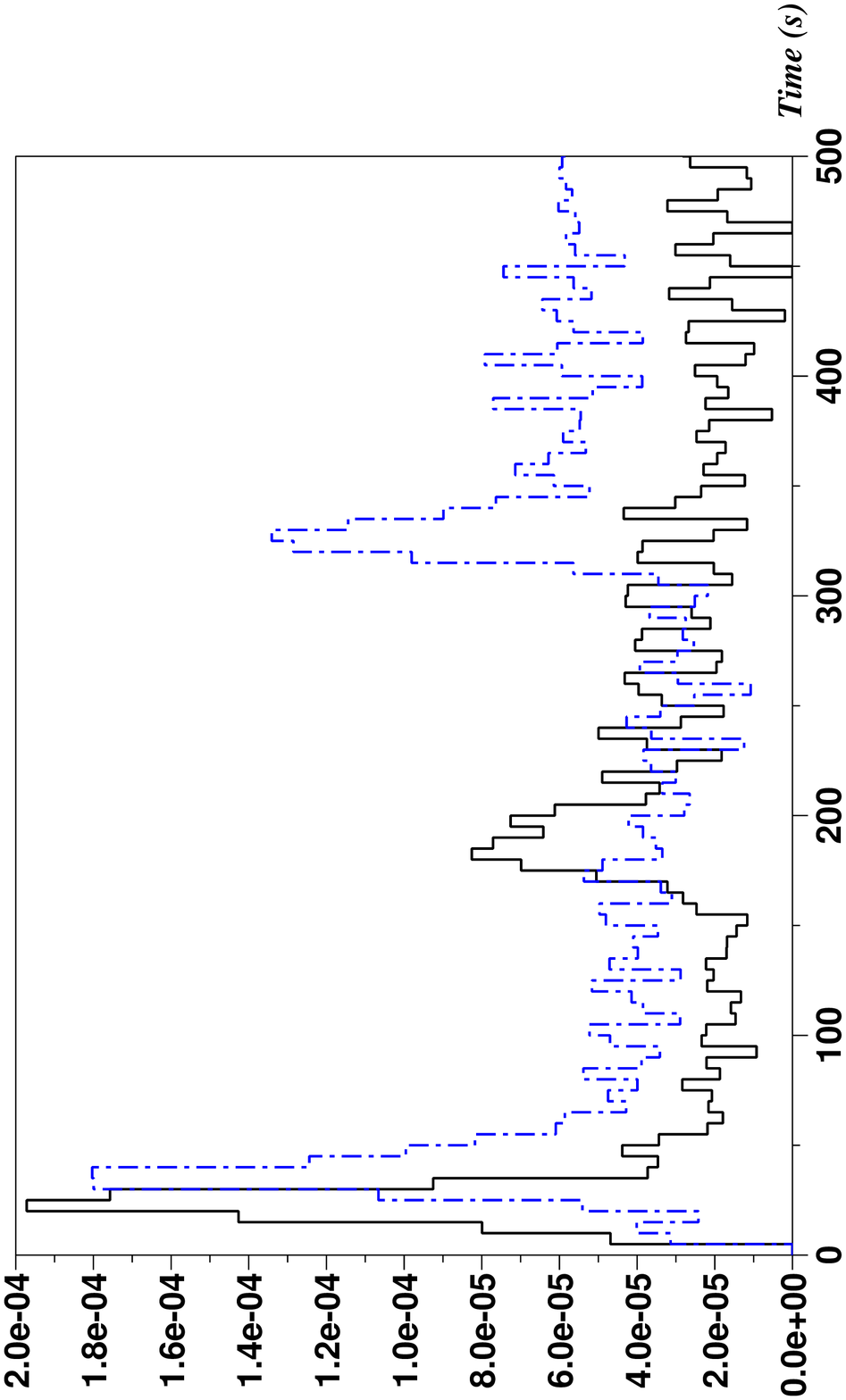}}}}}
{\subfigure[\footnotesize Denoising]{
\rotatebox{-90}{\resizebox{!}{5cm}{%
   \includegraphics{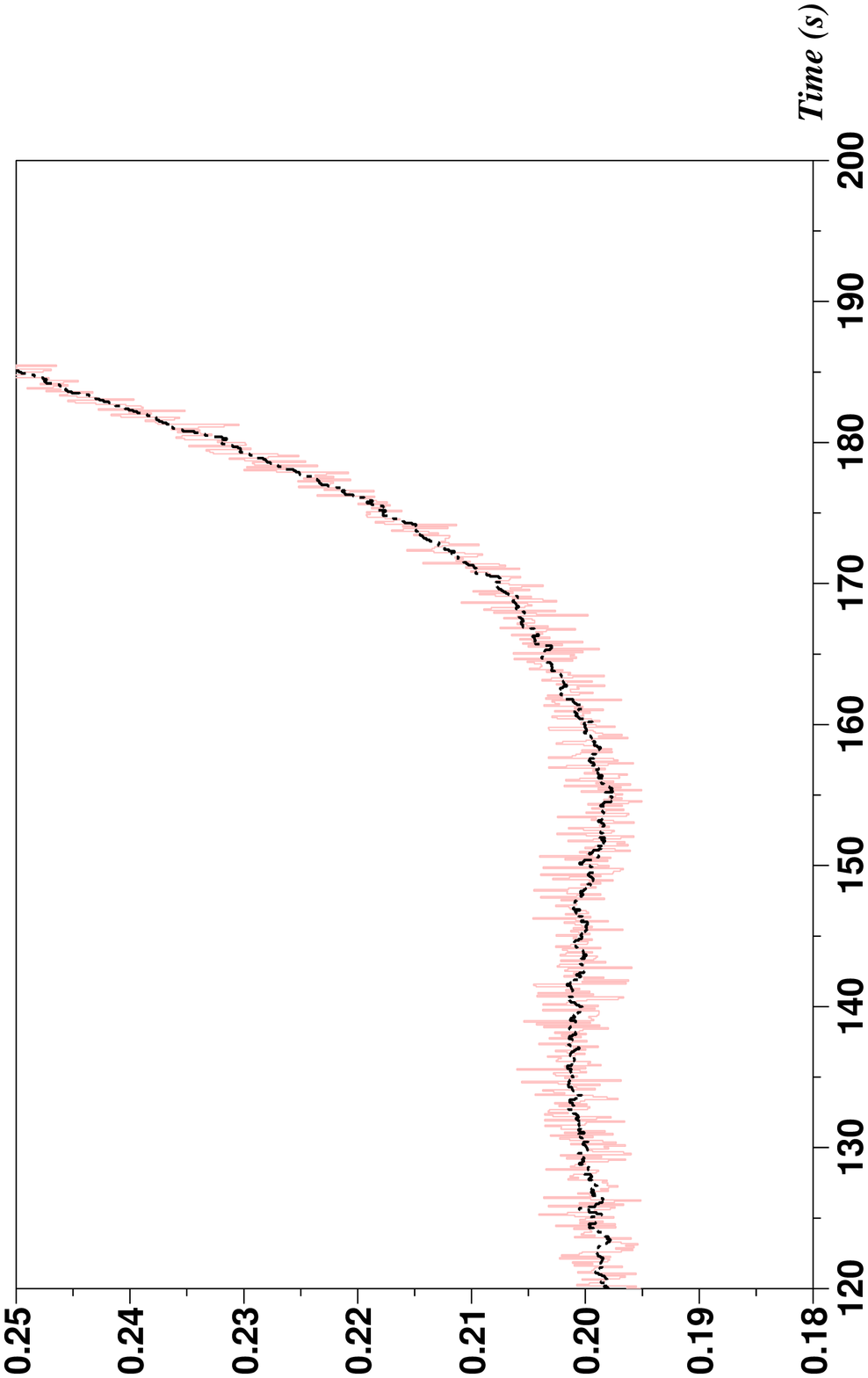}}}}}
 \caption{Simulations for the 3 tank system\label{fig_3cuves}}
\end{figure*}
\clearpage
\begin{figure*}[H]
\vspace{-0.5cm} \centering {\subfigure[\footnotesize The friction
law]{
\rotatebox{-90}{\resizebox{!}{5cm}{%
   \includegraphics{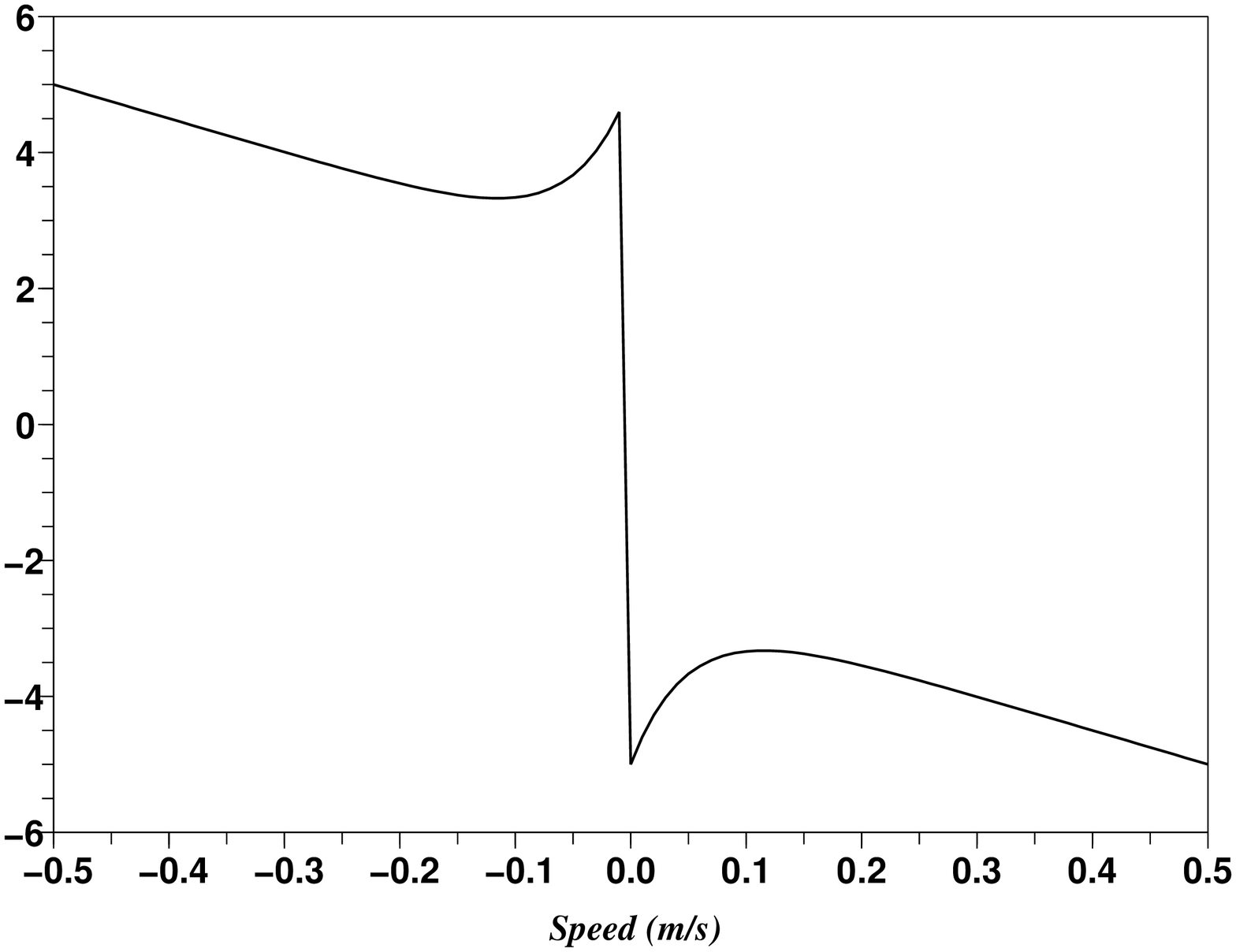}}}}}
\vspace{-0.2cm} {\subfigure[\footnotesize Friction
${\mathcal{F}}(\dot y)$ (--); Estimation of the whole unknown
effects \text{$[{\mathcal{G}}(\dot y)]_e$} (-
-)]{\rotatebox{-90}{\resizebox{!}{5cm}{\includegraphics{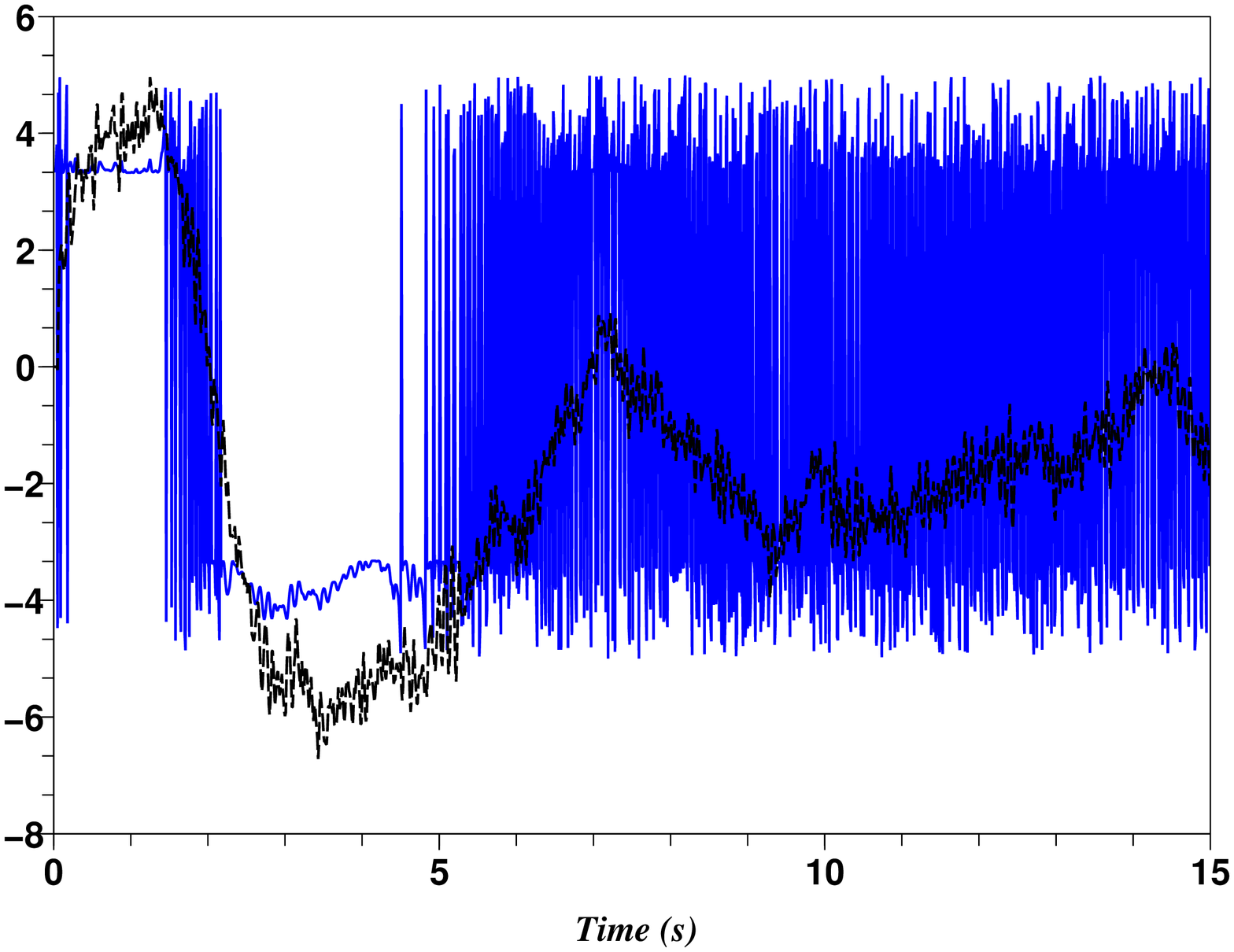}}}}}
\vspace{-0.3cm} {\subfigure[\footnotesize i-PID control]{
\rotatebox{-90}{\resizebox{!}{5cm}{%
   \includegraphics{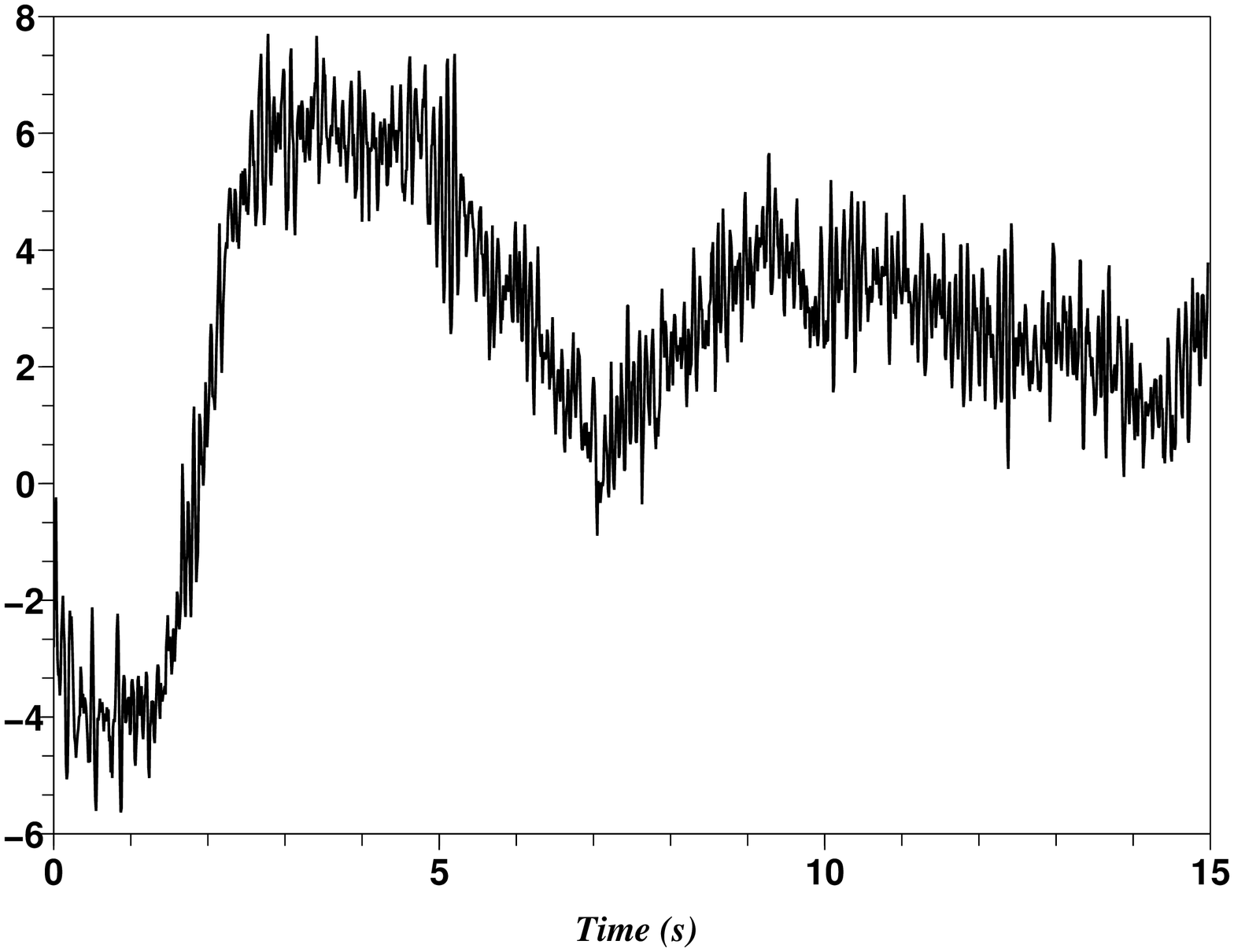}}}}}
{\subfigure[\footnotesize i-PID: output (--), reference (- -);
denoised output (. .)]{
\rotatebox{-90}{\resizebox{!}{5cm}{%
   \includegraphics{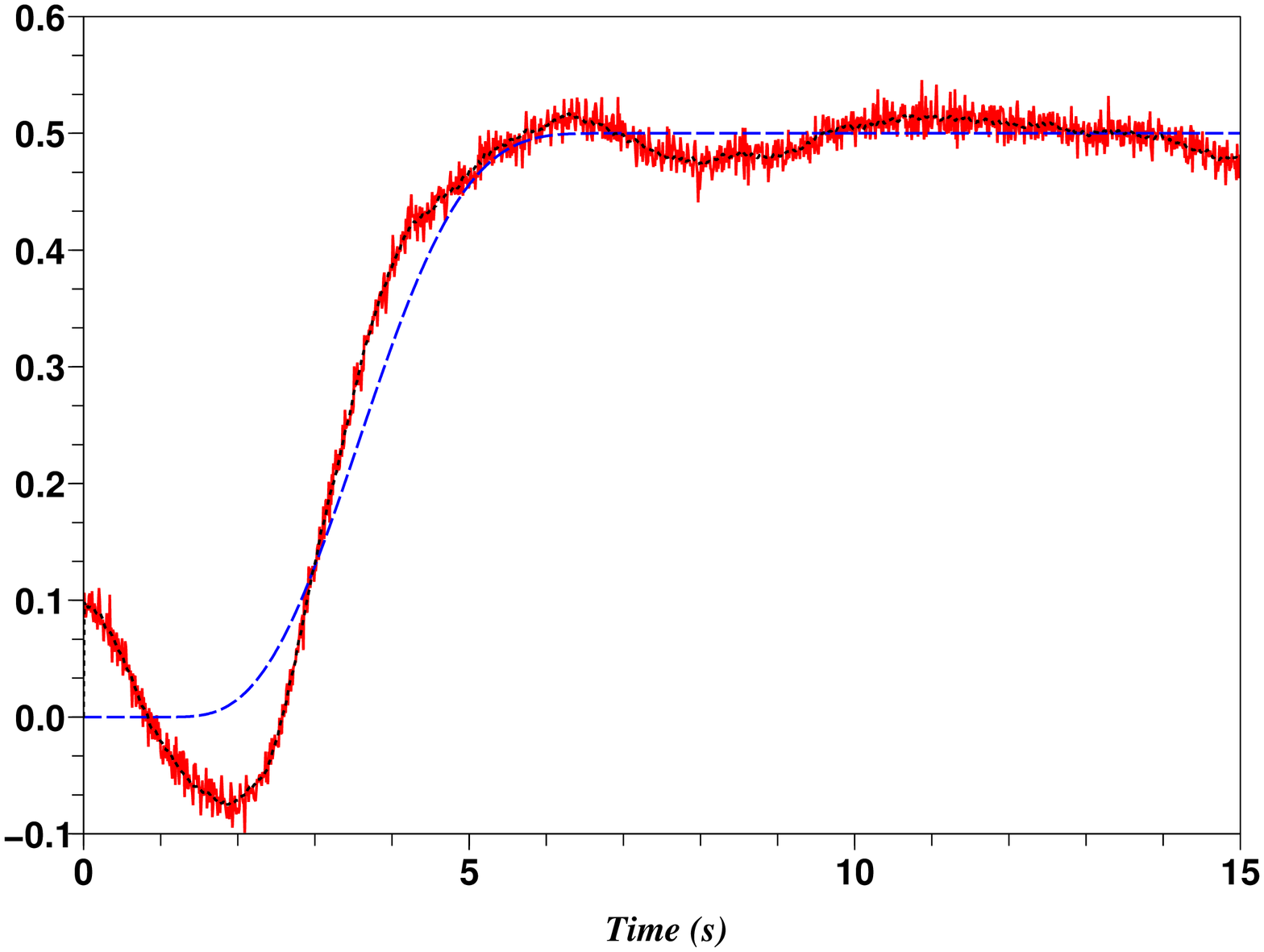}}}}}
 \vspace{-0.3cm}
{\subfigure[\footnotesize Flatness-based control + classic PID ]{
\rotatebox{-90}{\resizebox{!}{5cm}{%
   \includegraphics{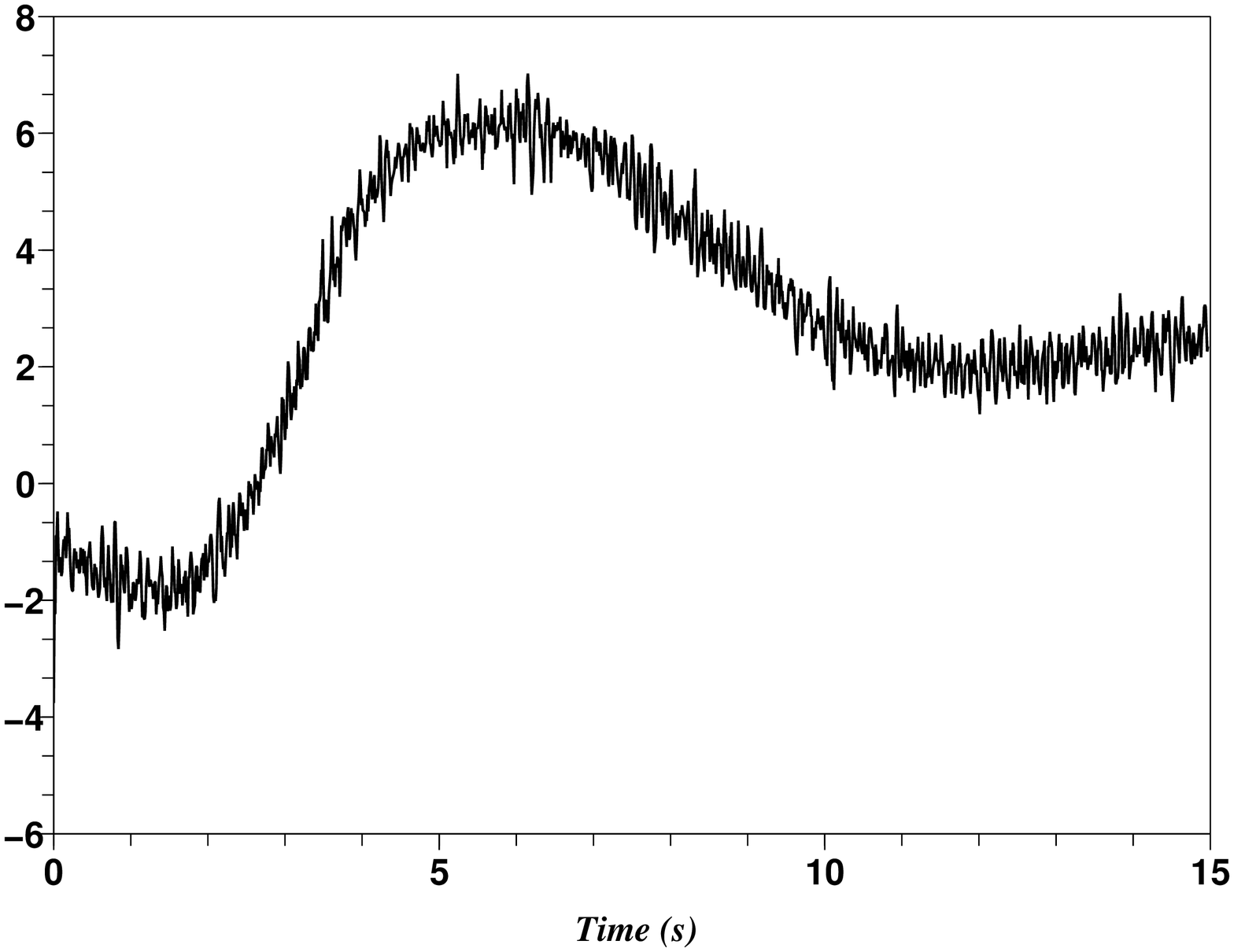}}}}}
{\subfigure[\footnotesize Flatness-based control + classic PID:
output (--), reference (- -), denoised output (. .)]{
\rotatebox{-90}{\resizebox{!}{5cm}{%
   \includegraphics{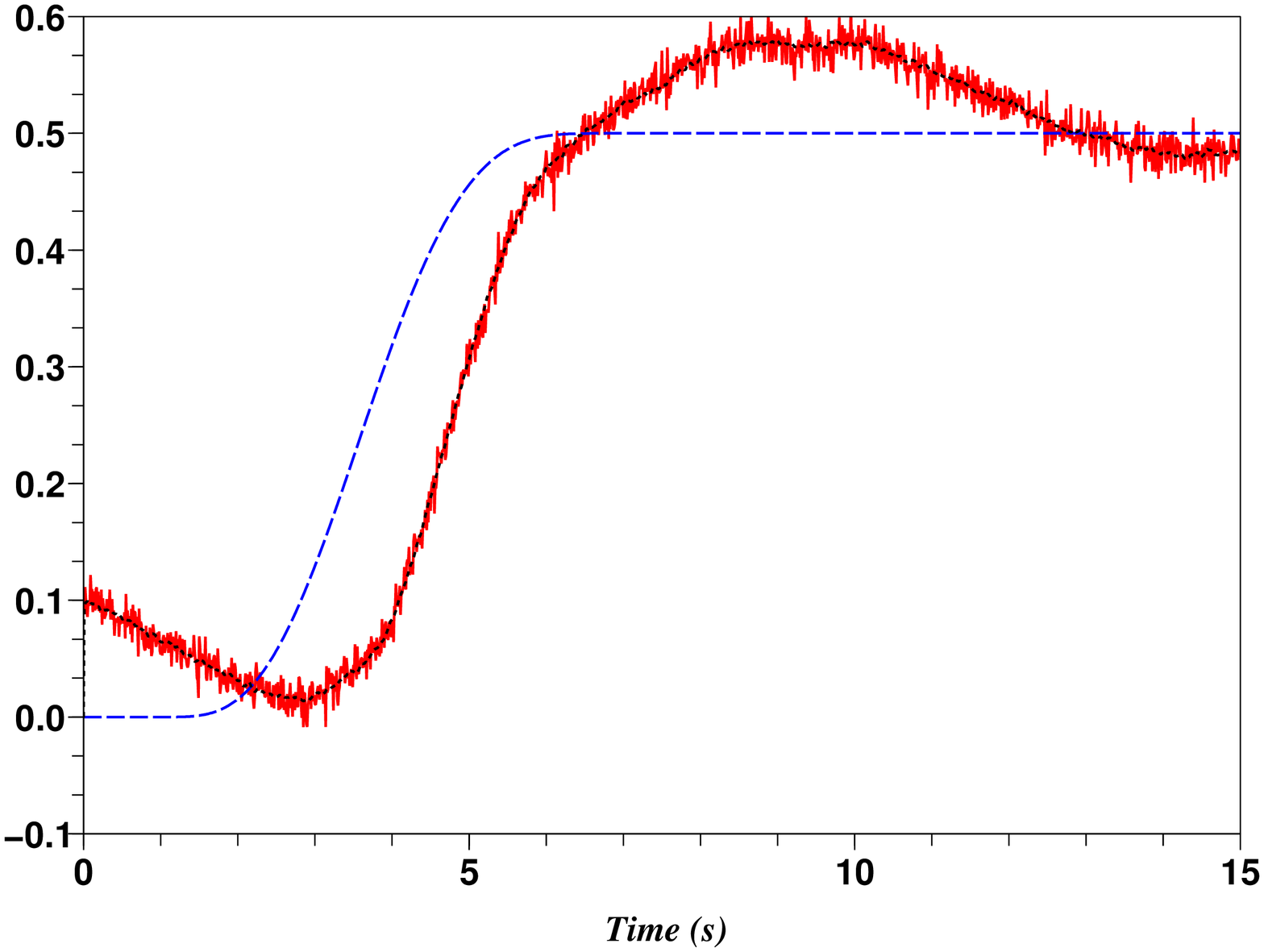}}}}}
 \vspace{-0.4cm}
{\subfigure[\footnotesize Non-flatness-based control + PID]{
\rotatebox{-90}{\resizebox{!}{5cm}{%
   \includegraphics{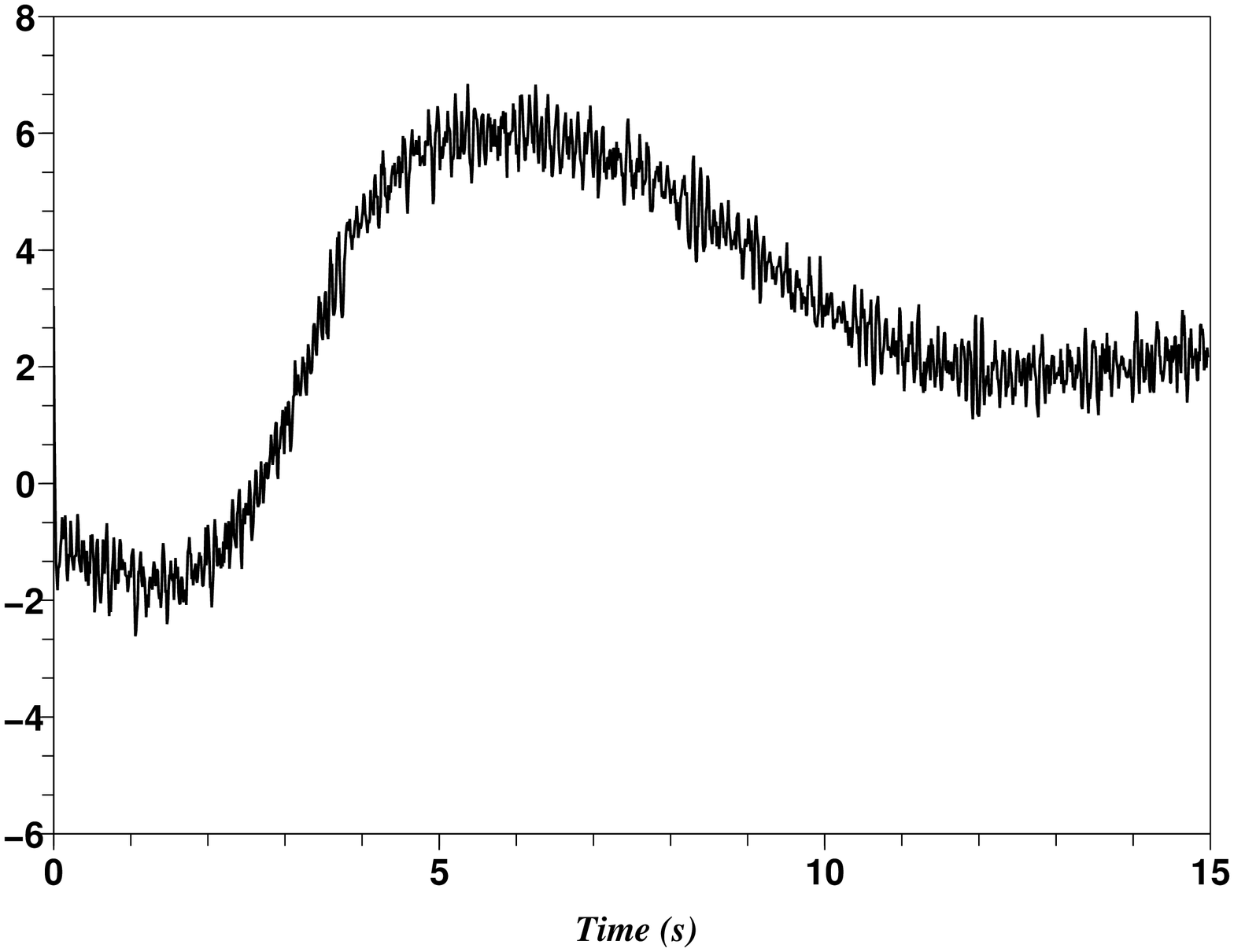}}}}}
{\subfigure[\footnotesize Non-flatness-based control: output (--),
reference (- -), denoised output (. .)]{
\rotatebox{-90}{\resizebox{!}{5cm}{%
   \includegraphics{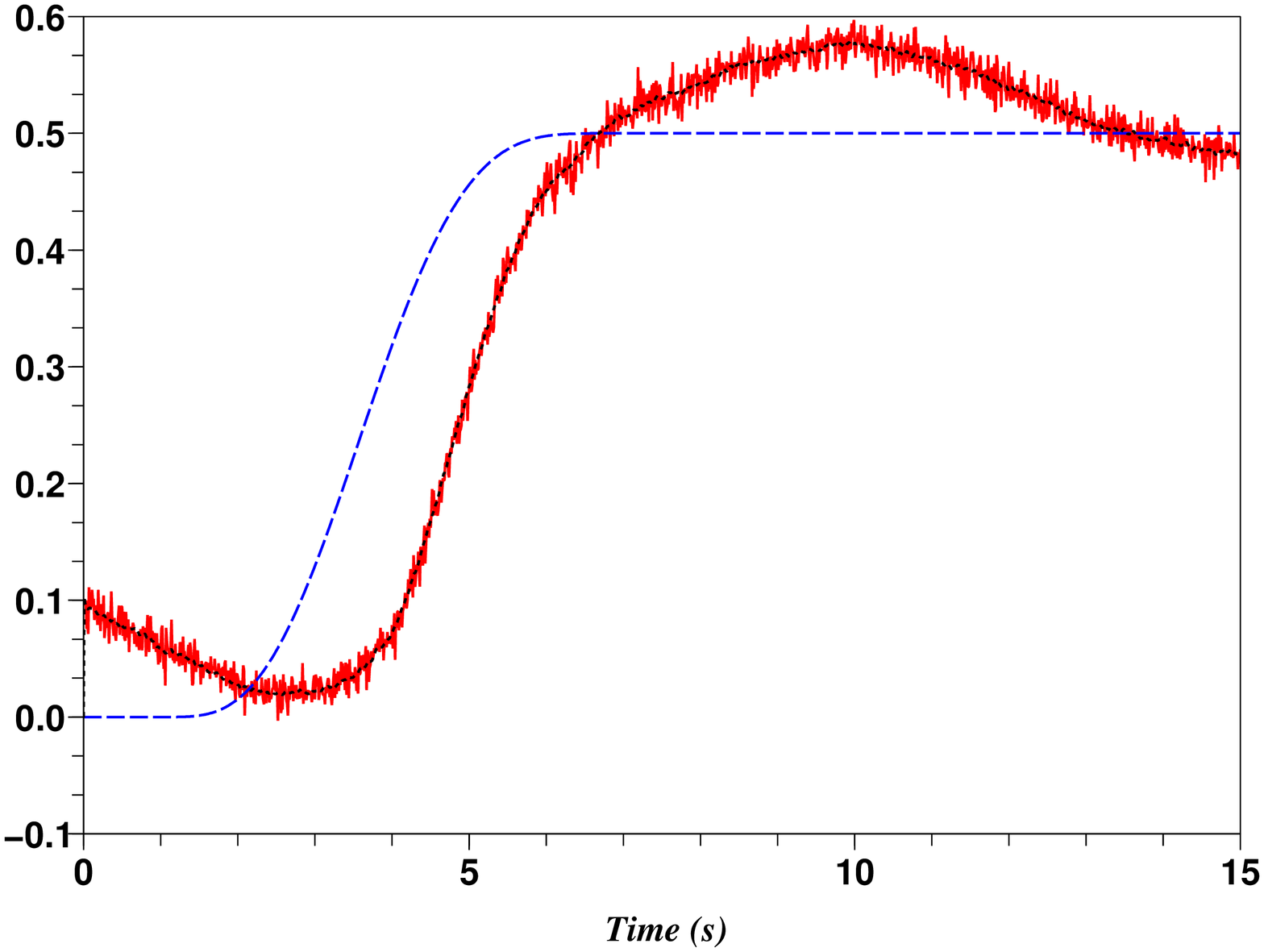}}}}}
%
 \caption{The spring with unknown with nonlinearity, friction and damping \label{fig_SLF}}
\end{figure*}

\begin{figure*}[H]
\centering {\subfigure[\footnotesize $u$ (--); $u^\star$ (- -)]{
\rotatebox{-90}{\resizebox{!}{5cm}{%
   \includegraphics{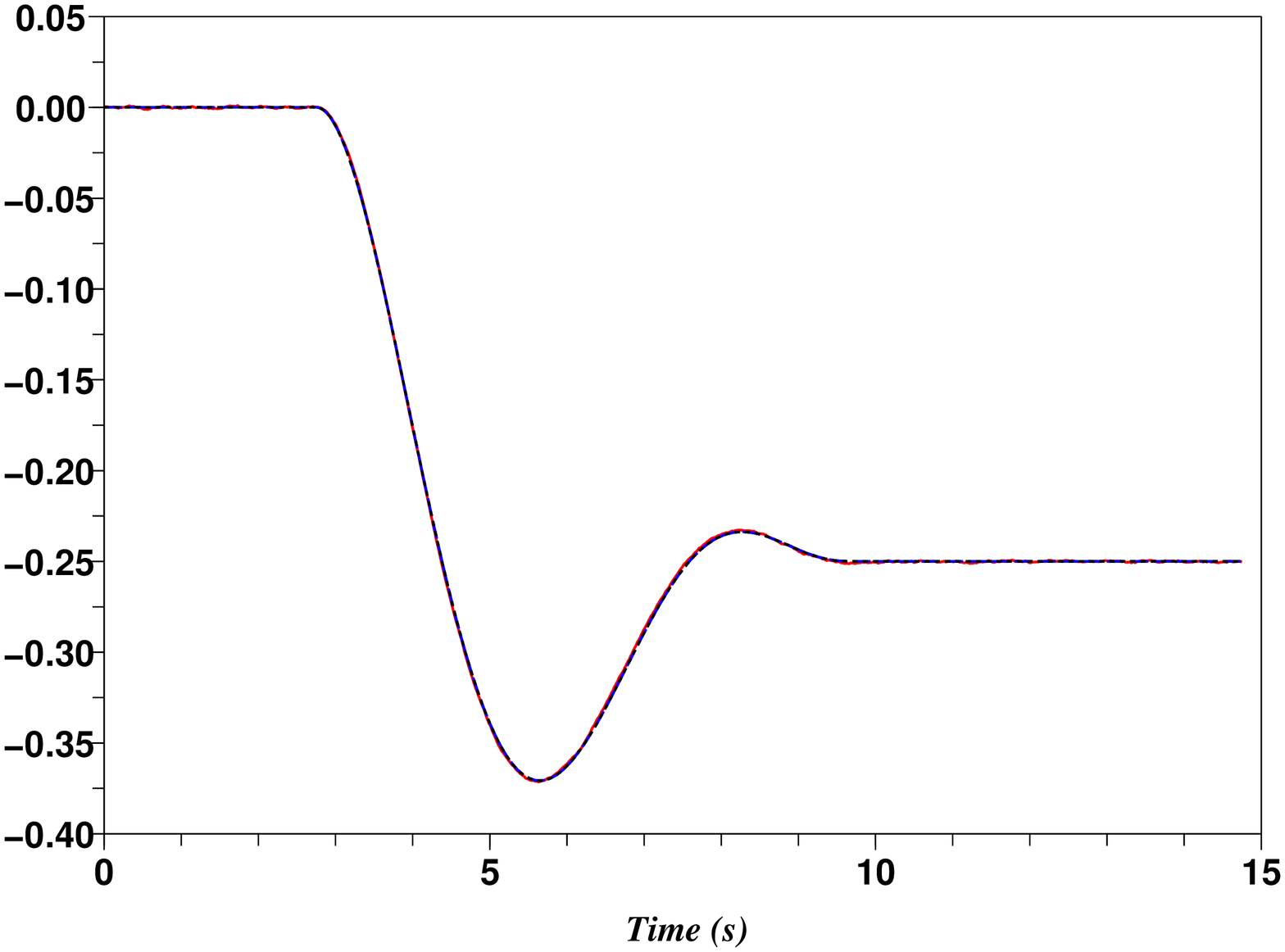}}}}}
{\subfigure[\footnotesize Output (--); reference (- -); denoised
output (. .)]{
\rotatebox{-90}{\resizebox{!}{5cm}{%
   \includegraphics{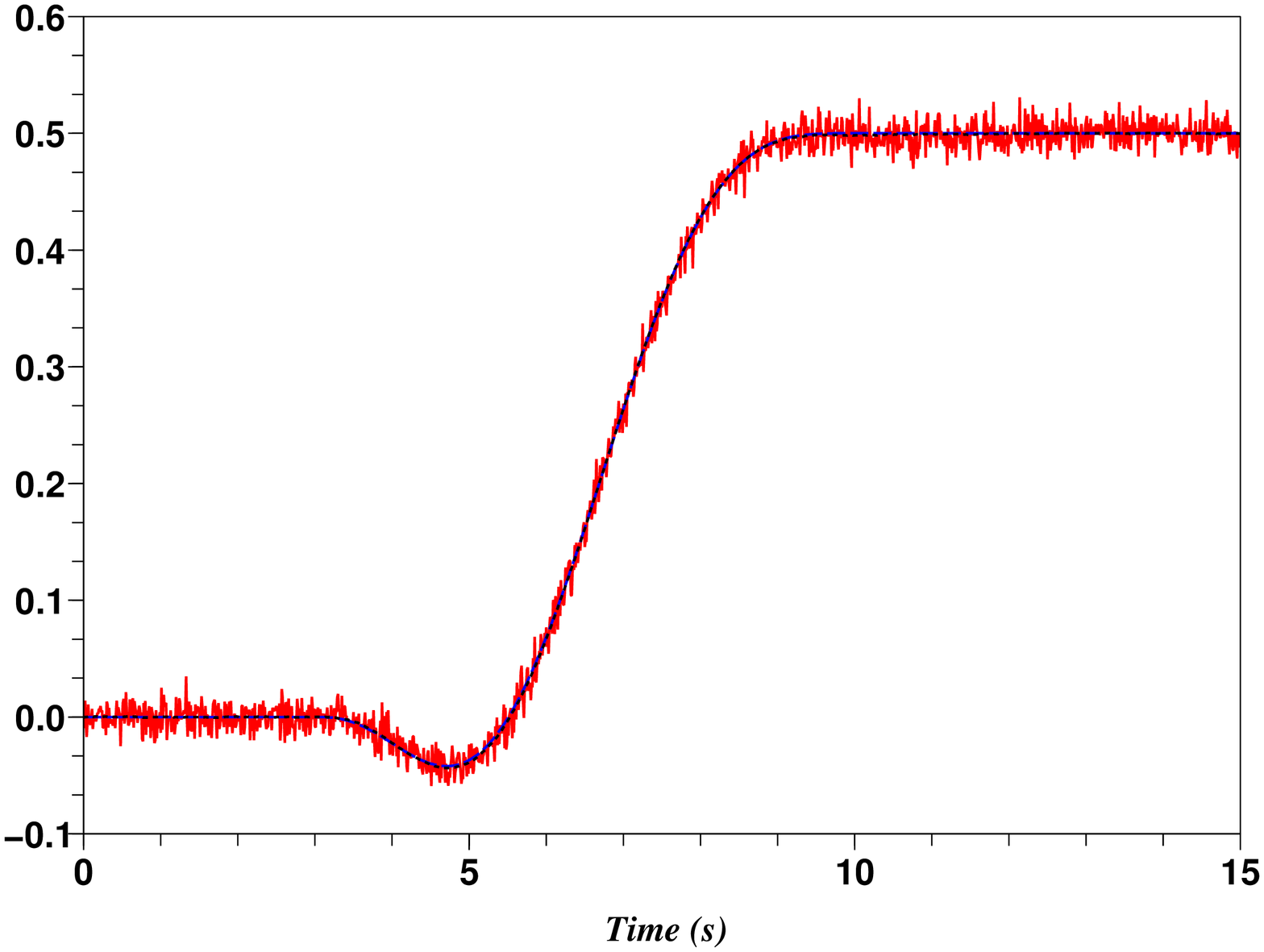}}}}}
\centering {\subfigure[\footnotesize $u$ (--); $u^\star$ (- -)]{
\rotatebox{-90}{\resizebox{!}{5cm}{%
   \includegraphics{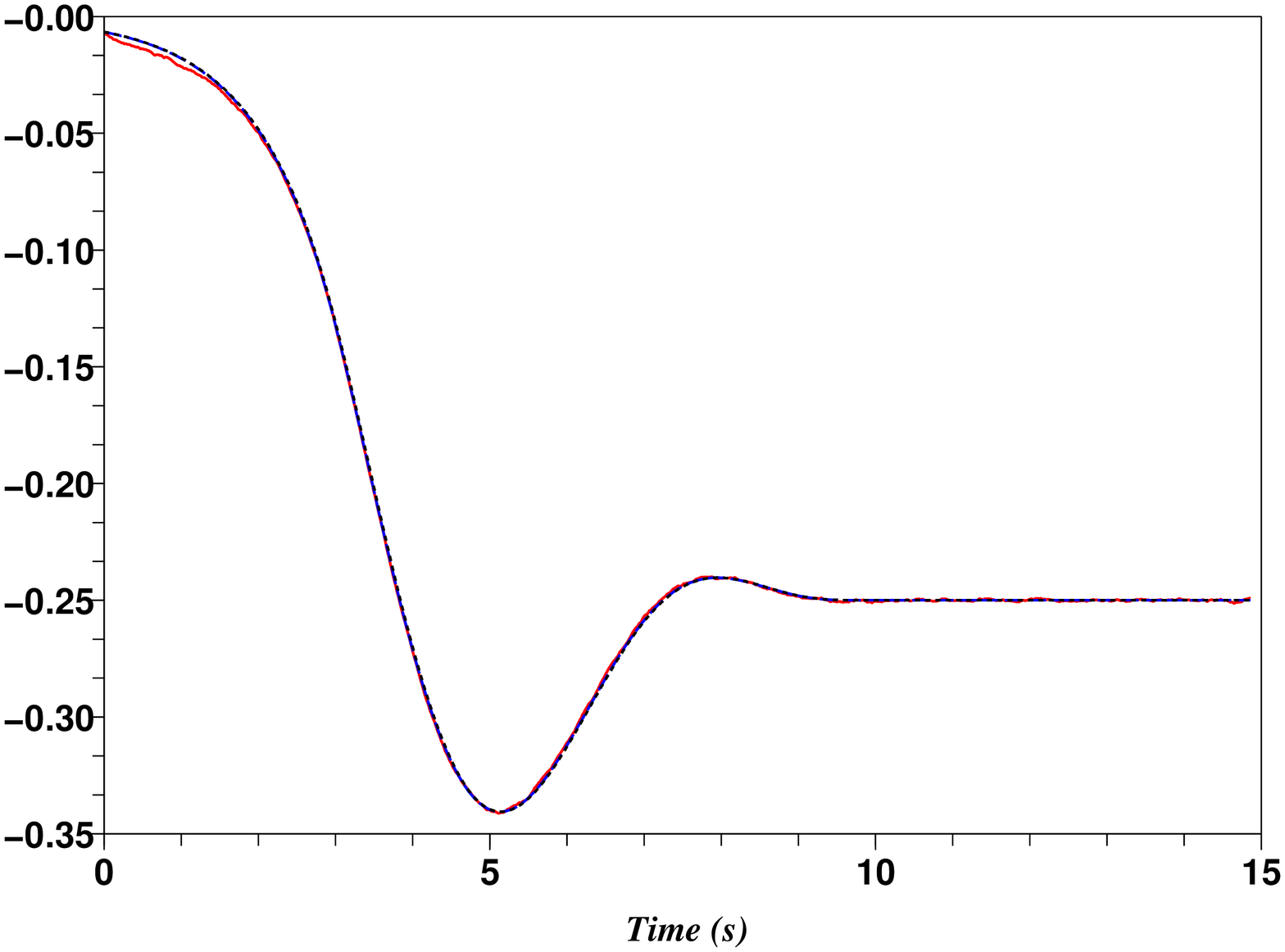}}}}}
{\subfigure[\footnotesize Output (--); reference (- -); denoised
output (. .)]{
\rotatebox{-90}{\resizebox{!}{5cm}{%
   \includegraphics{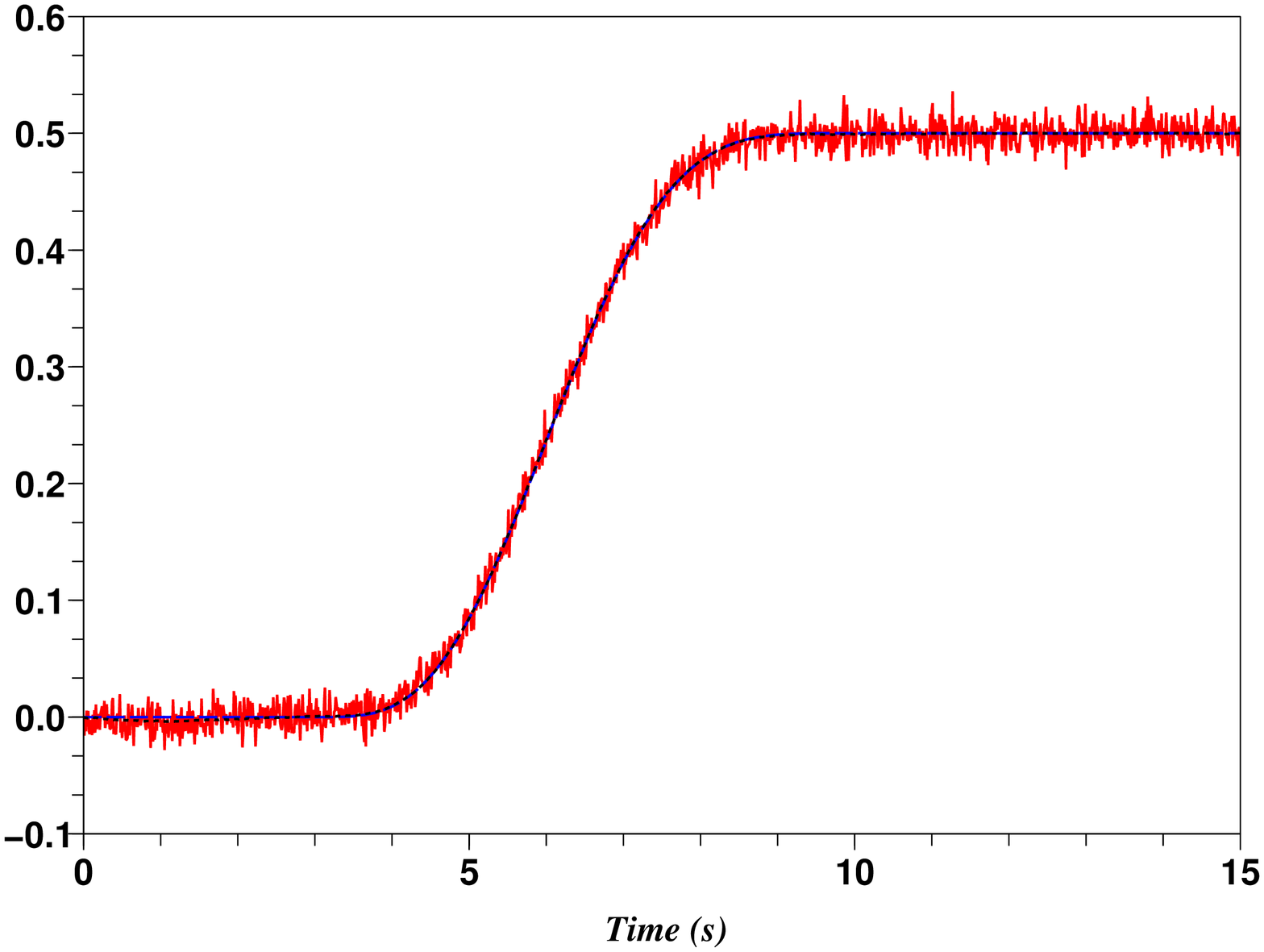}}}}}
\caption{Non-minimum phase system \label{pnm0}}
\end{figure*}

\begin{figure*}[H]
\centering {\subfigure[\footnotesize $u$ (--); $u^\star$ (- -)]{
\rotatebox{-90}{\resizebox{!}{5cm}{%
   \includegraphics{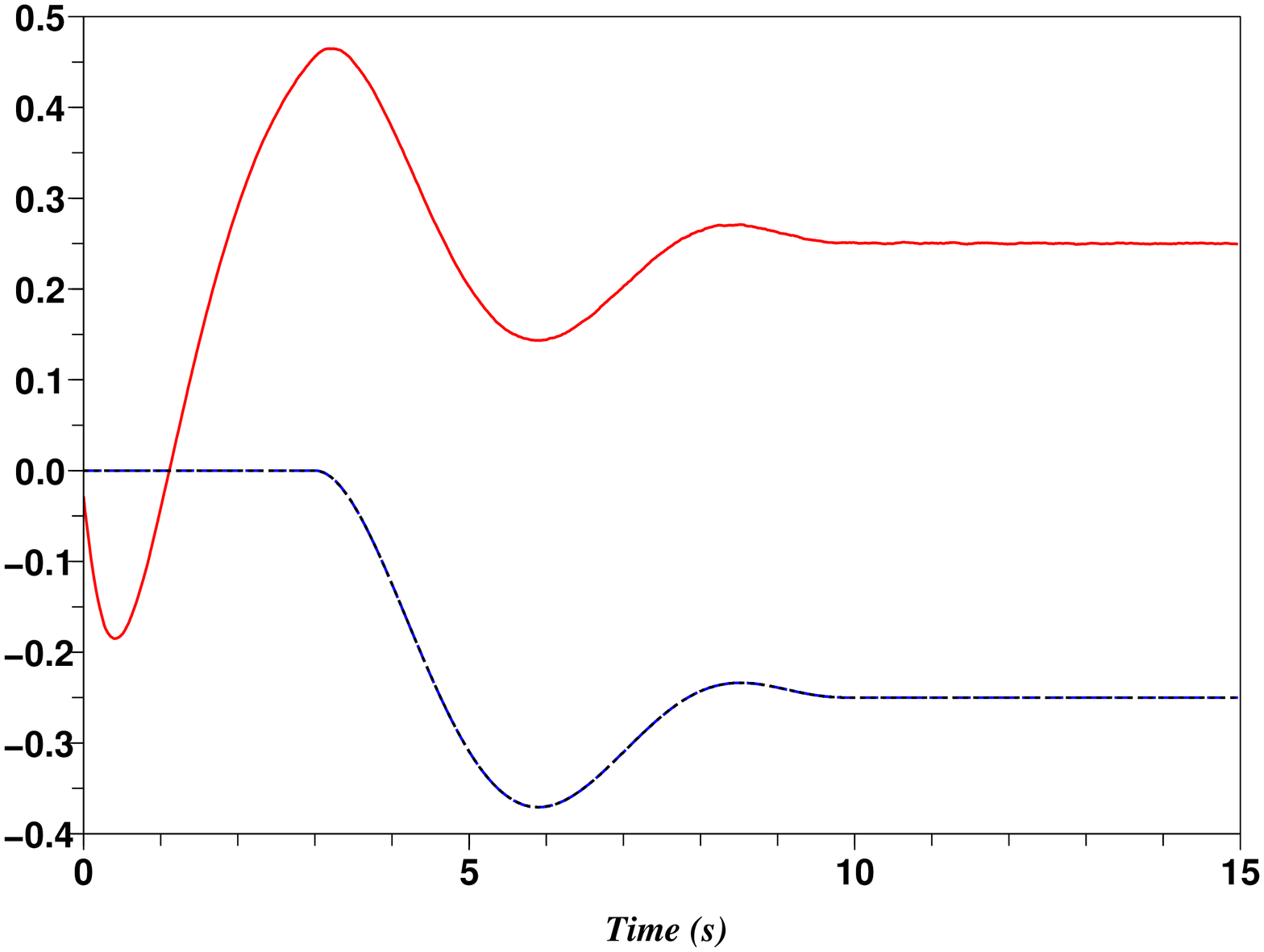}}}}}
{\subfigure[\footnotesize Output (--); reference (- -); denoised
output (. .)]{
\rotatebox{-90}{\resizebox{!}{5cm}{%
   \includegraphics{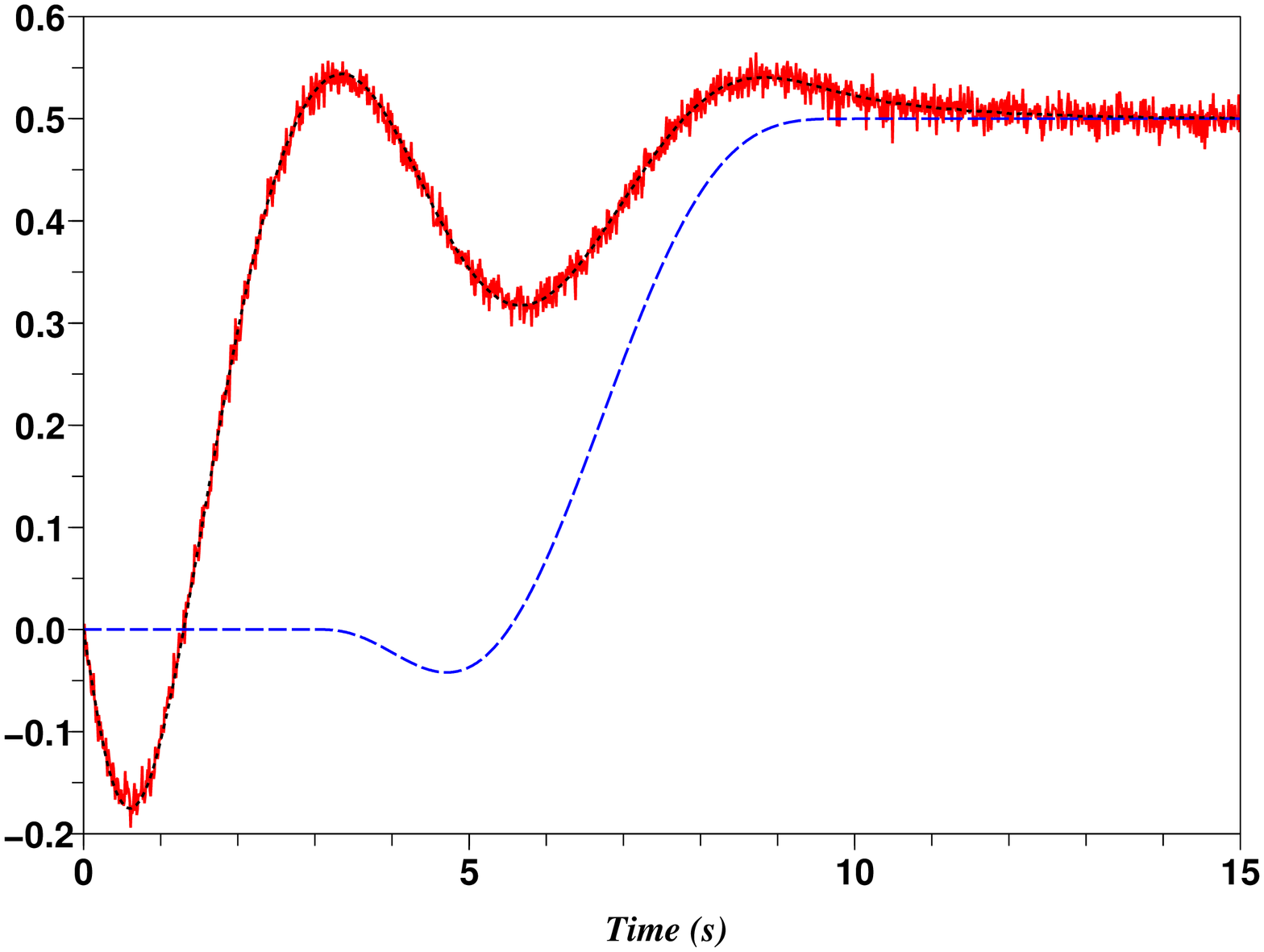}}}}}
{\subfigure[\footnotesize Perturbation (- -); estimated perturbation
(--)]{
\rotatebox{-90}{\resizebox{!}{5cm}{%
   \includegraphics{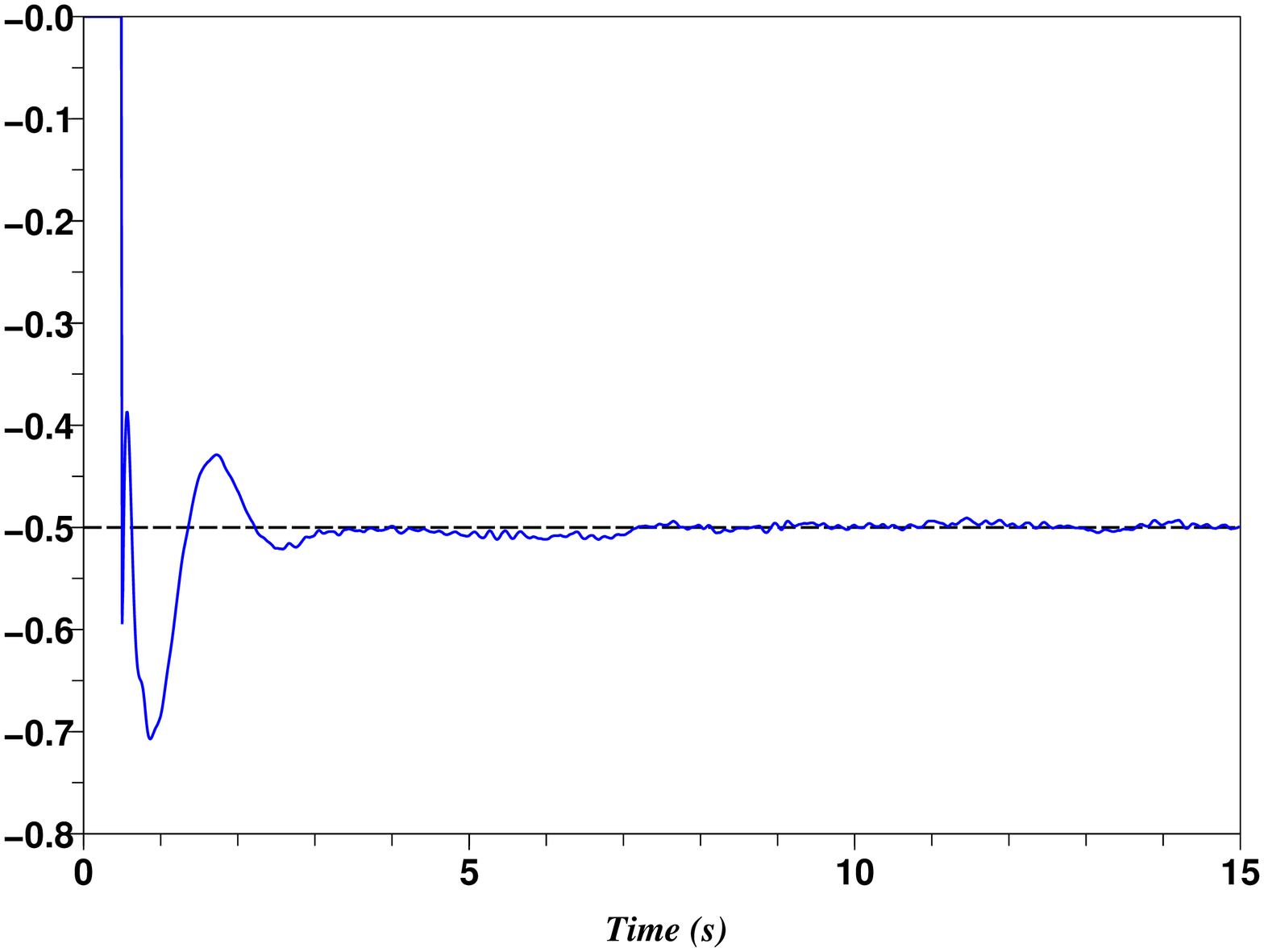}}}}}
\centering {\subfigure[\footnotesize $u$ (--); $u^\star$ (- -);
$u^\star_{\text{pert}}$(. .)]{
\rotatebox{-90}{\resizebox{!}{5cm}{%
   \includegraphics{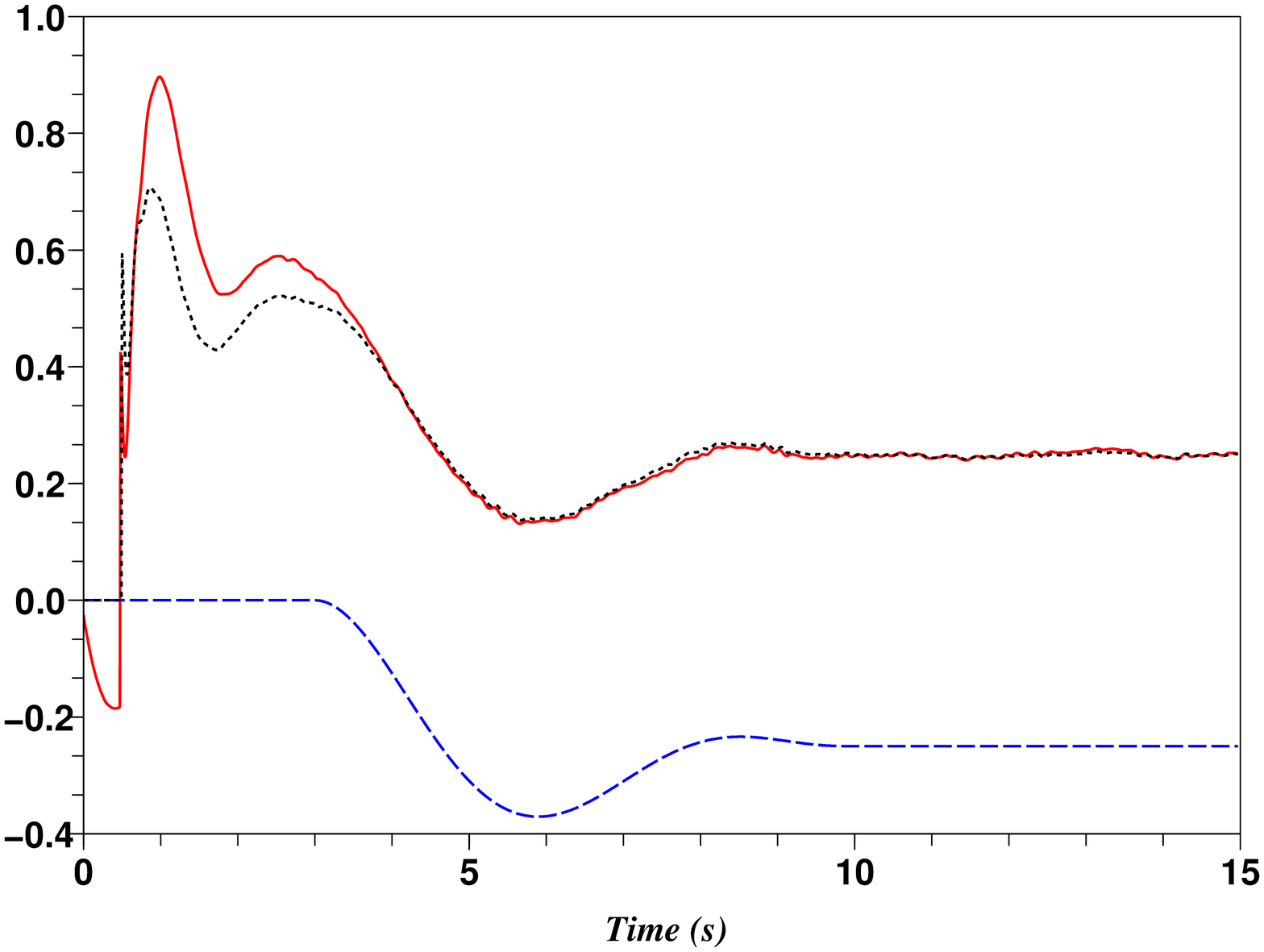}}}}}
{\subfigure[\footnotesize Output (--); reference (- -); denoised
output (. .)]{
\rotatebox{-90}{\resizebox{!}{5cm}{%
   \includegraphics{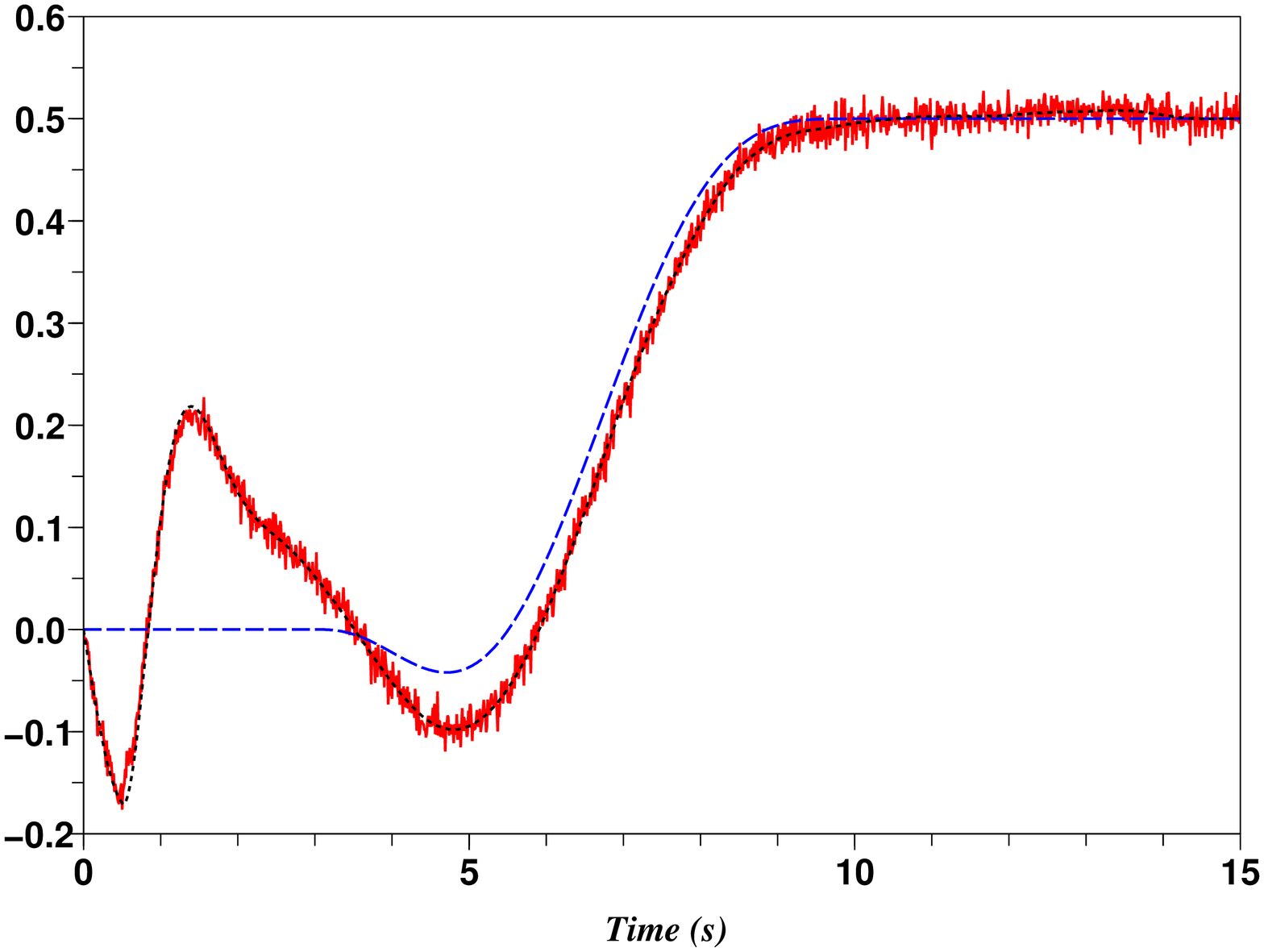}}}}}
 \caption{The non-minimum phase system where the first effect is not modeled \label{pnm1} }
\end{figure*}

\begin{figure*}[H]
\centering {\subfigure[\footnotesize $u$ (--); $u^\star$ (- -)]{
\rotatebox{-90}{\resizebox{!}{5cm}{%
   \includegraphics{pnm1su.eps}}}}}
{\subfigure[\footnotesize Output (--); reference (- -); denoised
output (. .)]{
\rotatebox{-90}{\resizebox{!}{5cm}{%
   \includegraphics{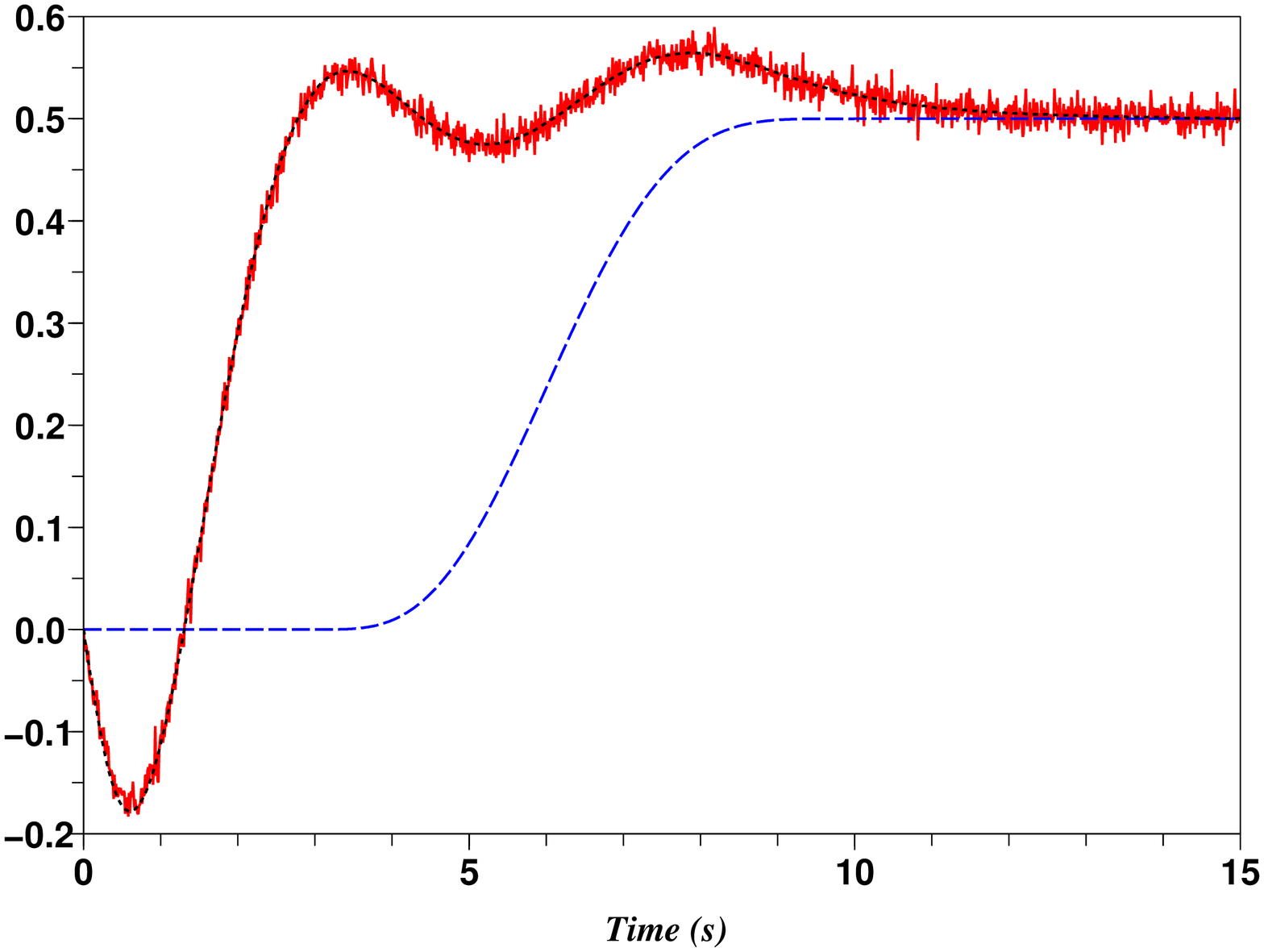}}}}}
{\subfigure[\footnotesize Perturbation (- -); estimated perturbation
(--)]{
\rotatebox{-90}{\resizebox{!}{5cm}{%
   \includegraphics{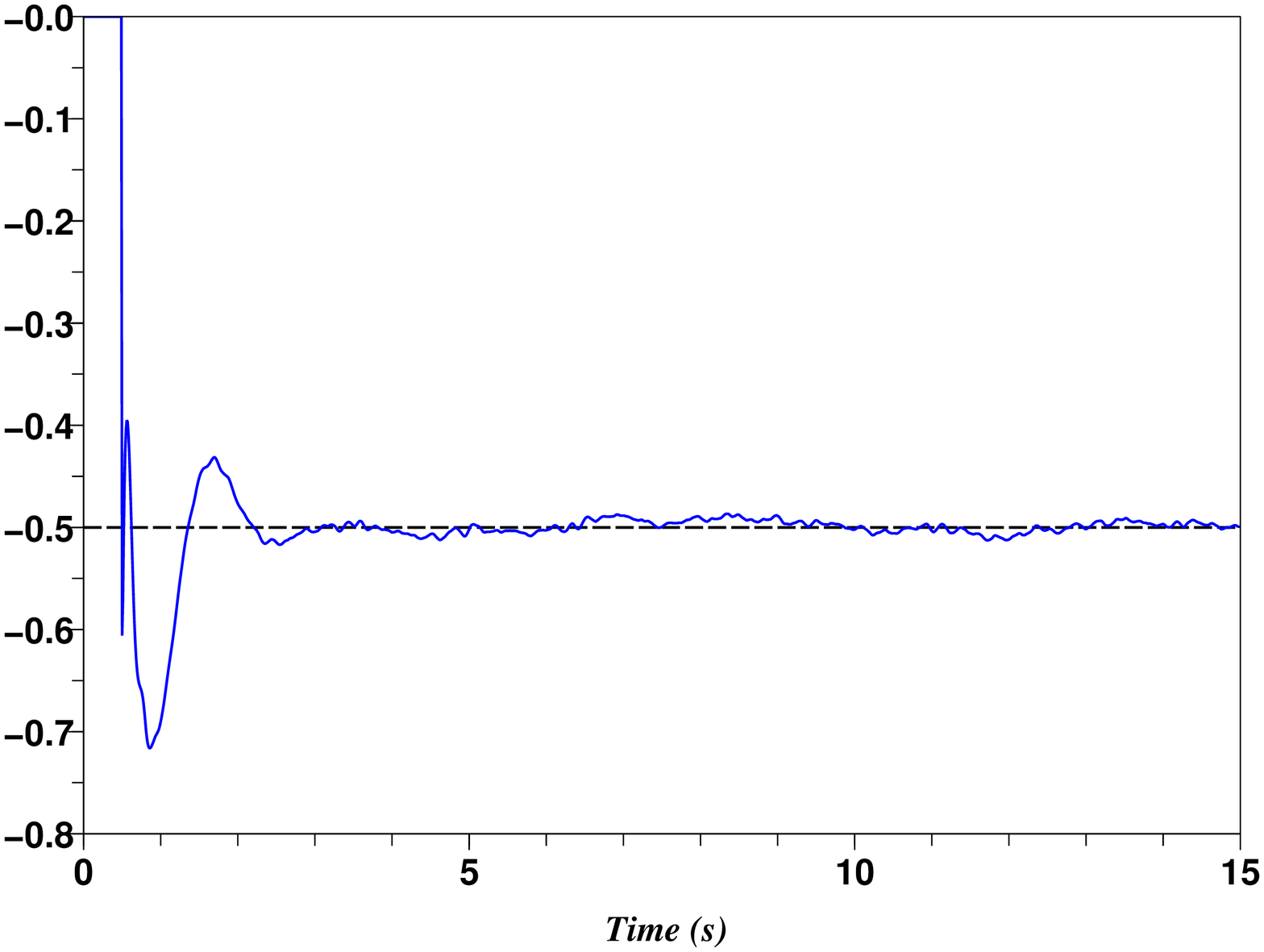}}}}}
\centering {\subfigure[\footnotesize $u$ (--); $u^\star$ (- -); $u^\star_{\text{pert}}$(. .)]{
\rotatebox{-90}{\resizebox{!}{5cm}{%
   \includegraphics{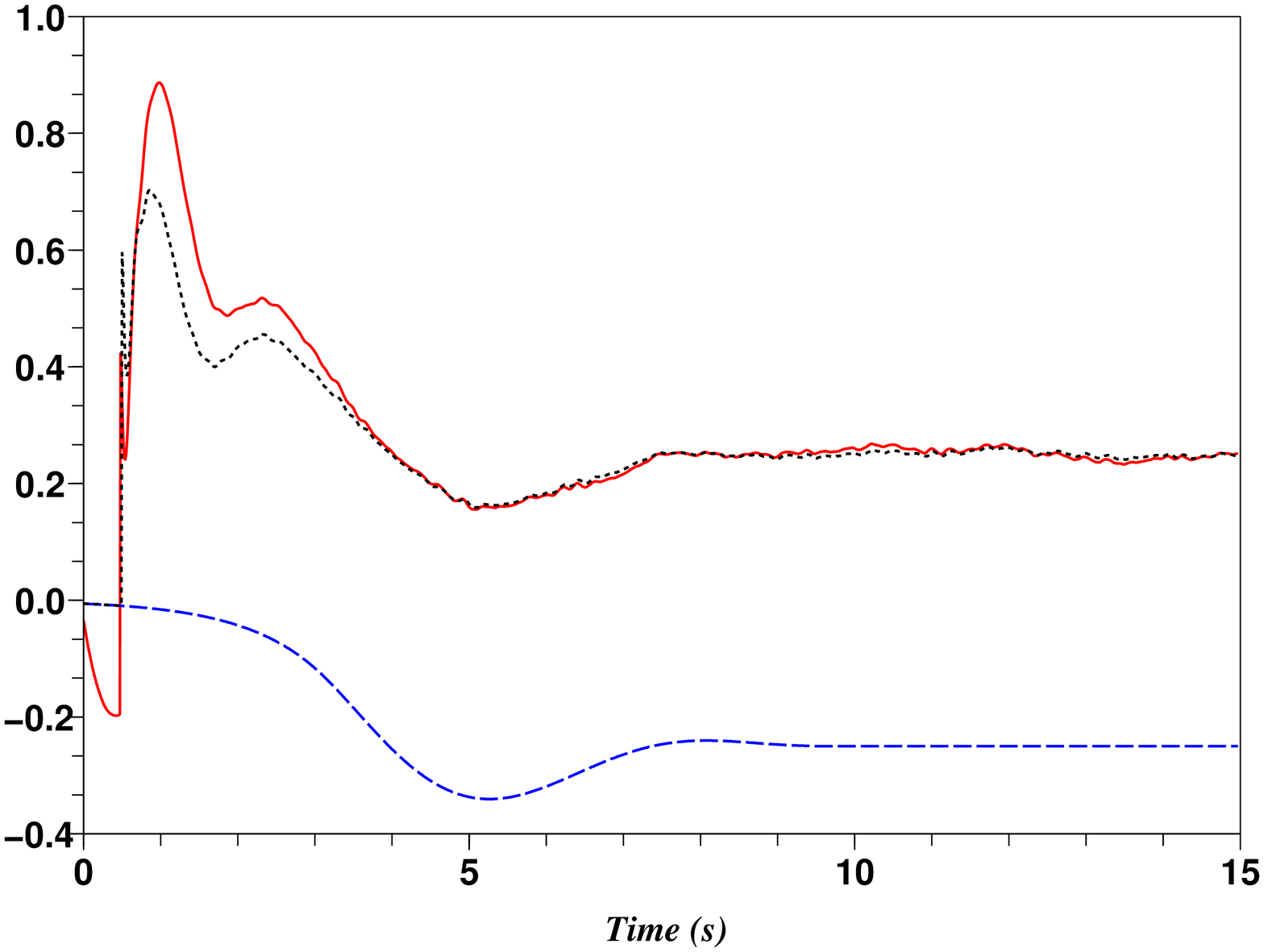}}}}}
{\subfigure[\footnotesize Output (--); reference (- -); denoised
output (. .)]{
\rotatebox{-90}{\resizebox{!}{5cm}{%
   \includegraphics{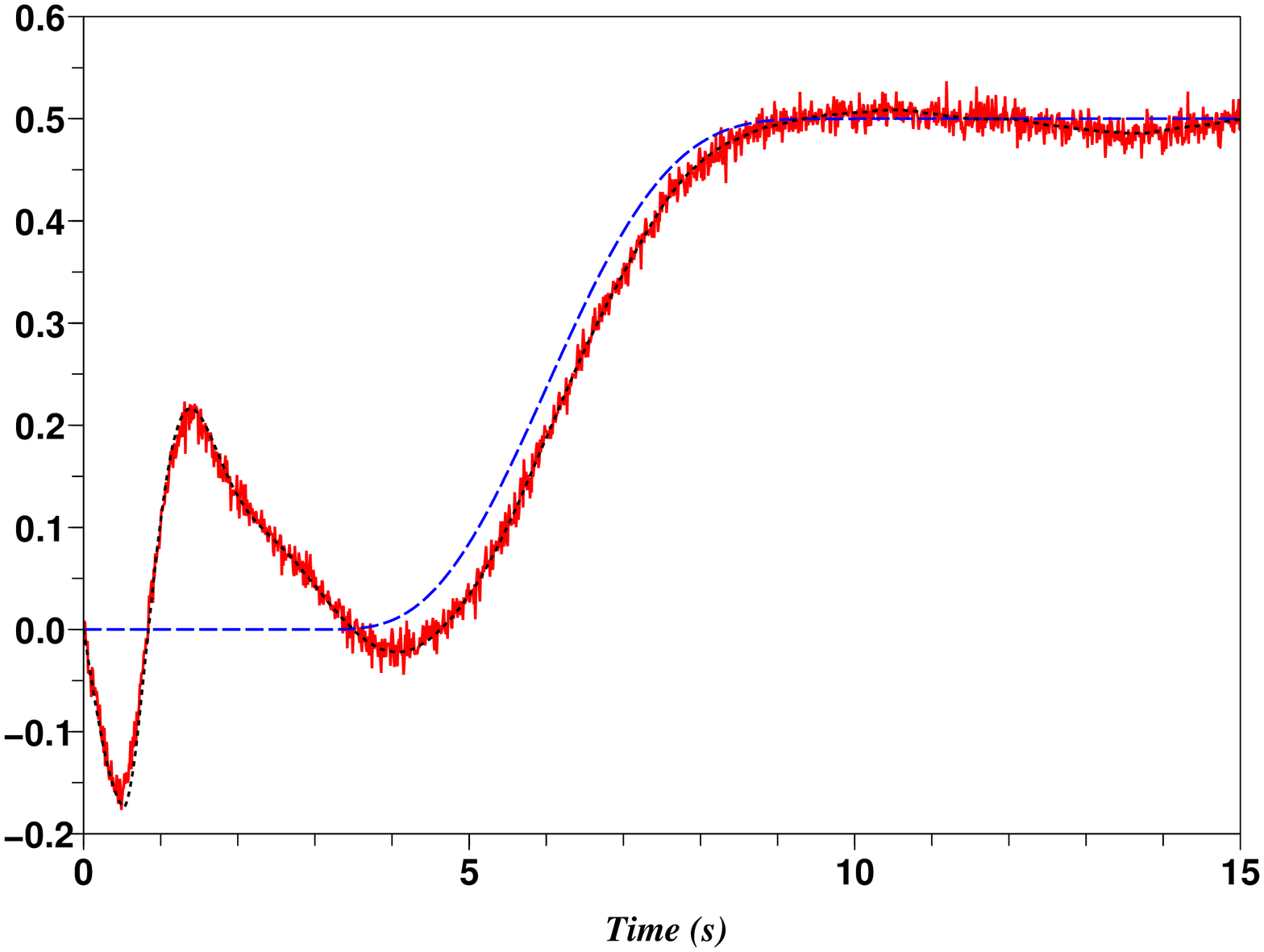}}}}}
 \caption{The non-minimum phase system where the first effect is not modeled \label{pnm2}}
\end{figure*}
\begin{figure*}[H]
\centering {\subfigure[\footnotesize $u$ (--); $u^\star$ (- -)]{
\rotatebox{-90}{\resizebox{!}{5cm}{%
   \includegraphics{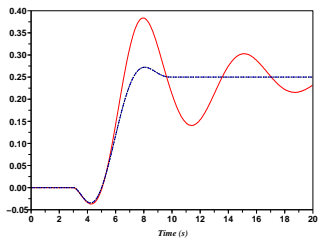}}}}}
{\subfigure[\footnotesize Output (--); reference (- -); denoised
output (. .)]{
\rotatebox{-90}{\resizebox{!}{5cm}{%
   \includegraphics{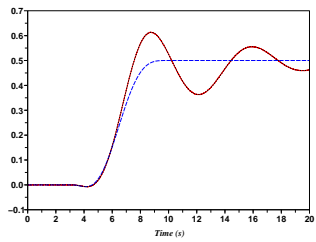}}}}}
{\subfigure[\footnotesize Perturbation (- -); estimated perturbation
(--)]{
\rotatebox{-90}{\resizebox{!}{5cm}{%
   \includegraphics{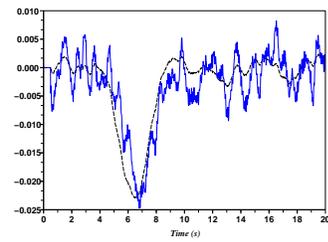}}}}}
\centering {\subfigure[\footnotesize $u$ (--); $u^\star$ (- -); $u^\star_{\text{pert}}$(. .)]{
\rotatebox{-90}{\resizebox{!}{5cm}{%
   \includegraphics{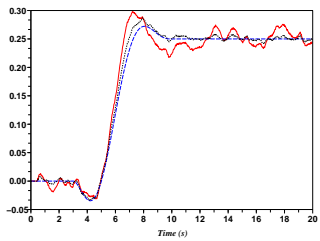}}}}}
{\subfigure[\footnotesize Output (--); reference (- -); denoised
output (. .)]{
\rotatebox{-90}{\resizebox{!}{5cm}{%
   \includegraphics{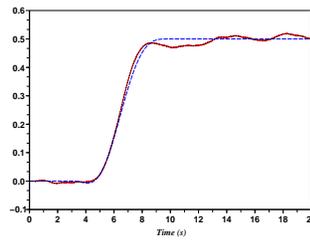}}}}}
 \caption{The non-minimum phase system where the second effect is not modeled \label{pnm3}}
\end{figure*}

\end{document}